\documentclass[11pt,a4paper]{article}
\usepackage[T1]{fontenc}
\usepackage{lmodern}
\usepackage[english]{babel}
\usepackage[margin=24mm]{geometry}
\usepackage{amsmath,amssymb,amsthm,mathtools,mathrsfs}
\usepackage{microtype,booktabs,xcolor,verbatim,graphicx}
\usepackage{enumitem,tikz,tikz-cd,pgfplots}
\usepackage{float}
\usepackage{caption,placeins}
\usetikzlibrary{arrows.meta}
\pgfplotsset{compat=1.18}
\usepackage[colorlinks=true,linkcolor=black,citecolor=black,urlcolor=blue!40!black]{hyperref}
\hypersetup{pdftitle={Graph continuity and spectral flow for singular Sturm--Liouville families with varying domains},pdfauthor={Xijun Hu; Alessandro Portaluri; Li Wu}}
\numberwithin{equation}{section}
\newtheorem{theorem}{Theorem}[section]
\newtheorem{proposition}[theorem]{Proposition}
\newtheorem{lemma}[theorem]{Lemma}
\newtheorem{corollary}[theorem]{Corollary}
\theoremstyle{definition}

\theoremstyle{remark}
\newtheorem{remark}[theorem]{Remark}
\newtheorem{example}[theorem]{Example}
\newtheoremstyle{mainresults}
  {12pt}        % Space above
  {12pt}        % Space below
  {\itshape}    % Statement in italics
  {}            % No indentation
  {\bfseries}   % Bold heading
  {.}           % Punctuation
  {\newline}    % Statement starts on a new line
  {\thmname{#1}\thmnumber{ #2}\thmnote{ --- #3}}

\theoremstyle{mainresults}
\newtheorem*{introA}{Theorem A}
\newtheorem*{introB}{Theorem B}
\newtheorem*{introC}{Theorem C}
\newtheorem*{introD}{Theorem D}
\newtheorem*{introE}{Theorem E}
\newtheorem*{introF}{Theorem F}

\theoremstyle{plain}
\newcommand{\R}{\mathbb R}
\newcommand{\C}{\mathbb C}

\newcommand{\D}{\mathcal D}
\newcommand{\B}{\mathcal B}

\newcommand{\F}{\mathcal F}
\renewcommand{\sf}{\operatorname{sf}}
\newcommand{\Mas}{\mu_{\mathrm{CLM}}}

\newcommand{\itriple}{\iota}
\newcommand{\nm}{n_-}
\newcommand{\np}{n_+}
\newcommand{\ind}{\operatorname{ind}}

\newcommand{\ran}{\operatorname{ran}}
\newcommand{\gr}{\operatorname{gr}}
\newcommand{\Lag}{\operatorname{Lag}}
\newcommand{\Span}{\operatorname{span}}

\newcommand{\Tr}{\operatorname{Tr}}

\newcommand{\ip}[2]{\langle #1,#2\rangle}

\setlist{itemsep=3pt,topsep=5pt}

\title{Graph continuity and spectral flow for singular\\
Sturm--Liouville families with varying domains}
\author{Xijun Hu\thanks{Partially supported by the National
Natural Science Foundation of China (Grant number: 12521001) and the
Taishan Scholars Climbing Program of Shandong Province
(Grant number: TSPD20240802).}
\and Alessandro Portaluri\thanks{Partially supported by Progetto di
Ricerca GNAMPA--INdAM.}
\and Li Wu\thanks{Partially supported by the National Natural Science
Foundation of China, grant 12171281.}}
\date{}

\begin{document}
\maketitle

\begin{abstract}
We develop a spectral-flow theory for gap-continuous paths of closed
symmetric Fredholm operators with varying domains on a Hilbert space. Any
continuous family of finite-dimensional graph complements forms a
symplectic boundary bundle characterizing norm-resolvent-continuous
self-adjoint extensions via continuous Lagrangian subspaces, even when the
minimal kernel dimension varies; for injective minimal operators, interval
spectral flow equals minus the Maslov index of the boundary conditions and
adjoint-kernel Cauchy data. For singular Sturm--Liouville systems with
varying domains, a quantitative Volterra criterion verifies minimal gap
continuity from homogeneous solutions, and interior cutting extends the
formula when uncut minimal operators have nontrivial kernels.

Applied to quadratic regularizations of Bessel-type inverse-square
potentials across $-1/4$ and throughout the limit-circle range, this
theory yields joint minimal gap continuity, norm-resolvent-continuous
families realizing every singular left-endpoint condition with the right
Dirichlet condition fixed, and non-convergence of two-ended regular
Dirichlet realizations for each fixed supercritical parameter. For matrix
Bessel-type operators with lower-order perturbations and noncommuting
gyroscopic Hill systems, it gives an explicit logarithmic Morse-index
count with remainder uniform over all self-adjoint boundary conditions.

In the planar gravitational $N$-body problem, these results establish
uniform Floquet Morse-index asymptotics with bounded remainder for
elliptic relative equilibria as the eccentricity tends to one. For the
collinear Euler three-body family, they unconditionally show that the
first antiperiodic degeneracy curve tends to zero while all remaining
fixed-label periodic and antiperiodic branches tend to $1/8$, with no
curves accumulating in $(0,1/8)$. In the logarithmic eccentricity scale,
they resolve the bulk locations and spacings to subleading order and show
that near $1/8$ neighboring triples of branches approach common limiting
rays with distinct quadratic displacements.
\end{abstract}

\setcounter{tocdepth}{1}
\begingroup
\small
\tableofcontents
\endgroup

\section{Introduction and main results}

We develop a boundary-bundle spectral-flow theory for gap-continuous
families of closed symmetric Fredholm operators with varying domains on a
Hilbert space, together with quantitative criteria for singular
Sturm--Liouville operators and applications to planar $N$-body elliptic
relative equilibria. When an endpoint becomes singular, continuous
symplectic coordinates in a graph complement compare boundary conditions
across the limit: a fixed physical boundary condition can rotate while the
Cauchy-data space of homogeneous solutions converges, turning Maslov
rotation into logarithmic Morse-index counts with uniformly bounded remainders and sharp Euler
degeneracy-curve asymptotics down to the Hardy threshold.

Spectral flow was introduced by Atiyah, Patodi and Singer \cite{APS76} and
formulated operator-theoretically by Phillips \cite{Phi96}; gap-topology
spectral flow for unbounded self-adjoint Fredholm operators is developed
in \cite{BLP05}, and its relation to the Maslov index and symplectic paths
in \cite{CLM94,RS95,BF98,ZL99,LZ00,Fur04}. Boundary-space spectral-flow
formulas appear for fixed maximal domains \cite{BZ05}, varying symplectic
Banach domain bundles \cite{BZ18}, second-order elliptic boundary value
problems \cite{LS18}, and half-line heteroclinic and homoclinic orbits
\cite{HP17}; graph-based $K$-theoretic family indices for Fredholm
extensions of semi-Fredholm operators are studied in \cite{Pro24}.
Related boundary-space results on interlacing and perturbation theory
appear in \cite{BCLS26,LS24}, renormalized oscillation theory for singular
Hamiltonian systems in \cite{HS22}, and the singular Sturm--Liouville
boundary-quotient theory for fixed minimal domains in \cite{HPW25} (see
also \cite{Zet05}).

Building on this literature and extending \cite{HPW25} to varying minimal
domains, our results advance both the abstract and differential-operator
theory. On the abstract side (Theorem~A), we start solely from a
gap-continuous path of closed symmetric Fredholm operators on a fixed
Hilbert space $H$, without an ambient Banach domain bundle or an auxiliary
invertible reference extension: arbitrary continuous graph complements in
$H\oplus H$ characterize all norm-resolvent-continuous self-adjoint
extensions via continuous Lagrangian paths even when the minimal kernel
dimension varies, and under minimal injectivity scalar shifts and
symplectic shears give a direct interval spectral-flow formula including
endpoint crossings. On the analytical side (Theorem~B), whereas abstract
theories take minimal graph continuity as an input, our Volterra
$L^2$-frame criterion establishes gap continuity of singular
Sturm--Liouville minimal realizations directly from homogeneous solutions
over arbitrary parameter spaces, and interior cutting removes the
uncut-injectivity restriction.

In the planar gravitational $N$-body problem, a central configuration is a
mass arrangement in which the gravitational force on each body is
proportional to its mass times its position relative to the center of
mass. Each planar central configuration generates a family of elliptic
relative equilibria, parameterized by the eccentricity $e\in[0,1)$, in
which the configuration evolves by a common rotation and dilation. The
Meyer--Schmidt symplectic reduction \cite{MS05} reduces the linearized
variational equation to a matrix Sturm--Liouville system with a gyroscopic
term whose coefficients develop an inverse-square singularity as
$e\uparrow1$.

For homographic solutions of the planar three-body problem, Mart\'inez,
Sam\`a and Sim\'o analyzed near-circular resonances \cite{MSS06a} and
near-collision Floquet asymptotics under two nonvanishing connection
conditions \cite[Theorem~1]{MSS06b} (verified for each $\delta>1/8$ and
small matching cutoff in Lemma~\ref{lem:euler-mss-matching}). For general
central configurations, Hu and Ou proved a logarithmic Morse-index lower
bound in the spiral regime ($\lambda_{\min}(R)<-1/8$)
\cite[Theorem~1.1]{HO16}, and in the strict non-spiral regime obtained
near-collision index formulas under nondegeneracy and collision
nondegeneracy \cite[Theorem~1.2]{HO16}. Hu, Ou and Yu obtained a
physical-time logarithmic ratio limit at a collision or parabolic end of a
doubly asymptotic orbit when the limiting configuration at that end
satisfies the spiral condition \cite[Theorem~1.1(a)]{HOY21} (see
\eqref{eq:euler-physical-end}). For collinear Euler configurations, Long
and Zhou \cite{ZL17} determined the global curve structure, and strict
separation and sector assignments within each antiperiodic pair for
$0<e<1$ were proved in \cite[Theorem~1.1 and Remark~1.2]{HWZ26}. As
$e\uparrow1$, working in boundary-bundle coordinates, we sharpen the
$N$-body count to a uniform asymptotic equality with bounded remainder
over all Floquet multipliers, determine all limiting Euler endpoints
unconditionally, and resolve the bulk and near-threshold curve spacings.

\paragraph{Boundary bundles and spectral flow.}
For a densely defined closed operator $A$ on a real separable Hilbert space
$H$, write $\D(A)$ for its domain, $A^*$ for its adjoint, and
$\gr A=\{(u,Au):u\in\D(A)\}$. If $P_M$ denotes the orthogonal
projection onto a closed subspace $M$, gap continuity means norm continuity
of $P_{\gr A}$. For a symmetric Fredholm family $A_s$, choose
continuous finite-dimensional graph complements $E_s$ and set
\begin{equation}\label{eq:intro-canonical-E}
 \gr A_s^*=\gr A_s\oplus E_s,
 \qquad C_s=P_{E_s}(\ker A_s^*\times\{0\}),
\end{equation}
where $P_{E_s}:\gr A_s^*\to E_s$ is the projection along $\gr A_s$.
Continuity of subspaces is understood in the gap topology.
The boundary form $\widehat\omega_s$ is the restriction to $E_s$ of
$\widetilde\omega((u,v),(x,y))=\langle u,y\rangle-\langle v,x\rangle$.
Under the natural identification of $E_s$ with $\D(A_s^*)/\D(A_s)$,
$C_s$ represents the Cauchy data of $\ker A_s^*$.
We use $\sf$ for spectral flow with upward
crossings positive, and $\Mas$ for the Cappell--Lee--Miller Maslov index
with the endpoint conventions of Section~\ref{sec:abstract}.

We denote by $\mathcal{CF}^{\mathrm{sym}}(H)$ the space of
densely defined, closed, symmetric Fredholm operators on $H$,
equipped with the gap topology, and by
$\mathcal{CF}^{\mathrm{sa}}(H)$ its self-adjoint subspace. 

\par\begin{samepage}
\begin{introA}[Boundary bundles and varying-domain spectral flow]
Let $A:[0,1]\to\mathcal{CF}^{\mathrm{sym}}(H)$,
$s\mapsto A_s$, be a gap-continuous path, and let $E_s$ be any
continuous family of finite-dimensional subspaces satisfying
$\gr A_s^*=\gr A_s\oplus E_s$. With $C_s$ as in
\eqref{eq:intro-canonical-E}, the following hold.
\begin{enumerate}[label=(\roman*)]
\item The spaces $(E_s,\widehat\omega_s)$ form a continuous
symplectic bundle of constant rank $-2\ind A_s$.

\item Lagrangian subspaces $N_s\subset E_s$ determine self-adjoint
extensions $B_s$ through
$\gr B_s=\gr A_s\oplus N_s$.
The family $B_s$ is norm-resolvent continuous if and only if
$N_s$ varies continuously.

\item If $\ker A_s=\{0\}$ for every $s$, then $\ker A_s^*$ varies
continuously as a subspace of $H$, and the Cauchy-data spaces
$C_s$ form a continuous Lagrangian family in $E_s$. For every continuous family of
Lagrangian subspaces $N_s\subset E_s$, the realizations in (ii) satisfy
\begin{equation}\label{eq:intro-sf-formula}
\sf(B_s;s\in[0,1])
=-\Mas(N_s,C_s,\widehat\omega_s;s\in[0,1]).
\end{equation}
\end{enumerate}
\end{introA}
\end{samepage}
\begin{remark}
For each $s$, the boundary construction gives a Lagrangian Cauchy-data
space even when $\ker A_s$ is nontrivial. Injectivity ensures that
these spaces vary continuously and that the boundary quotient detects
the entire kernel of $B_s$.
Under this hypothesis,
\[
\ker B_s\simeq N_s\cap C_s.
\]
In general the intersection identifies with $\ker B_s/\ker A_s$.
For example, $A_s=B_s=(s-\tfrac12)I$ on $\mathbb R$, $0\le s\le1$,
has spectral flow $+1$ although the boundary quotient is zero throughout.
\end{remark}

Theorem~A is proved in Section~\ref{sec:abstract}
(Theorems~\ref{thm:canonical} and~\ref{thm:sf});
Example~\ref{ex:robin-boundary-family} and
Figure~\ref{fig:robin-boundary-crossing} work out the bundle objects
explicitly for a moving Robin boundary condition.

\paragraph{General Sturm--Liouville families.}
To apply Theorem~A to singular differential operators, one must establish
gap continuity and, for the index formula, obtain an injective minimal
family. Theorem~B provides an adjoint-kernel criterion for gap continuity
and injectivity, together with a separate interior-cut construction for
general gap-continuous families whose uncut minimal operators may fail to
be injective.
The minimal operator of a differential expression is the closure
of its restriction to compactly supported smooth test functions.

\begin{introB}[Sturm--Liouville families: graph continuity and interior cutting]
Let $X$ be a topological space and consider, in
$H=L^2((0,1);\mathbb R^n)$,
\[
 \tau_vy=-(P_vy'+Q_vy)'+Q_v^Ty'+R_vy,\qquad v\in X,
\]
where $P_v,Q_v$ are locally $C^1$, $R_v$ is locally continuous,
$P_v=P_v^T$ is invertible and $R_v=R_v^T$.
Both endpoints may be singular. Let $T_v$ be the minimal realizations.
\begin{enumerate}[label=(\roman*)]
\item For arbitrary $X$, suppose $\dim\ker T_v^*=2n$ and these kernels vary continuously
in the gap topology of subspaces of $H$. Locally in $v$, suppose also that
$P_v^{-1},Q_v,R_v$ vary continuously in $C^0(I)$ on a fixed
compact interval $I\Subset(0,1)$ with nonempty interior.
Then the minimal realizations $T_v$ are gap-continuous,
injective Fredholm operators of index $-2n$, and
their maximal graph domains embed compactly into $H$.

\item Let $X=[0,1]$ and suppose $v\mapsto T_v$ is gap-continuous and
Fredholm. Fix $c\in(0,1)$ and a compact interval
$I\Subset(0,1)$ with $c\in\operatorname{int}I$.
Locally in $v$, suppose $P_v,Q_v$ vary continuously in $C^1(I)$,
$R_v$ varies continuously in $C^0(I)$, and $P_v^{-1}$ is uniformly
bounded on $I$. Let
$\widehat T_v=T_v^-\oplus T_v^+$ be the direct sum of the minimal
realizations on $(0,c)$ and $(c,1)$.
Then $\widehat T_v$ is gap-continuous, injective and Fredholm, with
\[
 \ind\widehat T_v=\ind T_v-2n.
\]

\item Under the hypotheses of (ii), let $\widetilde T_v$ be a gap-continuous family
of self-adjoint extensions of $T_v$. In a continuous graph complement
$\widehat E_v$ for the cut family, let $\widehat N_v$ represent the
domain of $\widetilde T_v$, including the matching condition at $c$, and let
$\widehat C_v$ be the Cauchy-data space of $\widehat T_v$.
Then
\[
 \sf(\widetilde T_v;v\in[0,1])
 =-\Mas(\widehat N_v,\widehat C_v,\widehat\omega_v;v\in[0,1]).
\]
\end{enumerate}
\end{introB}
Part~(i) is proved in Section~\ref{sec:volterra}
(Theorem~\ref{lem:volterra}), and parts~(ii)--(iii) in
Section~\ref{sec:sl} (Theorem~\ref{thm:cut}). The two parts are
independent: part~(i) allows an arbitrary parameter space $X$ (used for
$(s,r)$ in Theorem~C) and assumes all $2n$ homogeneous solutions lie in
$L^2$, which automatically yields a compact Volterra right inverse and
hence Fredholmness without a separate hypothesis. This applies directly to
the scalar Bessel models of
Sections~\ref{sec:bessel}--\ref{sec:renormalization} without cutting; for
$N$-body counts min--max reduces to scalar comparisons, while the coupled
Euler family combines Bessel boundary coordinates on its limit-circle
component with weighted energy estimates on its limit-point component.
Under~(ii), the uncut minimal operator may have a nontrivial kernel;
Appendix~\ref{sec:coupled} illustrates how cut coordinates handle
transmission and coupled outer boundary conditions.

\begin{remark}[Relation to half-line index formulas]
Under asymptotic hyperbolicity assumptions, Hu and Portaluri \cite{HP17}
express spectral flow for heteroclinic, homoclinic and halfclinic problems
(including varying half-line Lagrangian boundary conditions) via stable
and unstable Lagrangian spaces. Here the minimal and maximal graph domains
themselves may vary as an endpoint becomes singular. Under the
separated-boundary hypotheses of Theorem~\ref{thm:trace-index} (matching
across $c$ with continuous half-interval Lagrangian Cauchy traces
$W_v^\pm(c)\subset\R^{2n}$), the formula reduces to
\[
 \sf(\widetilde T_v;v\in[0,1])
 =-\Mas\bigl(W_v^+(c),W_v^-(c);\R^{2n},\omega_0\bigr),
\]
directly paralleling \cite{HP17} without requiring asymptotic
hyperbolicity.
\end{remark}

\paragraph{Bessel regularization and singular boundary selection.}
The Bessel operator $-d^2/dt^2+q/t^2$ on $(0,1)$ is the local model for
our inverse-square singularities \cite{Zet05,HPW25}. Its spectral theory
is governed by two thresholds (Remark~\ref{rem:thresholds}): Weyl's
limit-circle/limit-point value $q=3/4$ (for fixed $q\ge3/4$, the
regularized minimal operators do not converge in the gap topology to the
singular minimal operator as $s\downarrow0$,
Proposition~\ref{prop:threshold}), and within $q<3/4$ the Hardy threshold
$q=-1/4$, separating semibounded realizations with finite Morse index
($q\ge-1/4$) from oscillatory realizations with infinite Morse index
($q<-1/4$). Explicit frames for its quadratic regularization establish gap
continuity across $q=-1/4$ throughout $q<3/4$ and prescribe every limiting
singular left-endpoint condition with $y(1)=0$ fixed.

\begin{introC}[Bessel regularization and selection of singular boundary conditions]
In $L^2(0,1)$, consider
\[
 \ell_{s,r}=-\frac{d^2}{dt^2}
 +\frac{r|r|-1/4}{t^2+s^2},\qquad 0\le s\le1,\quad r<1,
\]
and let $G_{s,r}$ be its minimal realization.
\begin{enumerate}[label=(\roman*)]
\item The map $(s,r)\mapsto G_{s,r}$ is jointly gap-continuous on
$[0,1]\times(-\infty,1)$; every $G_{s,r}$ is symmetric Fredholm,
and its maximal graph domain embeds compactly in $L^2(0,1)$.

\item There exists a jointly norm-resolvent continuous family of
self-adjoint realizations $A_{s,r,\xi}$ of $\ell_{s,r}$,
indexed by $\xi\in\mathbb{RP}^1$, all satisfying $y(1)=0$. At $s=0$, this family exhausts the separated self-adjoint
realizations with that right-endpoint condition.

\item For every fixed $r<0$, the realizations of $\ell_{s,r}$ with
Dirichlet conditions at both endpoints have no norm-resolvent
limit as $s\downarrow0$.
\end{enumerate}
\end{introC}
These statements are established in Sections~\ref{sec:bessel}
and~\ref{sec:renormalization} (Theorems~\ref{thm:bessel2}
and~\ref{thm:boundary-selection}, and Proposition~\ref{prop:dirichlet-cycle}),
where we also show that continuous boundary selection preserves spectral
flow along paths with invertible singular endpoint operators;
Example~\ref{ex:quadratic-boundary-selection} illustrates both selection
and Dirichlet cycling explicitly.

\paragraph{Logarithmic counts for planar elliptic relative equilibria.}
For each supercritical scalar component the physical Dirichlet line rotates
at an explicit logarithmic rate, and bounding a matrix perturbation
between scalar potentials leaves this rotation unchanged in both
comparison problems. Because the corresponding Cauchy-data lines of the
homogeneous equations admit limits, Theorem~A converts the rotation
directly into a negative-index count; the critical comparison bounds the
negative indices of all subcritical and critical components, while
changing endpoint conditions alters the index by at most a bounded amount. We state the resulting count for planar elliptic
relative equilibria first.

Write $n_-(T)$ for the negative Morse index, with eigenvalues counted
with multiplicity, and set $x_+=\max\{x,0\}$ and
$L(e)=\log(1/(1-e))$.

For a planar $N$-body central configuration, let $d=2N-4$,
$J=\left(\begin{smallmatrix}0&-I_{d/2}\\I_{d/2}&0\end{smallmatrix}\right)$,
and let $R=R^T$ be its essential $d\times d$ matrix.
For eccentricity $e\in[0,1)$ and $\omega\in\mathbb C$, $|\omega|=1$, set
\begin{equation}\label{eq:intro-ere-operator}
 A_{R,e,\omega}=-\frac{d^2}{d\theta^2}-2J\frac{d}{d\theta}
                +\frac{R}{1+e\cos\theta},\qquad -\pi<\theta<\pi,
\end{equation}
with domain $\{y\in H^2((-\pi,\pi);\mathbb C^d):
y(\pi)=\omega y(-\pi),\ y'(\pi)=\omega y'(-\pi)\}$, and define
\[
 C(R)=\frac1\pi\sum_{j=1}^{d}
       \sqrt{\left(-2\lambda_j(R)-\frac14\right)_+}.
\]

\begin{introD}[Floquet counts and spectral flow for elliptic relative equilibria]
For fixed $R$, uniformly in $|\omega|=1$,
\[
 n_-(A_{R,e,\omega})=C(R)L(e)+O(1),\qquad e\uparrow1.
\]
For fixed $e_0\in[0,1)$ and $\omega$, the spectral flow along
$\varepsilon\in[e_0,e]$, $e_0\le e<1$, is
$-C(R)L(e)+O(1)$ as $e\uparrow1$.
\end{introD}
This follows from Theorem~\ref{thm:ere-uniform-floquet} and the
Morse-index difference formula in the proof of Theorem~\ref{thm:global-sf}.
Complex Floquet conditions use complex multiplicities; at $\omega=\pm1$
these agree with the real Morse indices.

\begin{remark}[Generality of the counting formula]\label{rem:intro-count-generality}
Theorem~\ref{thm:gyroscopic-master} proves the same Morse-index estimate for
\[
 \mathcal G_e=-\frac{d^2}{d\theta^2}-2K\frac{d}{d\theta}
             +\frac{L}{1+e\cos\theta}+B_e(\theta),\qquad -\pi<\theta<\pi,
\]
with coefficient $\pi^{-1}\sum_j\sqrt{(-2\lambda_j(L)-1/4)_+}$,
where $K^T=-K$, $L^T=L$ (not necessarily commuting), and $B_e$ is
measurable, symmetric and uniformly bounded.
For fixed $K,L$, the $O(1)$ remainder is uniform over all
self-adjoint endpoint conditions, including parameter-dependent ones,
and over perturbations with a common bound. Boundary data use
$y'+Ky$. The Floquet case above follows by taking $K=J$, $L=R$
and $B_e=0$; translating by $t=\theta+\pi$ converts $(-\pi,\pi)$ with
$1+e\cos\theta$ to the $(0,2\pi)$ and $1-e\cos t$ convention of
Theorem~\ref{thm:gyroscopic-master}. Critical eigenvalues $\lambda_j(L)=-1/8$ contribute only to $O(1)$,
and the counting formula requires no limit-circle restriction on the
remaining eigenvalues.
\end{remark}
In particular, for the same elliptic relative equilibria, Theorem~D
refines the single-eigenvalue logarithmic lower bound of Hu and Ou
\cite[Theorem~1.1]{HO16} to an asymptotic equality accounting for every
supercritical eigenvalue with a bounded remainder uniform in the Floquet
multiplier (see Section~\ref{subsec:ere-brake} for a direct formula
comparison).

\paragraph{Euler degeneracy curves.}
For the collinear three-body Euler family, the essential matrix is
$R_\delta=\operatorname{diag}(-\delta,2\delta+3)$, where the physical
mass parameter in \eqref{eq:euler-mass-parameter} lies in $[0,7]$; we
consider the extended operator family for all $\delta>0$. Following
Long--Zhou \cite{ZL17}, let $g_j(e)$ and $h_j^\pm(e)$ denote the
periodic ($\omega=1$) and antiperiodic ($\omega=-1$) degeneracy curves
in the $(\delta,e)$-plane (where the monodromy matrix acquires the
Floquet multiplier $+1$ or $-1$, equivalently where the corresponding
boundary value problem has a nontrivial kernel), ordered as in
Proposition~\ref{prop:euler-interlacing}.

At collision ($e=1$), both Friedrichs realizations have Morse index one
for small positive $\delta$, and Lemma~\ref{lem:euler-friedrichs-gap}
excludes zero from their spectra for every $0<\delta\le1/8$, so this
index persists throughout the interval. For each such fixed $\delta$,
requiring the weak component to vanish at the collision end yields a
restricted family that is uniformly bounded below and converges in norm
resolvent to the corresponding Friedrichs realization as $e\uparrow1$.
Freeing this endpoint value in each sector with a terminal Neumann
condition increases the sector index from one to two while preserving
trivial kernel, so the periodic and antiperiodic problems both have
Morse index three and trivial kernel for $e$ sufficiently close to one
(Appendix~\ref{app:euler-critical-index}). Counting crossings with multiplicity relative to their
values at small positive mass accounts for exactly one antiperiodic
crossing and no periodic crossing on $(0,1/8]$, with that antiperiodic
crossing lying in $(0,a)$ for every fixed $a\in(0,1/8)$ as $e\uparrow1$.
Together with the divergence of both indices at every fixed $\delta>1/8$,
this yields
\[
 h_1^-(e)\longrightarrow0,\qquad
 g_j(e),\,h_j^+(e)\longrightarrow\tfrac18\quad(j\ge1),\qquad
 h_j^-(e)\longrightarrow\tfrac18\quad(j\ge2)
\]
for fixed labels (Corollary~\ref{cor:euler-first-pair}).
Mart\'inez, Sam\`a and Sim\'o identified the first tongue's two
endpoints numerically \cite[Section~7]{MSS06b}. Hu and Ou obtained the
endpoint limits for all fixed-label periodic and antiperiodic curves
under their Numerical Result~A \cite[Theorem~5.9]{HO16}.

Above the threshold, Theorem~D and the jump multiplicities give the
bulk distribution. For fixed $1/8<a<b<\infty$, as $e\uparrow1$,
\[
 \#\{j:a<g_j(e)<b\}
 =\frac{L(e)}{2\pi}
   \left(\sqrt{2b-\frac14}-\sqrt{2a-\frac14}\right)+O(1).
\]
The count of antiperiodic labels $(j,\pm)$ in $(a,b)$ has twice this
leading term.

The resulting $O(L^{-1})$ location error is of the same order as the
bulk gaps (Theorem~\ref{thm:euler-density}). To resolve the gaps we use
time-reversal symmetry on $[0,2\pi]$ (Appendix~\ref{app:reversible}) to
split the periodic ($NN,DD$) and antiperiodic ($ND,DN$) problems on
$[0,\pi]$ according to Neumann ($N$) or Dirichlet ($D$) conditions for the
weak component at $0$ and $\pi$. In the continuous Bessel boundary
coordinates across $\delta=1/8$, the subleading bulk locations and gaps
are governed by two limiting phase functions of
$\sigma=\sqrt{2\delta-1/4}$: a limiting $NN$ slow phase $B$ and an inner
Bessel phase shift $a$ between the weak $N$ and $D$ conditions.

\begin{introE}[Euler bulk locations and gaps]
There exist a continuous function $B:[0,\infty)\to\R$ with $B(0)=0$ and a
smooth function $a:[0,\infty)\to[0,\pi/2)$ with $a(0)=0$ such that, as
$e\uparrow1$:
\begin{enumerate}[label=(\roman*)]
\item Put $L=L(e)$ and $k_j=2\pi j$.
Uniformly for integers $j$ with $j/L$ in a fixed compact
subinterval of $(0,\infty)$, the periodic bulk locations satisfy
\[
 g_j(e)=\frac18+\frac12\left(\frac{k_j}{L}\right)^2-\frac{2k_j}{L^2}B\!\left(\frac{k_j}{L}\right)+o(L^{-1}).
\]
\item Write $s_j=\sqrt{2g_j(e)-1/4}$, $a_j=a(s_j)$ and $L=L(e)$.
Uniformly over the bulk labels in (i),
\begin{equation}\label{eq:intro-euler-splitting}
 g_j-h_j^+=\frac{2s_ja_j}{L}+o(L^{-1}),\qquad
 h_j^+-h_j^-=\frac{2s_j(\pi-2a_j)}{L}+o(L^{-1}).
\end{equation}
\end{enumerate}
\end{introE}
Here $a(\sigma)$ and $B(\sigma)$ are the inner Bessel phase shift and
limiting $NN$ slow phase constructed in Lemma~\ref{lem:bessel-inner-phase},
Remark~\ref{rem:euler-critical-chart}, and
Theorem~\ref{thm:euler-uniform-localization}; part~(i) follows from
$s_j=k_j/L+O(L^{-1})$ (Theorem~\ref{thm:euler-density}) and uniform
continuity of $B$, while~(ii) follows from
Corollary~\ref{cor:euler-gap-regimes}.

\paragraph{The critical end of the Euler curves.}
Theorem~E applies only when $j/L$ stays in a compact subinterval of
$(0,\infty)$, where $o(L^{-1})$ absorbs $O(L^{-2})$. For fixed $j\ge2$,
one has $j/L\to0$ and $\sigma\to0$, so Theorem~F(i) does not follow from
\eqref{eq:intro-euler-splitting} and requires phase control relative to
$\sigma$ near $\delta=1/8$.
\par
\begin{figure}[!tbp]
 \centering
 \includegraphics[width=\textwidth]{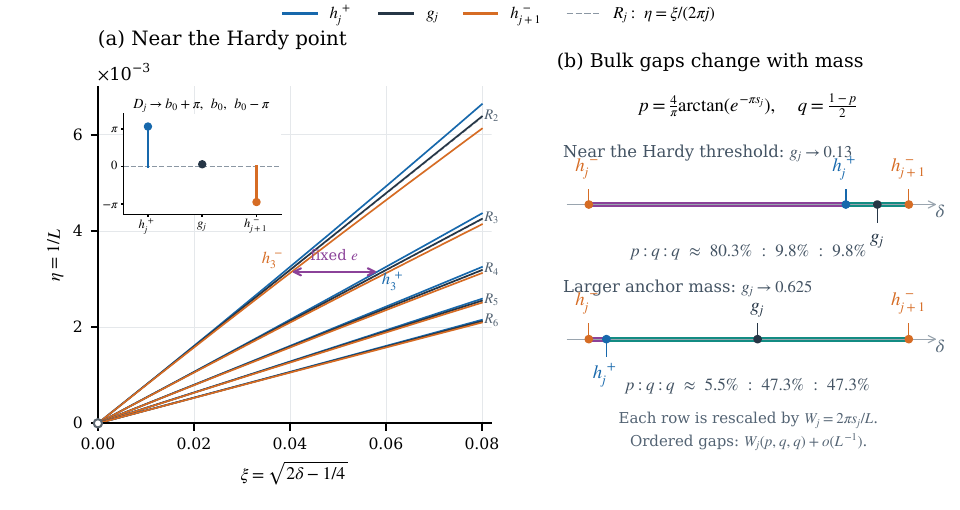}
 \caption{Euler curve geometry in the critical and bulk regimes.
 (a) Quadratic models $\eta=\xi/k_j+2c\xi^2/k_j^2$ for $j=2,\ldots,6$
 (Theorem~F), with dashed reference rays $R_j$; the inset shows the
 limits of the normalized ray deviations $D_j$ in
 \eqref{eq:euler-scaled-ray-deviation}. The horizontal arrow compares
 $h_3^-$ and $h_3^+$ at fixed $e$, which approach adjacent rays; the
 axes have unequal scales.
 (b) Leading bulk gap fractions (Theorem~E) at reference values $g_j\to0.13$ and
 $g_j\to0.625$, within compact positive $j/L$ ranges. Each row is
 normalized by $W_j=2\pi s_j/L$, where $s_j=\sqrt{2g_j-1/4}$; $p$
 denotes the leading normalized width of $h_j^+-h_j^-$, and $q$ that of
 each adjacent gap around $g_j$. Both panels show the asymptotic models.}\label{fig:euler-quantile}
\end{figure}

\begin{samepage}
\begin{introF}[Critical Euler gaps and curve geometry]
Put $L=L(e)$.
\begin{enumerate}[label=(\roman*)]
\item For each fixed $j\ge2$, as $e\uparrow1$,
\[
 h_j^+-h_j^-\sim2\pi^2(2j-1)/L^2.
\]

\item The phase functions $a$ and $B$ of Theorem~E admit finite right
derivatives $a_1=a'(0^+)>0$ and $b_0=B'(0^+)$. Set
$\xi=\sqrt{2\delta-1/4}$, $\eta=L^{-1}$, and
\[
 c=b_0+a_1,\ b_0,\ b_0-a_1
\]
for the branches $h_j^+,g_j,h_{j+1}^-$, respectively.
Along each of these branches,
\[
 \eta=\frac\xi{k_j}+\frac{2c\xi^2}{k_j^2}
       +o\!\left(\frac{\xi^2}{k_j^2}\right),
 \qquad (\xi,\eta)\to(0,0),
\]
uniformly over $j\ge2$ and the three branch choices.
\end{enumerate}
\end{introF}
\end{samepage}
Parts~(i) and~(ii) are proved in Corollaries~\ref{cor:euler-gap-regimes}
and~\ref{cor:euler-weighted-rays}, with (ii) locating the three branches
about the ray $\eta=\xi/(2\pi j)$; closed-form values of $a(\sigma)$,
$a_1$, and $b_0$ appear in Appendix~\ref{app:scalar-phases},
\eqref{eq:euler-inner-shift-value} and \eqref{eq:euler-b0}.
Figure~\ref{fig:euler-quantile} illustrates Theorems~E and~F, and
computed curves and pair widths for the coupled Euler system appear in
Figure~\ref{fig:euler-numerical-overview}.

\paragraph{Organization of the paper.}
Sections~\ref{sec:abstract}--\ref{sec:sl} establish Theorems~A and~B
(boundary-bundle spectral flow, the Volterra $L^2$-frame criterion, and
interior cutting). Sections~\ref{sec:bessel} and~\ref{sec:renormalization}
construct the uniform Bessel frame across $s=0$ and $r=0$, yielding
Theorem~C and $a(\sigma)$. Section~\ref{sec:counting} handles lower-order
perturbations and noncommuting gyroscopic Hill systems via
$O(t)=\exp(-Kt)$. Section~\ref{sec:ere-applications} proves Theorem~D,
determines all Euler endpoints unconditionally, and derives Theorems~E
and~F.

In the appendices, Appendices~\ref{app:quadratic-proof}
and~\ref{app:reversible} give the Bessel estimates and half-period
reduction; Appendix~\ref{app:scalar-phases} evaluates $a(\sigma)$, $a_1$,
and $b_0$ and verifies $d_g,e_g\ne0$ in \cite{MSS06b} for each
$\delta>1/8$ and small matching cutoff;
Appendix~\ref{app:euler-critical-index} computes the Euler sector indices
for each fixed $0<\delta\le1/8$ and $e$ sufficiently close to one;
Appendix~\ref{app:domain-graph-example} gives the regular boundary
comparison; and Appendix~\ref{app:boundary-applications} treats a Hill
family and a two-ended transmission problem.
\FloatBarrier

\section{Continuous boundary coordinates and the index reduction}
\label{sec:abstract}

\paragraph{Notation.}
\label{sec:intro-notation}

All operators are densely defined. For a closed operator
$A$ on a Hilbert space $H$, $\D(A)$ is its domain, $A^*$ its
adjoint and $\gr A=\{(u,Au):u\in\D(A)\}$ its graph in
$H\oplus H$. The graph norm is
$(\|u\|^2+\|Au\|^2)^{1/2}$. Throughout, $\mathcal{CF}^{\mathrm{sym}}(H)$ and
$\mathcal{CF}^{\mathrm{sa}}(H)$ carry the gap topology.
For a closed subspace $M$, $P_M$ is its orthogonal projection, except when a nonorthogonal decomposition
$G^{\max}=G\oplus E$ is explicitly used: then $P_E$ denotes the
projection onto $E$ along $G$, defined on $G^{\max}$.
The operator gap is $\|P_{\gr A}-P_{\gr B}\|$; gap continuity
means norm continuity of $P_{\gr A}$. Put $J(x,y)=(-y,x)$ on
$H\oplus H$. For finite-dimensional quadratic forms, $n_-$ and $n_+$
denote the negative and positive inertia.
The boundary map $\gamma_A:\D(A^*)\to\beta(A)$ is the quotient
map; we abbreviate it to $\gamma_s$ for $A=A_s$.
$\Lag(V,\omega)$ is the set of Lagrangian subspaces of the
symplectic space $(V,\omega)$.

For self-adjoint $T$, $\sigma(T)$ is its spectrum and
$n_-(T)=\dim\operatorname{ran}\mathbf1_{(-\infty,0)}(T)$ is its
negative index, where $\mathbf1_{(-\infty,0)}(T)$ is the spectral
projection; the value $+\infty$ is allowed. The numbers
$\lambda_j(Q)$ are the eigenvalues of a symmetric matrix $Q$,
repeated with multiplicity, and $x_+=\max\{x,0\}$.
$\sf$ counts negative-to-positive crossings positively, while
$\Mas$ denotes the Cappell--Lee--Miller Maslov
index with the Green-form convention used in
\eqref{eq:sf}. Norms without subscripts are Hilbert,
Euclidean matrix or operator norms according to the argument;
$\|\cdot\|_2$ is an $L^2$ norm and $\|\cdot\|_{\rm HS}$ a
Hilbert--Schmidt norm.

In the regularization sections $s\ge0$ is the cutoff
parameter, $q\in\R$ a scalar inverse-square coefficient, $Q$ its
constant symmetric matrix analogue, and $B_s(t)$ a lower-order
matrix remainder. The signed parameter $r$ is defined by
\eqref{eq:q}; $\xi=[a:b]\in\mathbb{RP}^1$ is a real boundary
line, with $(a,b)\ne(0,0)$ considered up to nonzero scaling.
In the ERE section $N$ is the number of bodies, $d=2N-4$ the
essential configuration dimension, $e\in[0,1)$ the eccentricity,
and $\omega\in\C$, $|\omega|=1$, a Floquet multiplier.
We denote the Euler mass parameter by $\delta$ and the distinct Lagrange
mass parameter by $\beta$. Similarly, $L(e)=\log(1/(1-e))$ denotes the
logarithmic eccentricity scale (not an operator), and $C(R)$ in
\eqref{eq:ere-C} is the numerical counting coefficient, whereas $C(A)$
in \eqref{eq:boundary-facts} is the Cauchy-data subspace; local notation
is reintroduced whenever its meaning changes.

The abstract boundary theory is formulated on a real Hilbert space,
complexified when spectral parameters enter. Floquet operators for
general $\omega$ act on complex $L^2$ spaces with complex multiplicity,
which coincides with the real multiplicity at $\omega=\pm1$. Unless
stated otherwise, $O(1)$ denotes a bound independent of the limiting
parameter, with any further uniformity in mass or boundary data
specified in the corresponding theorem.

\subsection{Boundary quotient, graph complement, and Cauchy data}

The following facts establish the algebraic and symplectic equivalence between
the boundary quotient theory of \cite[Section 2]{HPW25} (see also \cite{Fur04,BZ18})
and concrete graph complements for a single operator, including the case of a nontrivial
minimal kernel. Here $\ind A=\dim\ker A-\dim\ker A^*$. For $A\in\mathcal{CF}^{\mathrm{sym}}(H)$, the
\emph{Green form} on $\D(A^*)$ is
\[
 \mathcal G_A(u,v)
 :=\langle u,A^*v\rangle-\langle A^*u,v\rangle,
 \qquad u,v\in\D(A^*).
\]
Equivalently, it is the restriction to $\gr A^*$ of
the ambient symplectic form
\[
 \widetilde\omega((u,x),(v,y))
 :=\langle u,y\rangle-\langle x,v\rangle,
\]
since
$\mathcal G_A(u,v)
=\widetilde\omega((u,A^*u),(v,A^*v))$.
The next lemma shows that this form descends to a
symplectic form on $\D(A^*)/\D(A)$.

\begin{lemma}\label{lem:boundary-facts}
Let $A\in\mathcal{CF}^{\mathrm{sym}}(H)$, and set
$G=\gr A$ and $G^{\max}=\gr A^*$. Equip $H\oplus H$ with the standard symplectic form
$\widetilde\omega((u_1,y_1),(u_2,y_2))=\ip{u_1}{y_2}-\ip{y_1}{u_2}$.
\begin{enumerate}[label=(\roman*)]
\item The Green form induces a symplectic form on $\beta(A)=\D(A^*)/\D(A)$ with
\begin{equation}\label{eq:boundary-facts}
 \dim\beta(A)=-2\ind A,\qquad C(A)=\gamma\ker A^*\in\Lag(\beta(A),\omega_A),
\end{equation}
where $\gamma:\D(A^*)\to\beta(A)$ is the quotient map. The radical of $\widetilde\omega$
on $G^{\max}$ is $G$. For any direct sum $G^{\max}=G\oplus E$, the map
\[
 \gamma_E:E\longrightarrow\beta(A),
 \qquad
 \gamma_E(u,A^*u)=u+\D(A),
\]
is a symplectic isomorphism from
$(E,\widetilde\omega|_E)$ onto $(\beta(A),\omega_A)$.
Equivalently, $\gamma_E=\gamma\circ(\pi_1|_E)$, where
$\pi_1:H\oplus H\to H$, $\pi_1(u,y)=u$, is the
first-coordinate projection.
\item The space of homogeneous solutions $\mathcal K=\ker A^*\oplus 0 \subset G^{\max}\cap(H\oplus 0)$
is isotropic in $H\oplus H$. Its projection along $G$ onto $E$,
\begin{equation}\label{eq:cauchy-proj}
 C = P_E(\mathcal K),
\end{equation}
is a Lagrangian subspace of $(E,\widetilde\omega|_E)$, and satisfies $\gamma_E(C)=C(A)$.
\item There is a bijective correspondence between self-adjoint extensions of $A$ and
Lagrangian subspaces $N\subset E$ (or $\Lambda=\gamma_E(N)\subset\beta(A)$), given by
$\gr A_N = G\oplus N$. Every such extension is Fredholm with index zero, and
\begin{equation}\label{eq:kernel-boundary}
 \ker A_N/\ker A\simeq C(A)\cap\Lambda \simeq C\cap N.
\end{equation}
\end{enumerate}
These relations are summarized by the commutative diagram:
\begin{equation}\label{eq:cd-boundary}
\begin{tikzcd}[row sep=large, column sep=large]
\mathcal K = \ker A^*\oplus 0 \ar[r, hook] \ar[d, "\pi_1"'] & G^{\max} = G\oplus E \ar[r, "P_E"] \ar[d, "\pi_1"] & (E, \widetilde\omega|_E) \ar[d, "\gamma_E"', "\simeq"] \\
\ker A^* \ar[r, hook] & \D(A^*) \ar[r, "\gamma"'] & (\beta(A), \omega_A)
\end{tikzcd}
\end{equation}
\end{lemma}
\begin{proof}
Since $A$ is symmetric, $\ker A\subset\ker A^*$, which gives $\ind A=\dim\ker A-\dim\ker A^*\le 0$;
hence the codimension $-2\ind A\ge 0$ is nonnegative.
Equip $\D(A)$ and $\D(A^*)$ with their graph norms. The inclusion
$j:\D(A)\to\D(A^*)$ has closed range and satisfies $A=A^*j$.
By Fredholm index additivity \cite{Kat80}, $\ind j=\ind A-\ind A^*=2\ind A$.
Since $j$ is injective, $\dim(G^{\max}/G)=\dim(\D(A^*)/\D(A))=-\ind j=-2\ind A$,
so any direct complement $E$ has dimension $-2\ind A$.

We first prove the graph statements: the symplectic structure of $E$ in (i) and the Cauchy data in (ii).
A vector $(u,A^*u)\in G^{\max}$ belongs to the radical of $\widetilde\omega$ on $G^{\max}$
if and only if $\ip{u}{A^*v}-\ip{A^*u}{v}=0$ for all $v\in\D(A^*)$, which characterizes
$u\in\D(A^{**})=\D(A)$ and $A^*u=Au$, i.e., $(u,A^*u)\in G$.
Hence $\operatorname{rad}(G^{\max},\widetilde\omega)=G$.
Because $G\cap E=\{0\}$, the restricted form $\widetilde\omega|_E$ is nondegenerate,
making $(E,\widetilde\omega|_E)$ a symplectic vector space of dimension $-2\ind A$.

For (ii), $\mathcal K=\ker A^*\oplus 0\subset H\oplus 0$ is isotropic in $H\oplus H$
because $\widetilde\omega|_{H\oplus 0}\equiv 0$.
Since $G$ is the radical of $\widetilde\omega$ on $G^{\max}$, the projection $P_E$ along $G$
preserves the symplectic form: for any $k_1,k_2\in\mathcal K$,
\[
 \widetilde\omega(P_Ek_1,P_Ek_2)=\widetilde\omega(k_1,k_2)=0.
\]
Hence $C=P_E(\mathcal K)$ is isotropic in $E$.
The kernel of $P_E|_{\mathcal K}$ is $\mathcal K\cap G=\ker A\oplus 0$, which gives
\[
 \dim C=\dim\ker A^*-\dim\ker A=-\ind A=\tfrac12\dim E.
\]
Therefore $C$ is Lagrangian in $E$, without an injectivity assumption.

With (i)--(ii) established on the graph complement $E$, the diagram \eqref{eq:cd-boundary}
commutes directly because $\pi_1(G)=\D(A)=\ker\gamma$.
Consequently $\gamma_E=\gamma\circ\pi_1|_E:(E,\widetilde\omega|_E)\to(\beta(A),\omega_A)$ is a canonical
symplectic isomorphism that maps $C$ onto $C(A)=\gamma\ker A^*$.
Under this isomorphism, the remaining quotient statements in (i) and the self-adjoint extension
correspondence in (iii) follow directly. Indeed, $(u,0)\in\gr A_N = G\oplus N$ if and only if
$u\in\ker A_N$; projecting along $G$ yields $P_E(u,0)\in N\cap C$, and the kernel of $u\mapsto P_E(u,0)$
on $\ker A_N$ is $G\cap\mathcal K = \ker A\oplus 0$, which yields $\ker A_N/\ker A \simeq C\cap N \simeq C(A)\cap\Lambda$.
Each extension is Fredholm with index $\ind A + \dim N = \ind A - \ind A = 0$.
\end{proof}
These facts concern one operator. For a continuous family, the index formula
uses the Maslov index of the pair $(N_s,C_s)$ in the graph complements
$E_s\subset H\oplus H$.

\subsection{Continuous graph complements}

\begin{theorem}\label{thm:canonical}
Let $X$ be a topological space and let $v\mapsto A_v\in\mathcal{CF}^{\mathrm{sym}}(H)$, $v\in X$,
be gap-continuous, and set $G_v=\gr A_v$ and
$G_v^{\max}=\gr A_v^*$. Suppose that
\[
 G_v^{\max}=G_v\oplus E_v,
\]
where $E_v$ is a continuous family of finite-dimensional subspaces
of $H\oplus H$. Let $P_{E_v}:G_v^{\max}\to E_v$ be the projection
along $G_v$. Then:
\begin{enumerate}[label=(\roman*)]
\item $(E,\widehat\omega)$, with
$\widehat\omega_v=\widetilde\omega|_{E_v}$, is a continuous
symplectic bundle of locally constant rank $-2\ind A_v$. Moreover,
$\pi_1|_{E_v}$ is injective and
$\D(A_v^*)=\D(A_v)\oplus\pi_1(E_v)$.
\item For Lagrangian subspaces $N_v\subset E_v$, the corresponding
self-adjoint realizations $B_v$ with $\gr B_v=G_v\oplus N_v$ are
gap-continuous if and only if $N_v$ is continuous in $H\oplus H$.
\item If $\ker A_v=\{0\}$ for every $v$, then $\ker A_v^*$ varies
continuously as a subspace of $H$, and
$C_v=P_{E_v}(\ker A_v^*\times\{0\})$
is a continuous Lagrangian subbundle of $E$.
\end{enumerate}
\end{theorem}

\begin{proof}
The adjoint identity $G_v^{\max}=JG_v^\perp$ proves continuity of the
maximal graphs. Since $G_v^{\max}=G_v\oplus E_v$, Lemma~\ref{lem:boundary-facts}(i)
gives $\dim E_v = -2\ind A_v$, which is locally constant.
Because $E_v \subset G_v^{\max}$ is a subspace of an operator graph and graphs meet $\{0\}\times H$ trivially,
the projection $\pi_1|_{E_v}:E_v\to H$ is injective.
Applying $\pi_1$ to $G_v^{\max} = G_v \oplus E_v$ yields the domain decomposition $\D(A_v^*) = \D(A_v) \oplus \pi_1(E_v)$.
The restricted forms $\widehat\omega_v=\widetilde\omega|_{E_v}$
are everywhere nondegenerate by Lemma~\ref{lem:boundary-facts}(i).
For any local continuous frame $\xi_1(v),\dots,\xi_m(v)$ of $E_v$ in $H\oplus H$,
the matrix entries $\Omega_{ij}(v)=\widetilde\omega(\xi_i(v),\xi_j(v))$
are continuous scalar functions of $v$, so $(E,\widehat\omega)$ is a continuous
symplectic vector bundle, and a continuous symplectic Gram--Schmidt process produces
local Darboux trivializations, proving (i).

The graph of the extension in (ii) is the direct sum $\gr B_v=G_v\oplus N_v$.
Because $G_v$ and $E_v$ are continuous direct summands of $G_v^{\max}$,
gap continuity of $B_v$ is equivalent to continuity of $N_v$ in $H\oplus H$.
This proves the equivalence without requiring a separately chosen domain frame.

For (iii), Lemma~\ref{lem:boundary-facts}(ii) ensures that each fiber
$C_v=P_{E_v}(\ker A_v^*\oplus 0)$ is a Lagrangian subspace of $E_v$.
To establish continuity, let $\mathcal H_0=H\oplus 0$.
A vector $e\in E_v$ belongs to $C_v$ if and only if $e+g\in\mathcal H_0$
for some $g\in G_v$, which yields the geometric description
\[
 C_v = E_v \cap (G_v + \mathcal H_0).
\]
Since $\ran A_v^*=H$ by the closed-range theorem, the sum satisfies
\[
 E_v + (G_v + \mathcal H_0) = G_v^{\max} + \mathcal H_0 = H\oplus H,
\]
showing that $E_v$ and $G_v+\mathcal H_0$ are transverse in $H\oplus H$.
The subspaces $G_v+\mathcal H_0=H\oplus\ran A_v$ vary continuously:
since $\ker A_v=\{0\}$, the second projection $(u,A_vu)\mapsto A_vu$ on
$G_v$ is an injective Fredholm map of locally constant cokernel
dimension, so its ranges vary continuously.
Because the intersection has constant dimension $\dim C_v=-\ind A_v$,
$C_v$ varies continuously as a subbundle of $E_v$.
\end{proof}

\begin{remark}
For example, $E_v^{\rm orth}=G_v^{\max}\cap G_v^\perp$ satisfies the
complement hypothesis, with continuity given by the orthogonal
projection formula
\begin{equation}\label{eq:projection-difference}
 P_{E_v^{\rm orth}}=P_{G_v^{\max}}-P_{G_v}.
\end{equation}
Under the injectivity hypothesis in (iii), surjectivity of $A_v^*$ and
gap continuity imply continuity of $\mathcal K_v=\ker A_v^*\times\{0\}$,
and $\ker A_v=\{0\}$ gives $G_v\cap\mathcal K_v=\{0\}$. Taking $V_v$ to
be the orthogonal complement of $G_v\oplus\mathcal K_v$ inside
$G_v^{\max}$ yields an adapted continuous complement
$E_v=\mathcal K_v\oplus V_v$ for which $C_v=\mathcal K_v$. Over a
parameter interval, every finite-dimensional vector bundle is trivial
and continuous symplectic frames exist globally.
\end{remark}

The theorem assumes gap continuity of $A_v$; establishing it for a
differential family is the separate task of Sections~\ref{sec:volterra}
and~\ref{sec:bessel}. Without injectivity, (i)--(ii) remain valid;
continuity of the Cauchy data requires additional information.

\subsection{Index conventions and the index formula in boundary coordinates}

We use the spectral-flow convention in which a negative-to-positive
crossing has sign $+1$ and zero belongs to the nonnegative side.
For a $C^1$-family of self-adjoint Fredholm operators $B(s)$ on a common domain, the crossing form at
a crossing $s$ (where $\ker B(s)\ne\{0\}$) is the quadratic form on $\ker B(s)$ defined by
\[
 \Gamma_B(s)[u] = \ip{u}{\dot B(s)u},
\]
following \cite{RS95,BLP05}. For regular crossings this means
\begin{equation}\label{eq:sf-convention}
 \sf(B;[0,1])=-\nm(\Gamma_B(0))
 +\sum_{s\in(0,1)}\operatorname{sign}\Gamma_B(s)
 +\np(\Gamma_B(1)).
\end{equation}
For a pair of $C^1$-curves of Lagrangian subspaces $(\lambda_s,\mu_s)$ in a symplectic
space $(V,\omega)$, the Robbin--Salamon crossing form at an intersection
$x\in\lambda_{s_0}\cap\mu_{s_0}$ is defined by
\[
 q_{s_0}(x) = \omega\bigl(x, \dot x_\mu(s_0) - \dot x_\lambda(s_0)\bigr)
\]
for curves $x_\lambda(s)\in\lambda_s$ and $x_\mu(s)\in\mu_s$ with $x_\lambda(s_0)=x_\mu(s_0)=x$,
following \cite{RS95,CLM94}. For regular crossings, the CLM convention formula is
\begin{equation}\label{eq:clm-convention}
 \Mas(\lambda,\mu,\omega;[0,1])=\np(q_0)
 +\sum_{s\in(0,1)}\operatorname{sign}q_s-\nm(q_1).
\end{equation}
The general continuous-path definitions and their homotopy and
additivity properties are taken from \cite{CLM94,RS95,BLP05,BZ18}.

For a continuous symplectic vector bundle $(E,\widehat\omega)$ of rank $2k$ over an interval $J$,
choosing a continuous symplectic frame field identifies each fiber
$(E_s,\widehat\omega_s)\cong(\R^{2k},\omega_0)$ and represents the Lagrangian subbundles
$N_s$ and $C_s$ as continuous curves $\widetilde N_s, \widetilde C_s\subset\Lag(\R^{2k},\omega_0)$,
defining the bundle Maslov index concretely by
\begin{equation}\label{eq:bundle-maslov}
 \Mas(N, C, \widehat\omega; J) := \Mas(\widetilde N, \widetilde C, \omega_0; J).
\end{equation}
This integer is independent of the choice of continuous frame (see \cite{BZ18}).
When the underlying symplectic form $\widehat\omega$ is clear from the context, we simply
write $\Mas(N, C; J)$ for brevity.

With the bundle Maslov index defined, we now state the main spectral-flow
formula for parameter-dependent self-adjoint realizations with varying domains.

\begin{theorem}[Index formula in boundary coordinates]\label{thm:sf}
Let $A_s$, $s\in[0,1]$, satisfy the hypotheses of Theorem~\ref{thm:canonical},
with $\ker A_s=\{0\}$ for every $s$, and choose any continuous family
of graph complements $E_s$ as in that theorem. Let $N_s\subset E_s$
be a continuous $\widehat\omega_s$-Lagrangian subbundle, and let $B_s$
be its self-adjoint realization, with $\gr B_s=\gr A_s\oplus N_s$.
If $C_s\subset E_s$ denotes the Cauchy data of Theorem~\ref{thm:canonical}(iii), then
\begin{equation}\label{eq:sf}
 \sf(B_s; s\in[0,1])=-\Mas(N_s,C_s,\widehat\omega_s; s\in[0,1]).
\end{equation}
\end{theorem}

To prove Theorem~\ref{thm:sf}, we first consider the special case of a positive path
where only the scalar shift $A+rI$ varies on a fixed domain (for which the boundary-quotient
formula was established in \cite[Section~2]{HPW25}).
In the language of direct graph complements, the scalar shift is represented by the
shear symplectomorphism
\begin{equation}\label{eq:shear}
 S_r(x,y) = (x, y+rx), \qquad S_r^{-1} = S_{-r},
\end{equation}
which preserves $\widetilde\omega$ on $H\oplus H$ and satisfies
\begin{equation}\label{eq:shear-graphs}
 \gr(A+rI) = S_r\gr A, \qquad \gr(A^*+rI) = S_r\gr A^*.
\end{equation}
Applying $S_r$ to the decomposition $\gr A^* = \gr A \oplus E$ automatically yields
\begin{equation}\label{eq:shear-decomp}
 \gr(A^*+rI) = \gr(A+rI) \oplus E_r, \qquad E_r = S_r E,
\end{equation}
and for any self-adjoint extension $B = A_N$ ($\gr B = \gr A \oplus N$),
$\gr(B+rI) = \gr(A+rI) \oplus N_r$ with $N_r = S_r N \subset E_r$.
The Cauchy data $C_r = P_{E_r}(\ker(A^*+rI)\times\{0\})\subset E_r$ pulls back
under $S_{-r}$ to $S_{-r}C_r = P_E\{(v,-rv) : v\in\ker(A^*+rI)\} \subset E$.

\begin{lemma}\label{lem:pencil}
Let $A\in\mathcal{CF}^{\mathrm{sym}}(H)$ be such that
$A+rI$ is injective and Fredholm for every $r$ in a compact
interval $J=[a,b]$, and let $B=A_N$ be a self-adjoint extension
with boundary condition $N\subset E$.
Then
\begin{equation}\label{eq:pencil}
 \sf(B+rI; r\in J)
 = -\Mas(N, S_{-r}C_r, \widehat\omega; r\in J)
 = -\Mas(N_r, C_r, \widehat\omega_r; r\in J).
\end{equation}
Under the canonical symplectic isomorphism $\gamma_E:(E,\widehat\omega)\to(\beta(A),\omega_A)$,
\eqref{eq:pencil} recovers the boundary-quotient formula
$-\Mas(\Lambda, \gamma\ker(A^*+rI), \omega_A; r\in J)$ of \cite[Section~2]{HPW25},
where $\Lambda=\gamma_E(N)$.
\end{lemma}
\begin{proof}
By Lemma~\ref{lem:boundary-facts}(iii), $\ker(B+rI)\simeq N\cap S_{-r}C_r$.
At a crossing $r_0$, let $u\in\ker(B+r_0I)$ and choose a smooth path
$u(r)\in\ker(A^*+rI)$ with $u(r_0)=u$.
Since $\gr A$ is the radical of $\widetilde\omega$ on $\gr A^*$ (Lemma~\ref{lem:boundary-facts}(i)),
the projection $e(r)=P_E(u(r),-ru(r))\in S_{-r}C_r$ satisfies
\[
 \widehat\omega(e(r_0), e(r))
 = \widetilde\omega\bigl((u,-r_0u), (u(r),-ru(r))\bigr)
 = (r_0-r)\ip{u}{u(r)}.
\]
Because $N$ is stationary ($\dot N=0$), the Robbin--Salamon crossing form is
\[
 q_{(N,S_{-r}C)}(r_0)[e] = \widehat\omega(e, \dot e(r_0)) = -\|u\|^2 = -\Gamma_B(r_0)[u].
\]
Since all crossings are regular with strictly opposite signs, the index conventions
\eqref{eq:sf-convention}--\eqref{eq:clm-convention} immediately yield the first equality in \eqref{eq:pencil}.
The second equality follows because $S_r:(E,\widehat\omega)\to(E_r,\widehat\omega_r)$ is a
symplectic isomorphism mapping $(N,S_{-r}C_r)$ to $(N_r,C_r)$.
\end{proof}

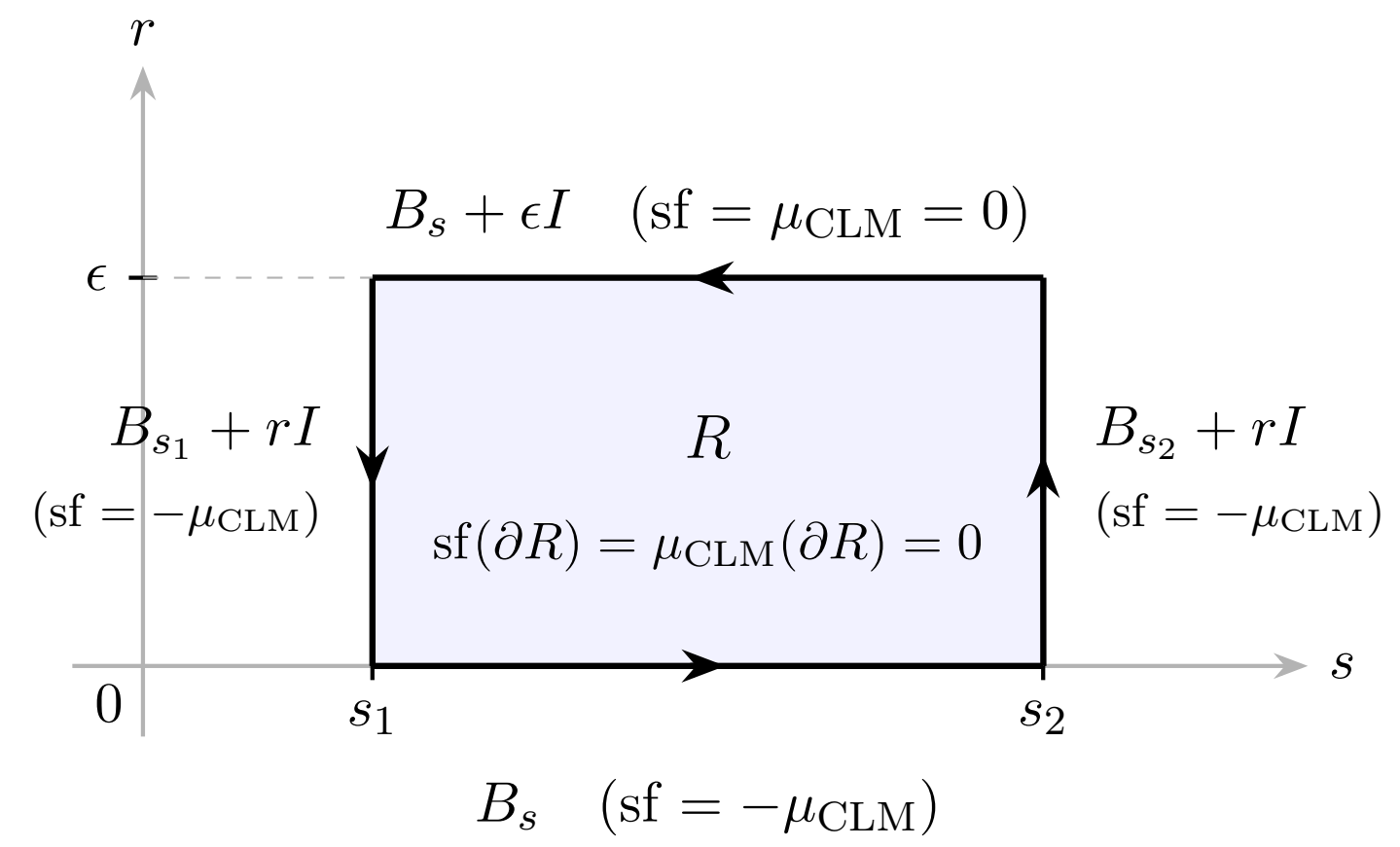
\begin{figure}[!tbp]
\centering
\begin{tikzpicture}[>=Stealth, scale=1.2]
  % Coordinate axes
  \draw[->, thick, gray!60] (-0.4, 0) -- (6.6, 0) node[right, black] {$s$};
  \draw[->, thick, gray!60] (0, -0.4) -- (0, 3.4) node[above, black] {$r$};

  % Axis ticks
  \draw[thick] (1.3, 0.08) -- (1.3, -0.08) node[below] {$s_1$};
  \draw[thick] (5.1, 0.08) -- (5.1, -0.08) node[below] {$s_2$};
  \draw[thick] (0.08, 2.2) -- (-0.08, 2.2) node[left] {$\epsilon$};
  \node[below left] at (0, 0) {$0$};

  % Corners
  \coordinate (BL) at (1.3, 0);
  \coordinate (BR) at (5.1, 0);
  \coordinate (TR) at (5.1, 2.2);
  \coordinate (TL) at (1.3, 2.2);

  % Shading
  \fill[blue!5] (BL) -- (BR) -- (TR) -- (TL) -- cycle;

  % Rectangle edges with counterclockwise orientation arrows
  \draw[very thick, ->] (BL) -- (3.3, 0);
  \draw[very thick] (3.1, 0) -- (BR);

  \draw[very thick, ->] (BR) -- (5.1, 1.2);
  \draw[very thick] (5.1, 1.0) -- (TR);

  \draw[very thick, ->] (TR) -- (3.1, 2.2);
  \draw[very thick] (3.3, 2.2) -- (TL);

  \draw[very thick, ->] (TL) -- (1.3, 1.0);
  \draw[very thick] (1.3, 1.2) -- (BL);

  % Dashed projection to r-axis
  \draw[dashed, gray!60] (0, 2.2) -- (TL);

  % Edge labels: short, clean, mathematical
  % Top
  \node[above=3pt] at (3.2, 2.2) {$B_s + \epsilon I \quad (\sf = \Mas = 0)$};

  % Bottom
  \node[below=18pt] at (3.2, 0) {$B_s \quad (\sf = -\Mas)$};

  % Left
  \node[left=6pt, align=right] at (1.3, 1.1) {$B_{s_1}+rI$\\[2pt] \footnotesize $(\sf = -\Mas)$};

  % Right
  \node[right=6pt, align=left] at (5.1, 1.1) {$B_{s_2}+rI$\\[2pt] \footnotesize $(\sf = -\Mas)$};

  % Center: region label and loop vanishing
  \node at (3.2, 1.3) {\large $R$};
  \node[below=2pt] at (3.2, 1.0) {\small $\sf(\partial R) = \Mas(\partial R) = 0$};

\end{tikzpicture}
\caption{The local homotopy rectangle $R = [s_1, s_2] \times [0, \epsilon]$ in the $(s,r)$-plane. Both families extend continuously over the contractible rectangle $R$, so after trivializing the symplectic bundle, homotopy invariance gives $\sf(\partial R) = \Mas(\partial R) = 0$. On the vertical edges $\sf = -\Mas$ holds by Lemma~\ref{lem:pencil}, and on the top edge both indices vanish by invertibility and transversality, establishing the formula $\sf = -\Mas$ along the bottom edge.}
\label{fig:homotopy-rectangle}
\end{figure}

\begin{proof}[Proof of Theorem~\ref{thm:sf}]
By path additivity, it suffices to prove the formula locally.
At any $s_0\in[0,1]$, because $B_{s_0}$ has discrete spectrum near zero,
we first choose $\epsilon>0$ such that $B_{s_0}+\epsilon I$ is invertible
and $A_{s_0}+rI$ is injective Fredholm for all $r\in[0,\epsilon]$.
Since invertibility and Fredholm injectivity are open conditions in the
gap topology \cite{Kat80}, there exists a neighborhood $[s_1,s_2]$ of $s_0$
such that $B_s+\epsilon I$ remains invertible for all $s\in[s_1,s_2]$ (ensuring
invertibility along the top edge) and $A_s+rI$ remains injective Fredholm on the
rectangle $R=[s_1,s_2]\times[0,\epsilon]$.

Since $R$ is contractible, the continuous symplectic bundle is trivializable over $R$.
Hence both the bundle Maslov index (independent of trivialization \cite{BZ18}) and the
spectral flow are homotopy invariant on $R$, so
\[
 \sf(\partial R)=0,\qquad \Mas(\partial R)=0.
\]
Along the boundary $\partial R$ (see Figure~\ref{fig:homotopy-rectangle}):
\begin{itemize}
\item On the vertical edges ($s\in\{s_1,s_2\}$ fixed, oriented with increasing $r$), Lemma~\ref{lem:pencil} gives $\sf = -\Mas$ directly.
\item On the top edge ($r=\epsilon$, oriented with increasing $s$), invertibility of $B_s+\epsilon I$ and triviality of the intersection $N\cap C=\{0\}$ (by \eqref{eq:kernel-boundary}) give $\sf(\text{top})=\Mas(\text{top})=0$.
\end{itemize}
Orienting each segment along the positive coordinate axes ($s$ increasing on bottom and top; $r$ increasing on left and right), the counterclockwise boundary decomposition is $\partial R = \text{bottom} + \text{right} - \text{top} - \text{left}$. Homotopy invariance $\sf(\partial R) = \Mas(\partial R) = 0$ then yields
\[
 \sf(\text{bottom}) = \sf(\text{left})-\sf(\text{right}) = -\Mas(\text{left})-(-\Mas(\text{right})) = -\Mas(\text{bottom}),
\]
which proves \eqref{eq:sf} locally. Summing over a finite partition of $[0,1]$ completes the proof.
\end{proof}

\begin{remark}[Geometric interpretation: infinite-dimensional Maslov index and reduction]
\label{rem:graph-boundary-geometry}
Theorem~\ref{thm:sf} has a coordinate-free geometric interpretation. Set $\mathcal H_0=H\times\{0\}$. More precisely, the graph path belongs to the
Fredholm--Lagrangian Grassmannian relative to $\mathcal H_0$,
\[
 \mathcal{FL}_{\mathcal H_0}(H\oplus H)
 :=
 \bigl\{L\in\Lag(H\oplus H,\widetilde\omega):
              (L,\mathcal H_0)\text{ is a Fredholm pair}\bigr\}.
\]
Indeed, self-adjointness of $B_s$ makes $\gr B_s$ Lagrangian,
while
\[
 \gr B_s\cap\mathcal H_0=\ker B_s\times\{0\},
 \qquad
 \gr B_s+\mathcal H_0=H\times\operatorname{ran}B_s.
\]
Since $B_s$ is Fredholm, this intersection is finite-dimensional
and this sum is closed and has finite codimension.
Gap continuity of $B_s$ therefore gives a continuous path
in $\mathcal{FL}_{\mathcal H_0}(H\oplus H)$.
The ambient Maslov index below is understood in this space,
with the same sign and endpoint conventions as in
Theorem~\ref{thm:sf}.
In the ambient Hilbert space $H\oplus H$ equipped with the standard symplectic form
$\widetilde\omega((u_1,v_1),(u_2,v_2))=\ip{u_1}{v_2}-\ip{v_1}{u_2}$,
the spectral flow $\sf(B_s)$ of the self-adjoint family $B_s$ can be viewed as an
infinite-dimensional Maslov index counting intersections of the moving operator
graphs $\gr B_s$ with the fixed horizontal Lagrangian subspace $\mathcal H_0$ \cite{RS95,Fur04}:
\[
 \sf(B_s) = -\Mas(\gr B_s, \mathcal H_0; H\oplus H).
\]
The minus sign reflects the fact that the moving Lagrangian occupies the first slot:
for a $C^1$ family on a common domain, at a crossing where
$\langle u,\dot B_su\rangle>0$ for every nonzero $u\in\ker B_s$,
the crossing form \eqref{eq:clm-convention} evaluates, with $x=(u,0)$, to
$q(x)=\widetilde\omega(x,-\dot x_{\gr B})=-\langle u,\dot B_su\rangle<0$.
For varying domains the graph crossing form is the intrinsic formulation.

By Theorem~\ref{thm:canonical}, the chosen graph complements form a
finite-dimensional symplectic subbundle $(E_s,\widehat\omega_s)$ in $H\oplus H$.
The adjoint graph $\gr A_s^*$ is coisotropic in $(H\oplus H,\widetilde\omega)$,
with $\gr A_s$ as its characteristic isotropic radical
(Lemma~\ref{lem:boundary-facts}(i)). Since $A_s$ is injective,
a complement adapted to the homogeneous solutions is obtained
by choosing a continuous complement $V_s$ such that
$\gr A_s^*=\gr A_s\oplus\mathcal K_s\oplus V_s$, and setting
\[
 E_s = \mathcal K_s \oplus V_s, \qquad \mathcal K_s = \ker A_s^* \times \{0\}.
\]
In this complement, the Cauchy data coincides with the homogeneous solution space
$\mathcal K_s$ itself (without applying $P_E$), while the operator graph
$\gr B_s = \gr A_s \oplus N_s$ projects along $\gr A_s$ to the boundary condition $N_s\subset E_s$.
Their intersection in $E_s$ is the physical kernel:
\[
 N_s \cap \mathcal K_s = \gr B_s \cap (\ker A_s^* \times \{0\}) = \ker B_s \times \{0\}.
\]
Thus Theorem~\ref{thm:sf} reduces the ambient infinite-dimensional
Maslov index of $(\gr B_s, \mathcal H_0)$ in $H\oplus H$ to the
finite-dimensional Maslov index of $(N_s, \mathcal K_s)$ in the subbundle $E_s$:
\[
 \sf(B_s) = -\Mas(\gr B_s, \mathcal H_0; H\oplus H) = -\Mas(N_s, \mathcal K_s, \widehat\omega_s; [0,1]).
\]
Geometrically, this is the coisotropic reduction of the ambient graph Maslov index
along $\gr A_s$ to the boundary-data bundle, where crossings are measured by intersections
between the boundary condition $N_s$ and the homogeneous solution space $\mathcal K_s$.
\end{remark}

\begin{remark}
\label{rem:degenerate-charts}
When an endpoint of a Sturm--Liouville problem passes from regular ($s>0$)
to singular ($s=0$), the classical evaluation map $\gamma_0(u)=(u'(0),u(0))$
may cease to be well defined in the graph topology. As a result, a continuous
Lagrangian section $s\mapsto N_s\subset E_s$ in the boundary bundle
$(E,\widehat\omega)$ can appear singular in unscaled classical coordinates
(for instance, with a Robin parameter blowing up like $O(1/s)$).
Theorem~\ref{thm:canonical}(ii) shows that this coordinate blow-up does not
affect operator continuity: as long as $N_s$ varies continuously in an
admissible graph complement $E_s\subset H\oplus H$, the realizations $B_s$
remain gap-continuous (hence norm-resolvent continuous) up to $s=0$.
\end{remark}

For an elementary derivation of the domain-to-graph correspondence,
see Appendix~\ref{app:domain-graph-example}.

\subsection{A Robin family in boundary coordinates}

\begin{figure}[!tbp]
\centering
\begin{tikzpicture}[
  x=0.85cm,y=0.85cm,
  >=stealth,
  cauchy/.style={blue!70!black,thick,dashed},
  boundary/.style={red!75!black,thick},
  every node/.style={font=\small}
]

% Common coordinate axes.
\foreach \a in {0,4.8,9.6}{
  \begin{scope}[shift={(\a,0)}]
    \draw[->,gray!65] (-1.65,0)--(1.7,0)
      node[right] {$x$};
    \draw[->,gray!65] (0,-1.65)--(0,1.7)
      node[above] {$p$};
  \end{scope}
}

% Before the crossing.
\begin{scope}
  \node at (0,2.45) {$s=\frac14,\quad\alpha_s=-\frac32$};
  \draw[cauchy] (-1.35,1.35)--(1.35,-1.35);
  \draw[boundary] (-0.95,1.425)--(0.95,-1.425);
  \fill (0,0) circle (1.5pt);
  \node at (0,-2.1) {$\ker B_{1/4}=\{0\}$};
\end{scope}

% At the crossing.
\begin{scope}[shift={(4.8,0)}]
  \node at (0,2.45) {$s=\frac12,\quad\alpha_s=-1$};
  \draw[red!75!black,line width=1.8pt]
    (-1.35,1.35)--(1.35,-1.35);
  \draw[cauchy] (-1.35,1.35)--(1.35,-1.35);
  \node[fill=white,inner sep=2pt] at (0.85,1.15)
    {$\ell_s=\ell_C$};
  \node at (0,-2.1) {$\dim\ker B_{1/2}=1$};
\end{scope}

% After the crossing.
\begin{scope}[shift={(9.6,0)}]
  \node at (0,2.45) {$s=\frac34,\quad\alpha_s=-\frac12$};
  \draw[cauchy] (-1.35,1.35)--(1.35,-1.35);
  \draw[boundary] (-1.45,0.725)--(1.45,-0.725);
  \fill (0,0) circle (1.5pt);
  \node at (0,-2.1) {$\ker B_{3/4}=\{0\}$};
\end{scope}

% Legend.
\draw[cauchy] (1.1,-3.05)--(2,-3.05);
\node[anchor=west] at (2.1,-3.05)
  {$\ell_C:\ p=-x$};

\draw[boundary] (5.8,-3.05)--(6.7,-3.05);
\node[anchor=west] at (6.8,-3.05)
  {$\ell_s:\ p=\alpha_s x$};

\end{tikzpicture}
\caption{The Robin family in reduced boundary coordinates
$(x,p)=(y(0),y'(0))$, after imposing $y(1)=0$.
The fixed line $\ell_C$ consists of the left boundary
values of solutions of $-y''=0$ satisfying the right
Dirichlet condition. The moving line $\ell_s$ represents
$y'(0)=\alpha_s y(0)$, with $\alpha_s=2s-2$.
The lines coincide only at $s=1/2$, where
$\ker B_{1/2}=\operatorname{span}\{1-t\}$.
This is the unique upward eigenvalue crossing, giving
$\sf(B_s;s\in[0,1])=1$.
The picture represents the two-dimensional symplectic
reduction of the four-dimensional boundary space $E$.}
\label{fig:robin-boundary-crossing}
\end{figure}
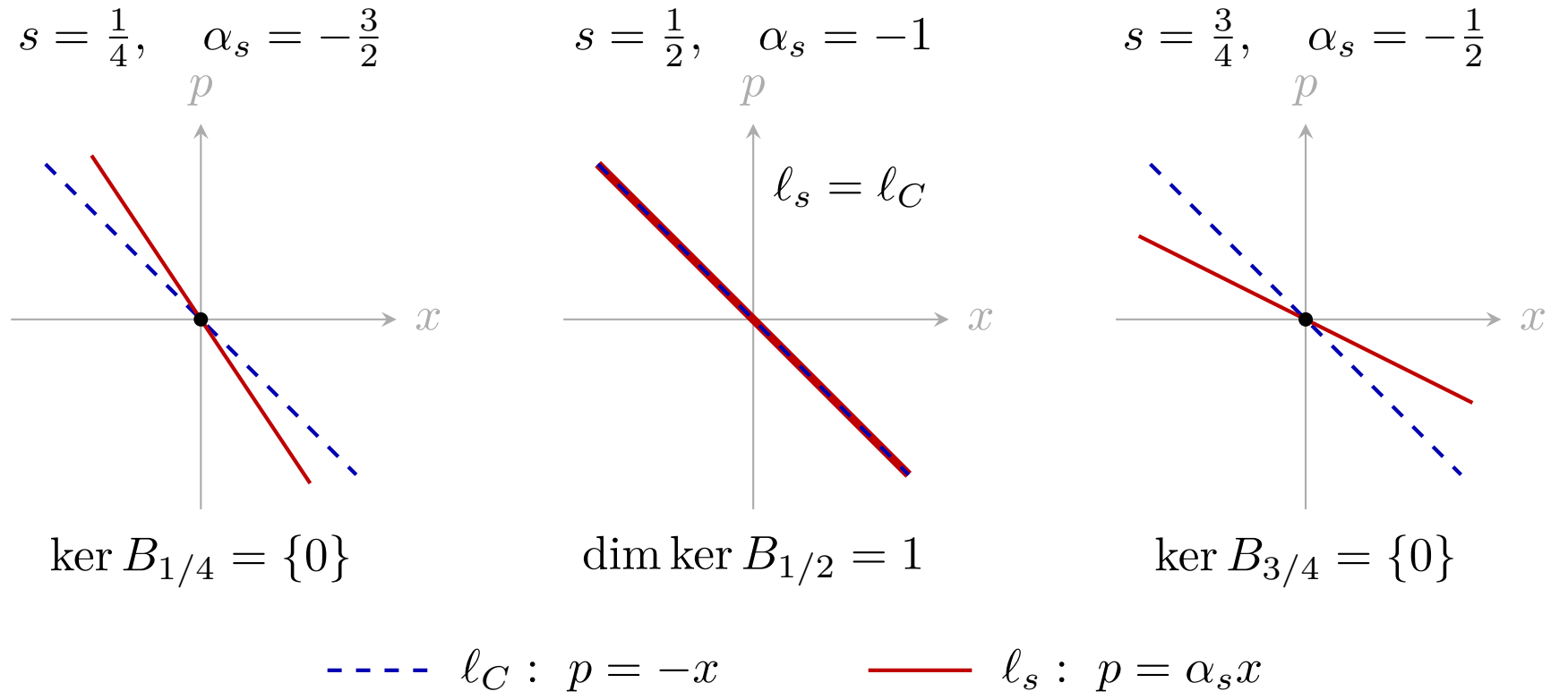

\begin{example}\label{ex:robin-boundary-family}
Let $H=L^2(0,1;\mathbb R)$ and take the constant minimal operator
$A_s=A$, $Au=-u''$, with domain
$\D(A)=\{u\in H^2(0,1):u(0)=u'(0)=u(1)=u'(1)=0\}$. Its adjoint is
$A^*u=-u''$ on $\D(A^*)=H^2(0,1)$, so $\ker A=\{0\}$,
$\ker A^*=\operatorname{span}\{1,t\}$, and $\ind A=-2$.
The surjective trace map $\tau u=(u(0),u'(0),u(1),u'(1))$ has kernel
$\D(A)$ and identifies $\beta(A)=\D(A^*)/\D(A)$ with $\mathbb R^4$; in
coordinates $(x_0,p_0,x_1,p_1)$, integration by parts gives the Green
form
\[
 \omega\bigl((x_0,p_0,x_1,p_1),(\widetilde x_0,\widetilde p_0,\widetilde x_1,\widetilde p_1)\bigr)
 =x_0\widetilde p_0-p_0\widetilde x_0-x_1\widetilde p_1+p_1\widetilde x_1.
\]

Let $\mathcal P_3$ be the cubic polynomials and set
$E_s=E=\{(p,-p''):p\in\mathcal P_3\}$. By existence and uniqueness of
cubic Hermite interpolation, $\gr A^*=\gr A\oplus E$, and
$(p,-p'')\mapsto\tau p$ identifies $(E,\widehat\omega)$ symplectically
with $(\mathbb R^4,\omega)$ (see Appendix~\ref{app:domain-graph-example}
for the explicit Hermite basis). For $\alpha_s=2s-2$, $s\in[0,1]$, the
Robin realizations
\[
 B_su=-u'',\qquad
 \D(B_s)=\{u\in H^2(0,1):u'(0)=\alpha_s u(0),\ u(1)=0\}
\]
correspond to the continuous Lagrangian subspaces $N_s\subset E$ with
trace coordinates $\tau(N_s)=\{(a,\alpha_s a,0,b):a,b\in\mathbb R\}$,
while the fixed Cauchy-data space $C_s=C=\{(u,0):u\in\ker A^*\}$ has
trace coordinates $\tau(C)=\{(a,b,a+b,b):a,b\in\mathbb R\}$.
By Theorem~A, continuity of $N_s$ implies norm-resolvent continuity of
$B_s$.

A nonzero vector lies in $N_s\cap C$ precisely when $b=\alpha_s a$ and
$a+b=0$, which occurs only at $s=1/2$ ($\alpha_{1/2}=-1$), with
$\ker B_{1/2}=\operatorname{span}\{1-t\}$. The quadratic forms
$\mathfrak b_s[u]=\int_0^1|u'|^2\,dt+\alpha_s|u(0)|^2$ share the common
form domain $\{u\in H^1(0,1):u(1)=0\}$ and satisfy
$\dot{\mathfrak b}_{1/2}[u_0]=2|u_0(0)|^2=6>0$ for the unit eigenfunction
$u_0(t)=\sqrt3(1-t)$ (cf.~\cite{LS24} for first-order eigenvalue
asymptotics under varying boundary conditions). There is thus a unique
upward crossing (Figure~\ref{fig:robin-boundary-crossing}), giving
\[
 \sf(B_s;s\in[0,1])=-\Mas(N_s,C,\widehat\omega;s\in[0,1])=1.
\]
\end{example}
\FloatBarrier

\section{A quantitative Volterra criterion for graph continuity}\label{sec:volterra}

Convergence of singular potentials away from the endpoints does not
control the varying minimal domains in operator or form norm. To obtain
gap continuity without evaluating endpoint traces, we invert the minimal
operator $T_v$ on its range $\mathcal R_v=(\ker T_v^*)^\perp$ by a
Volterra integral operator $K_v$ built directly from an $L^2$ symplectic
fundamental frame $F_v$, allowing both endpoints to be singular.
Expressing the inverse graph as the transverse intersection
$\gr T_v^{-1}=\gr K_v\cap(\mathcal R_v\oplus H)$ then turns
$L^2$-continuity of $F_v$ into a quantitative gap estimate for $T_v$ and
yields compact embedding of its maximal graph domain.

\begin{theorem}[Quantitative Volterra criterion]\label{lem:volterra}
Let $T_v$ be the minimal realizations in $H=L^2((0,1),\mathbb R^n)$ of
\[
\tau_vy=-(P_vy'+Q_vy)'+Q_v^Ty'+R_vy,\qquad v\in X,
\]
where $X$ is a topological space, $P_v,Q_v$ are locally $C^1$,
$R_v$ is locally continuous, $P_v=P_v^T$ is invertible and $R_v=R_v^T$.
Both endpoints may be singular; write $y_v^{[1]}=P_vy'+Q_vy$.
Suppose $\dim\ker T_v^*=2n$ and these kernels are gap-continuous in $H$.
Locally in $v$, assume that $P_v^{-1},Q_v,R_v$ are continuous in $C^0(I)$
for a fixed compact interval $I\Subset(0,1)$ with nonempty interior.

Then $T_v$ are gap-continuous, injective Fredholm operators of index $-2n$.
Their maximal graph domains embed compactly into $H$, and every closed
intermediate realization of $\tau_v-\lambda I_n$ is Fredholm for
$\lambda\in\mathbb C$ (after complexification).
Locally there is an $L^2$-continuous fundamental frame $F_v$ with
$\left(\begin{smallmatrix}-F_v^{[1]}\\F_v\end{smallmatrix}\right)$ symplectic, and
\begin{equation}\label{eq:quantitative-gap}
\bigl\|P_{\gr T_v}-P_{\gr T_{v'}}\bigr\|
\le C\,\|F_v-F_{v'}\|_{L^2}.
\end{equation}
In place of the kernel and local coefficient assumptions, one may
assume such a continuous symplectic frame directly.
\end{theorem}
\begin{proof}
We first construct the normalized frame from the kernel and local
coefficient assumptions. The adjoint is the maximal realization, with
$y,y_v^{[1]}\in AC_{\mathrm{loc}}$ and $\tau_vy\in H$ (no global $L^2$
condition on $y_v^{[1]}$ is needed). Fix $v_0$ and
$t_0\in\operatorname{int}I$, put $C_v=\ker T_v^*$, and let
$\mathcal S_v(I)$ denote the space of homogeneous solutions on $I$.
Restriction defines an isomorphism
\[
 \operatorname{res}_I:C_v\longrightarrow\mathcal S_v(I),
\]
since it is injective by uniqueness and both spaces have dimension $2n$.
Gap continuity of $C_v$ provides local continuous bases
\cite[Chapter~I, \S4.6]{Kat80}, whose restrictions form continuous bases
of $\mathcal S_v(I)$ because $\operatorname{res}_I$ is fixed, bounded,
and injective on $C_v$. Hence the image spaces $\mathcal S_v(I)$ vary
continuously in the gap topology, and the inverse restriction is
continuous and locally uniformly bounded by finite-dimensional matrix
inversion \cite[Chapter~I, \S4.4]{Kat80}. Letting $B_v$ be the position
component of the local fundamental matrix normalized by
$\left(\begin{smallmatrix}-B_v^{[1]}\\B_v\end{smallmatrix}\right)(t_0)=I_{2n}$,
local coefficient continuity makes $B_v$ continuous in $C^0(I)$, so
$F_v=(\operatorname{res}_I|_{C_v})^{-1}B_v$ is a global $L^2$-continuous
fundamental frame whose symplectic normalization persists throughout
$(0,1)$ by Green's identity. When the frame is given directly, this
preliminary step is omitted.

In either case, regard $F_v$ as an injective map $\R^{2n}\to H$
with range $C_v$.

\smallskip
\noindent\emph{Step 1: a Volterra right inverse.}
Put $J=\left(\begin{smallmatrix}0&-I_n\\ I_n&0\end{smallmatrix}\right)$
and define $K_v:H\to H$ by
\begin{equation}\label{eq:volterra}
(K_vf)(t)=\int_t^1F_v(t)JF_v(x)^Tf(x)\,dx,\qquad f\in H.
\end{equation}
The integral is well defined since $F_v^Tf\in L^1$;
no endpoint value of the frame is used.
Its square-integrable kernel makes $K_v$ Hilbert--Schmidt and hence compact on $H$.
The symplectic normalization of the frame gives
\[
F_v(t)JF_v(t)^T=0,\qquad
F_v^{[1]}(t)JF_v(t)^T=I_n.
\]
Local differentiation gives $K_vf,(K_vf)^{[1]}\in AC_{\mathrm{loc}}$ and
\[
\tau_vK_vf=f,\qquad f\in H.
\]
Thus $K_vf\in\mathcal D(T_v^*)$ and $T_v^*K_v=I_H$, so that
$\mathcal D(T_v^*)=\ker T_v^*\oplus\operatorname{ran}K_v$.
For $h\in C_c^\infty(0,1;\R^n)$, both $h$ and $K_v\tau_vh$ vanish
to the right of $\operatorname{supp}h$, and their difference solves the
homogeneous equation. Uniqueness for the regular initial-value problem
therefore gives $K_v\tau_vh=h$. Since $K_v$ is bounded, passage to the
minimal graph closure yields
\[
h=K_vT_vh,\qquad h\in\mathcal D(T_v).
\]
In particular, $\|h\|\leq\|K_v\|\|T_vh\|$, so $T_v$
is injective with closed range. As $T_v$ is densely
defined and closed,
\begin{equation}\label{eq:range}
\mathcal R_v:=\operatorname{ran}T_v=(\ker T_v^*)^\perp=C_v^\perp.
\end{equation}
For $f\in\mathcal R_v$, the identity $K_vT_vh=h$ shows that
$K_vf\in\mathcal D(T_v)$ and that $K_v$ restricts to the inverse of $T_v$.
Consequently $\operatorname{ind}T_v=-2n$ and
\begin{equation}\label{eq:graph-volterra}
\operatorname{gr}T_v
=\{(K_vf,f):f\in \mathcal R_v\}.
\end{equation}

\smallskip
\noindent\emph{Step 2: continuity and the gap estimate.}
Work near a fixed parameter $v_0$. The maps $F_v:\R^{2n}\to H$
have a locally uniform positive lower bound. The range-gap estimate
\cite[Proposition~3.16]{DW14}, applied to $F_v,F_{v'}$, and direct
comparison of the Volterra kernels therefore give
\[
 \|P_{\mathcal R_v}-P_{\mathcal R_{v'}}\|+\|K_v-K_{v'}\|_{\mathrm{HS}}
 \le C\|F_v-F_{v'}\|_{L^2}.
\]
By \eqref{eq:graph-volterra},
\[
 \operatorname{gr}T_v^{-1}
 =\operatorname{gr}K_v\cap(\mathcal R_v\oplus H).
\]
The minimum gap in Kato's sense satisfies
\[
 \gamma\bigl(\operatorname{gr}K_v,\mathcal R_v\oplus H\bigr)
 \ge (1+\|K_v\|^2)^{-1/2}
\]
by its definition \cite[Chapter~IV, (4.4)]{Kat80};
the two orders agree in a Hilbert space \cite[Chapter~IV, \S4.3]{Kat80}.
Thus these minimum gaps are locally bounded away from zero,
since $K_v$ is locally uniformly bounded.
Kato's intersection estimate \cite[Chapter~IV, Lemma~4.4]{Kat80},
applied in both directions, and invariance of the gap under interchange
of the two factors give
\[
 \operatorname{gap}(T_v,T_{v'})
 \le C\bigl(\|K_v-K_{v'}\|
       +\|P_{\mathcal R_v}-P_{\mathcal R_{v'}}\|\bigr)
 \le C\|F_v-F_{v'}\|_{L^2}.
\]
This proves \eqref{eq:quantitative-gap} and gap continuity.

\smallskip
\noindent\emph{Step 3: compactness and other realizations.}
The identity $h=K_vT_vh$ shows that the graph-norm embedding
$\mathcal D(T_v)\hookrightarrow H$ is compact. Since
$\dim\mathcal D(T_v^*)/\mathcal D(T_v)=4n$ by
Lemma~\ref{lem:boundary-facts}, the maximal graph domain and every
closed intermediate domain embed compactly into $H$ as well. Each
intermediate realization differs from $T_v$ by a finite-dimensional
domain extension and is therefore Fredholm; because a bounded spectral
shift is compact from its graph domain into $H$, it preserves both
Fredholmness and the index. In particular, every self-adjoint extension
has compact resolvent.
\end{proof}

The kernel condition and continuity near a regular interior point
supply the normalized frame; neither endpoint need be regular.
The direct frame criterion uses only $L^2$-continuity, rather than
convergence of the coefficients in a global multiplier or form topology. It is
unchanged under continuous symplectic changes of
frame, since these preserve the homogeneous solution
spaces and the normalization.

\begin{corollary}\label{cor:energy}
Under the hypotheses of Theorem~\ref{lem:volterra},
the shifted minimal operators $T_v-\lambda I_n$ are
jointly gap-continuous in $(v,\lambda)$ (after
complexification when $\lambda\in\C$). For real
$\lambda$ they remain injective Fredholm operators of
index $-2n$. Consequently, the abstract boundary
bundle and spectral-flow results apply to these
shifted families. In particular, a continuous family
of self-adjoint boundary conditions in the resulting
boundary bundle determines a norm-resolvent-continuous
family of self-adjoint realizations, to which
\eqref{eq:sf} applies.
\end{corollary}
\begin{proof}
Gap continuity is stable under norm-continuous bounded
perturbations, giving the first assertion. For real
$\lambda$, a solution $h\in\ker(T_v-\lambda I_n)$ satisfies
$h=\lambda K_vh$. Writing $c(t)=\int_t^1JF_v(x)^Th(x)\,dx$ gives
$c'=-\lambda JF_v^TF_vc$ and $c(1)=0$.
Since $F_v^TF_v\in L^1$, uniqueness for this integral equation gives
$c=0$ and hence $h=0$. The Fredholm
index is unchanged because the shift is compact from the
graph domain to $H$, and Theorems~\ref{thm:canonical} and \ref{thm:sf}
apply. For self-adjoint operators, gap continuity is
equivalent to norm-resolvent continuity at a nonreal
spectral parameter \cite{BLP05}.
\end{proof}

\section{Interior cutting and Cauchy trace coordinates}\label{sec:sl}

Like Section~\ref{sec:volterra}, the construction here avoids
singular endpoints by working on a regular interior subinterval and
establishing gap continuity through a transverse graph intersection in
$H\oplus H$. Here, however, we start from a general gap-continuous
minimal family that need not be injective, and use stronger $C^1$ local
regularity near an interior cut point $c$ to produce an injective cut
family and explicit Cauchy trace coordinates. Readers primarily
interested in the Bessel regularization and $N$-body applications may
proceed directly to Section~\ref{sec:bessel}: the scalar Bessel minimal
operators satisfy Section~\ref{sec:volterra} directly, the general
$N$-body counts at $e<1$ reduce to scalar comparison operators by the
min--max principle, and the coupled Euler family handles its limit-point
strong component through weighted energy estimates rather than an
interior cut.

Consider on $(a,b)$ the real matrix expression
\begin{equation}\label{eq:sl}
 \ell_su=-(P_su'+Q_su)'+Q_s^Tu'+R_su,
 \quad P_s=P_s^T,\quad R_s=R_s^T.
\end{equation}
Assume $P_s$ is invertible, $P_s,Q_s$ are locally $C^1$ and $R_s$
is locally continuous. Let $L_s$ be the minimal realization on $H=L^2((a,b),\R^n)$,
obtained by closing the expression on $C_c^\infty((a,b),\R^n)$, and let $L_s^*$
be the maximal realization. The ambient space $\mathcal{H}=H\oplus H$
carries the canonical Green symplectic form
\[
 \widetilde\omega\bigl((u,f),(v,g)\bigr) = \ip{u}{g}_H - \ip{f}{v}_H.
\]
By Green's identity \cite{Zet05,HPW25}, for maximal-domain elements with
quasi-derivatives $p_u = P_s u' + Q_s u$, Cauchy traces $z_u = (p_u, u)$, and
symplectic bracket $[u,v]_s = p_u^T v - u^T p_v$, one has
$\widetilde\omega\bigl((u,\ell_s u),(v,\ell_s v)\bigr) = [u,v]_s(b) - [u,v]_s(a)$
(with endpoint limits when necessary). In particular, $\gr L_s$ is an isotropic subspace of $\mathcal{H}$.

The abstract index theorem of Section~\ref{sec:abstract} requires
an injective minimal operator family, that is, $\ker A_s=\{0\}$.
If the uncut global minimal operator $L_s$ has a nontrivial
kernel, imposing vanishing Cauchy data at an interior regular
point $c$ makes the cut minimal operator injective: by ODE
uniqueness, every homogeneous solution with these data vanishes
on each subinterval. This restriction decreases the Fredholm
index by $2n$, where $n$ is the dimension of the dependent
variable. We establish the gap continuity and Fredholm index
formula of the cut minimal family directly through Cauchy
trace evaluation on the minimal graph.

Let $H_-=L^2((a,c),\R^n)$ and $H_+=L^2((c,b),\R^n)$, and let
$j:H\to H_-\oplus H_+$ be the restriction isometry. Write
\[
 \widehat L_s=L_s^- \oplus L_s^+
\]
for the direct sum of the minimal operators on the two subintervals.
No boundary matching is imposed in its maximal domain.

Fix a compact subinterval $I_1 \subset (a,b)$ with $c \in \operatorname{int} I_1$.
Locally in $s$, suppose
\begin{equation}\label{eq:local-coefficients}
 P_s,Q_s\longrightarrow P_{s_0},Q_{s_0}\quad\hbox{in }C^1(I_1),
 \qquad R_s\longrightarrow R_{s_0}\quad\hbox{in }C^0(I_1),
\end{equation}
and $P_s^{-1}$ is uniformly bounded on $I_1$.

\begin{theorem}[Continuity of the cut minimal family]\label{thm:cut}
Under \eqref{eq:local-coefficients}, if $s\mapsto L_s$ is a gap-continuous family of
symmetric Fredholm operators, then the cut minimal family $s\mapsto\widehat L_s$ is
gap-continuous, injective and Fredholm, with
\begin{equation}\label{eq:cut-index}
 \ind\widehat L_s=\ind L_s-2n.
\end{equation}
If $B_s$ is a gap-continuous family of self-adjoint extensions of $L_s$,
the matching realizations $jB_sj^{-1}$ have continuous Lagrangian data
$\widehat N_s$ in the admissible graph complement $\widehat E_s$.
Consequently, letting $\widehat C_s = P_{\widehat E_s}(\ker \widehat L_s^* \oplus 0)$ be the abstract Cauchy data of the kernel of the adjoint of the cut minimal operator, the index formula along any compact parameter interval is
\begin{equation}\label{eq:cut-sf}
 \sf(B_s)=-\Mas\bigl(\widehat N_s,\widehat C_s\bigr).
\end{equation}
\end{theorem}
\begin{remark}
In the proof we express the additional conditions created by
cutting at $c$ as constraints on a fixed ambient Hilbert space.
The ordinary Cauchy trace $T_su=z_u(c)$ is defined only on the
varying graph domain $\D(L_s)$, so the maps $T_s$ cannot be
compared directly in operator norm. To overcome this, Green's
identity and $2n$ fixed test functions supported on the left of
$c$ yield bounded maps
\[
\widetilde T_s:H\oplus H\longrightarrow\R^{2n}
\]
which depend continuously on $s$ in operator norm and satisfy
\[
\widetilde T_s(u,L_su)=M_sz_u(c),
\]
where $M_s$ is invertible. Hence the cut graph is represented by
\[
\gr\widehat L_s=\gr L_s\cap\ker\widetilde T_s.
\]
The regularity of $c$ allows arbitrary Cauchy data to be prescribed
by functions in $\D(L_s)$; equivalently,
$\widetilde T_s|_{\gr L_s}$ is surjective. This gives the
transversality needed to conclude that the intersections, and
therefore the cut graphs, vary continuously in the gap topology.
\end{remark}
\begin{proof}
Under the isometric splitting $j: H = L^2(a, b) \to L^2(a, c) \oplus L^2(c, b) = H_- \oplus H_+$,
the ambient space is $\mathcal{H} = H \oplus H \cong \mathcal{H}_- \oplus \mathcal{H}_+$.
Via Green's identity on $(a,c)$, the Cauchy trace at $c$ can be regarded as a bounded linear map
defined on the entire ambient space $\mathcal{H}$.
Specifically, choose $2n$ fixed smooth one-sided test functions
$\{v_i\}_{i=1}^{2n} \subset C_c^\infty((a, c], \R^n)$ supported in $I_1 \cap (a, c]$ whose Cauchy traces
$z_{v_i}(c^-)$ form a basis at the reference parameter (and hence nearby), and define
\[
 \widetilde T_s: \mathcal{H} \longrightarrow \R^{2n}, \qquad
 \widetilde T_s(u, f) = \Bigl( \widetilde\omega_{\mathcal{H}_-}\bigl((u, f), (v_i, \ell_s v_i)\bigr) \Bigr)_{i=1}^{2n} = \Bigl( \ip{(u, f)}{(\ell_s v_i, -v_i)}_{\mathcal{H}_-} \Bigr)_{i=1}^{2n}.
\]
Because each $v_i$ is fixed and smooth, the local regular condition \eqref{eq:local-coefficients}
ensures that $\ell_s v_i \to \ell_{s_0} v_i$ in $L^2(a,c)$.
Thus $s \mapsto \widetilde T_s$ is a norm-continuous family of bounded operators on the fixed
Hilbert space $\mathcal{H}$, and its kernel $\ker\widetilde T_s \subset \mathcal{H}$ is a closed subspace
of codimension $2n$ that varies continuously in the gap topology of $\mathcal{H}$.

By Green's identity on $(a,c)$, for any $u \in \D(L_s)$,
\[
 \widetilde T_s(u, L_s u) = \bigl( [u, v_i]_s(c^-) \bigr)_{i=1}^{2n} = M_s z_u(c),\qquad z_u(c)=(P_s(c)u'(c)+Q_s(c)u(c),u(c)),
\]
where $M_s$ is an invertible, continuously varying matrix of bracket coordinates.
Under the splitting $j$, an uncut function $u \in \D(L_s)$ belongs to the cut minimal domain
$\D(\widehat L_s) = \D(L_s^-) \oplus \D(L_s^+)$ if and only if its Cauchy trace vanishes at $c$: $z_u(c) = 0$.
The cut minimal graph is therefore the intersection
\begin{equation}\label{eq:cut-kernel}
 \gr\widehat L_s = \gr L_s \cap \ker\widetilde T_s.
\end{equation}
The transversality of this intersection in $\mathcal{H}$, $\gr L_s + \ker\widetilde T_s = \mathcal{H}$,
is equivalent to the surjectivity of the restriction $\widetilde T_s|_{\gr L_s} = M_s z_u(c)$ onto $\R^{2n}$.
Because $c$ is an interior regular point with $P_s(c)$ invertible, smooth functions in $\D(L_s)$
prescribe arbitrary 1-jets $(u'(c), u(c))$, so $\widetilde T_s(\gr L_s) = \R^{2n}$, establishing transversality.
Since both $\gr L_s$ and $\ker\widetilde T_s$ vary continuously in the gap topology of $\mathcal{H}$ and
meet transversally, their intersection $\gr\widehat L_s$ is gap-continuous in $\mathcal{H}$
\cite{Kat80}.

The codimension-$2n$ drop yields \eqref{eq:cut-index}, while ODE uniqueness on each subinterval gives $\ker\widehat L_s = \{0\}$.
Self-adjoint extensions matching across $c$ define continuous Lagrangian data with interface cancellation on the diagonal, so Theorem~\ref{thm:sf} directly yields \eqref{eq:cut-sf}.
\end{proof}

For families with $\dim\ker L_s^*=2n$, Theorem~\ref{lem:volterra}
verifies the gap-continuity and Fredholm hypotheses from continuity of
the kernels and of the coefficients near a regular interior point,
without requiring either endpoint to be regular.
The cut construction also requires the stronger local continuity
of $P_s,Q_s$ in $C^1$ and of $R_s$ in $C^0$ near the cut point,
as specified in \eqref{eq:local-coefficients}.

\begin{remark}[Piecewise smooth coefficients and transmission problems]
\label{rem:transmission}
Retain the gap-continuity and Fredholm hypotheses on $L_s$ and assume
the one-sided versions of \eqref{eq:local-coefficients} at $c$.
Although Theorem~\ref{thm:cut} is formulated for coefficients that are
$C^1$ across $c$, the same conclusions hold for transmission problems where
$P_s, Q_s, R_s$ are piecewise regular on each side of $c$ and may jump across
the interface. In that setting, the unbroken minimal operator
$L_s$ is defined on functions whose Cauchy traces match across $c$,
$z_u(c^-) = z_u(c^+)$, whereas the cut minimal operator $\widehat L_s = L_s^- \oplus L_s^+$
imposes vanishing Cauchy traces on both sides and its adjoint drops the matching constraint.
Because the functional $\widetilde T_s$ in \eqref{eq:cut-kernel} uses one-sided test functions
$\{v_i\} \subset C_c^\infty((a,c], \R^n)$ supported on the left subinterval, coefficient jumps
across $c$ do not enter its evaluation. Moreover, for any $u \in \D(L_s)$ the matching
condition $z_u(c^-) = z_u(c^+)$ implies that $\widetilde T_s(u, L_s u) = 0$ forces the Cauchy
traces to vanish on both sides of $c$. Consequently, the transversality of $\gr L_s$ and
$\ker\widetilde T_s$, the gap continuity of $\widehat L_s$, the index formula \eqref{eq:cut-sf},
and the trace-coordinate formulas of Theorem~\ref{thm:trace-index} extend without change to
piecewise regular transmission problems.
\end{remark}

\subsection{Explicit Cauchy trace coordinates}

To compute the Maslov index in \eqref{eq:cut-sf} and express boundary
conditions explicitly, we cut the interval at a regular interior point
$c$ so that each half-interval has one regular endpoint carrying ordinary
Cauchy traces. Decomposing each half-interval maximal graph into the
minimal graph, homogeneous solutions, and smooth test functions supported
near $c$ splits Green's form orthogonally into a singular-endpoint
bracket and the standard symplectic form at $c$, yielding concrete
trace coordinates. More precisely, as in \cite{HPW25}, on $(a,c]$ one has
\[
 \gr A_s^* = \gr A_s +
 \bigl((\ker A_s^*\oplus\{0\})
       + \gr(A_s^*|_Z)\bigr),
\]
where $Z=\{\chi(t)(\alpha+(t-c)\beta):\alpha,\beta\in\mathbb R^n\}$ is the
fixed $2n$-dimensional space of cut-off affine functions supported near $c$
(with a smooth cutoff $\chi$ equal to $1$ in a neighborhood of $c$ and
supported in $I_1\cap(a,c]$). Since $P_s(c)$ is invertible, the Cauchy trace
$z\mapsto z_c(z)=(P_s(c)\beta+Q_s(c)\alpha,\alpha)$ restricts to an
isomorphism $Z\xrightarrow{\sim}\mathbb R^{2n}$ for every $s$, and in
particular $\gr A_s\cap\gr(A_s^*|_Z)=\{0\}$.

The homogeneous solutions represent the Cauchy-data
space, while $Z$ provides the independent boundary directions at the
regular endpoint. This description makes it possible to express the
abstract boundary conditions, and hence the Maslov index, through
ordinary Cauchy traces. For any $u_1, u_2 \in \gr A_s^*$ with singular components $v_1, v_2 \in \ker A_s^*$, the fact that $Z$ vanishes near $a$ ensures that Green's identity decouples the boundary form without cross terms:
\[
 \omega(u_1, u_2) = [u_1, u_2](c) - [u_1, u_2](a^+) = \omega_0\bigl(z_c(u_1), z_c(u_2)\bigr) - [v_1, v_2]_a.
\]

The Green form splits into the singular-endpoint bracket and the
standard symplectic form on the regular Cauchy data. Since the latter
is nondegenerate, it suffices to choose a complement $V_s$ of the
radical of $-[\cdot,\cdot]_a$ on $\ker A_s^*$, on which this bracket
is nondegenerate. Identifying homogeneous solutions with their graph
vectors $(v,0)$, we obtain
\[
 E_s=V_s\oplus\gr(A_s^*|_Z),\qquad
 (E_s,\omega)\simeq
 (V_s,-[\cdot,\cdot]_a)\oplus(\mathbb R^{2n},\omega_0).
\]
Here the isomorphism sends the graph representative of $v+z$, with
$v\in V_s$ and $z\in Z$, to $(v,z_c(v+z))$. Its regular coordinate
is the total Cauchy trace, including the contribution of $v$.
The analogous right-hand identification is
\[
 (E_s^+,\omega)\simeq
 (\mathbb R^{2n},-\omega_0)\oplus(V_s^+,[\cdot,\cdot]_b).
\]
The complements can be chosen continuously in $s$. Indeed, injectivity
and Fredholmness imply continuity of the adjoint kernels $\ker A_s^*$, and the
regular-endpoint traces vary continuously by Theorem~\ref{thm:cut}.
On each adjoint kernel the singular-endpoint bracket agrees with
the bracket at $c$ and has constant rank $2(k-n)$, where
$k=\dim\ker A_s^*$. To see the rank assertion, in $\beta(A_s)$ let
$\mathcal Z_s=\gamma(Z)$ be the symplectic regular-trace subspace
of dimension $2n$. The Lagrangian $C(A_s)=\gamma\ker A_s^*$ has
dimension $k$. ODE uniqueness at $c$ gives
$C(A_s)\cap\mathcal Z_s^{\omega_{A_s}}=\{0\}$, hence
$C(A_s)+\mathcal Z_s=\beta(A_s)$. Projection along $\mathcal Z_s$
therefore maps $C(A_s)$ onto the symplectic space
$\mathcal Z_s^{\omega_{A_s}}$, of dimension $2(k-n)$.
The singular-endpoint form is its pullback, with the endpoint sign,
so it has the asserted rank. Its radical therefore has constant dimension,
so a continuous nondegenerate complement $V_s$ can be selected along
the parameter interval.

Applying this construction on both sides of $c$ gives a continuous
symplectic complement of the cut minimal graph,
\begin{equation}\label{eq:cut-explicit-E}
\begin{aligned}
\gr\widehat L_s^*&=\gr\widehat L_s\oplus\widehat E_s,\\
\widehat E_s&=E_s^-\oplus E_s^+\\
&=V_s^-\oplus V_s^+
\oplus\gr\bigl((L_s^-)^*|_{Z^-}\bigr)
\oplus\gr\bigl((L_s^+)^*|_{Z^+}\bigr).
\end{aligned}
\end{equation}
The half-interval trace maps then yield the symplectic identification
\begin{equation}\label{eq:trace-iso-abstract}
\begin{aligned}
\Phi_s:(\widehat E_s,\widehat\omega_s)&\xrightarrow{\sim}
\mathcal V_s
:=\bigl(V_s^-,-[\cdot,\cdot]_a\bigr)
\oplus(\R^{2n},\omega_0)
\oplus(\R^{2n},-\omega_0)
\oplus\bigl(V_s^+,[\cdot,\cdot]_b\bigr),\\
u&\longmapsto
\bigl(v_{u^-}^-,z_c(u^-),z_c(u^+),v_{u^+}^+\bigr).
\end{aligned}
\end{equation}

The dimensions of the singular boundary spaces follow from the
half-interval deficiency indices. Since $\ker L_s^\pm=\{0\}$,
Theorem~\ref{thm:canonical}(i) gives
\begin{equation}\label{eq:dim-V}
 m_\pm:=\dim V_s^\pm
 =2\bigl(\dim\ker(L_s^\pm)^*-n\bigr).
\end{equation}
In particular, $0\le m_\pm\le 2n$ are even, and
\begin{equation}\label{eq:dim-total}
 \dim\widehat E_s=\dim\mathcal V_s
 =m_-+4n+m_+=2\dim\ker\widehat L_s^*.
\end{equation}

For computations, choose continuous symplectic bases
$\{g_i^-\}_{i=1}^{m_-}$ of $V_s^-$ and
$\{g_j^+\}_{j=1}^{m_+}$ of $V_s^+$.
The nondegenerate endpoint brackets, together with the regular
Cauchy data at $c$, give the concrete trace map
\[
\begin{aligned}
 \Tr_s:\gr\widehat L_s^*&\longrightarrow
 \R^{m_-}\oplus\R^{2n}\oplus\R^{2n}\oplus\R^{m_+},\\
 (u,\widehat L_s^*u)&\longmapsto
 \Bigl(\bigl([u^-,g_i^-]_a\bigr)_i,\,
 z_c(u^-),\,z_c(u^+),\,
 \bigl([u^+,g_j^+]_b\bigr)_j\Bigr).
\end{aligned}
\]
The minimal graph has zero trace, and the restriction of $\Tr_s$
to $\widehat E_s$ is an isomorphism. Thus the abstract Cauchy-data
Lagrangian can be computed simply by evaluating $\Tr_s$ on
homogeneous solutions, without constructing a graph projection.
With these normalized frames the symplectic form is represented by the fixed block-diagonal matrix
\[
 \Omega=\Omega_a\oplus\omega_0\oplus(-\omega_0)\oplus\Omega_b,
\]
where, with $G_a=([g_i^-,g_j^-]_a)$ and
$G_b=([g_i^+,g_j^+]_b)$, one has
$\Omega_a=G_a^{-1}$ and $\Omega_b=-G_b^{-1}$.

We now combine the abstract index formula with these concrete
boundary coordinates and its specialization to separated problems.

\begin{theorem}[Trace coordinates and spectral flow]\label{thm:trace-index}
Let $B_s$ be a gap-continuous family of self-adjoint extensions of
$\widehat L_s$ as in Theorem~\ref{thm:cut}. Choose continuous
symplectic frames for the singular boundary components and let
\[
\Tr_s:\gr\widehat L_s^*\longrightarrow\R^{m_-+4n+m_+}
\]
be the resulting trace map. Its target carries the fixed symplectic
form
\[
\Omega=\Omega_a\oplus\omega_0\oplus(-\omega_0)\oplus\Omega_b.
\]
Set
\[
\mathcal L_s=\Tr_s(\gr B_s),\qquad
\mathcal K_s=\Tr_s(\ker\widehat L_s^*\times\{0\}).
\]
Then:
\begin{enumerate}[label=(\roman*)]
\item The spaces $\mathcal L_s$ and $\mathcal K_s$ are continuous
Lagrangian subspaces, and
\begin{equation}\label{eq:cut-concrete-clm}
\sf(B_s)=-\mu_{\rm CLM}(\mathcal L_s,\mathcal K_s;
\R^{m_-+4n+m_+},\Omega).
\end{equation}
\item Suppose $B_s$ imposes the matching condition
$z_u(c^-)=z_u(c^+)$ at the cut, its outer boundary conditions
are separated, and the corresponding half-interval solution spaces
have continuous Lagrangian Cauchy traces
$W_s^\pm(c)\subset\R^{2n}$. Then
\begin{equation}\label{eq:cut-clm}
\sf(B_s)=-\mu_{\rm CLM}\bigl(W_s^+(c),W_s^-(c);
\R^{2n},\omega_0\bigr).
\end{equation}
\end{enumerate}
\end{theorem}

\begin{proof}
\emph{(i)} By Theorem~\ref{thm:sf}, spectral flow is minus the
Maslov index of the extension and Cauchy-data Lagrangians in the
abstract boundary bundle of $\widehat L_s$. The trace map
annihilates $\gr\widehat L_s$ and induces, by the preceding
construction, a continuous symplectic isomorphism of this bundle
onto $(\R^{m_-+4n+m_+},\Omega)$. It sends the extension
Lagrangian to $\mathcal L_s$ and the Cauchy-data Lagrangian to
$\mathcal K_s$. Their continuity and Lagrangian character follow,
and symplectic invariance of the CLM index \cite{CLM94} gives
\eqref{eq:cut-concrete-clm}.

\smallskip
\noindent\emph{(ii)} The separated outer conditions select
half-interval homogeneous solution spaces whose traces at $c$ are
$W_s^-(c)$ and $W_s^+(c)$. Reducing the boundary pair by these
outer conditions leaves, in
$(\R^{2n}\oplus\R^{2n},\omega_0\oplus(-\omega_0))$,
the matching diagonal $\Delta_c$ and the product
$W_s^-(c)\times W_s^+(c)$. The reduction and diagonal identities
for the CLM index \cite{CLM94,HPW25} yield
\[
\mu_{\rm CLM}\bigl(\Delta_c,W_s^-(c)\times W_s^+(c)\bigr)
=\mu_{\rm CLM}\bigl(W_s^+(c),W_s^-(c)\bigr),
\]
with the indicated symplectic forms. Substituting into (i)
proves \eqref{eq:cut-clm}.
\end{proof}
\begin{remark}
The passage from the abstract Maslov formula to
\eqref{eq:cut-concrete-clm} relies on more than the existence of
boundary traces. The map $\Tr_s$ must annihilate
$\gr\widehat L_s$ and induce an isomorphism from the boundary
quotient $\gr\widehat L_s^*/\gr\widehat L_s$ onto the finite-dimensional
trace space. Green's identity identifies the induced symplectic form
with the boundary form, and choosing continuous symplectic frames for
the singular components makes its matrix $\Omega$ independent of
$s$. These properties ensure that the extension and Cauchy-data
Lagrangians, and hence their CLM index, are faithfully represented
by concrete traces without explicitly constructing a graph projection.
For regular Schr\"odinger operators on decomposed manifolds, an analogous
interface Maslov-index formula for transmission conditions is developed
in \cite{CJM17}.
\end{remark}

\section{Bessel regularization through the critical parameter}\label{sec:bessel}

We consider the singular inverse-square equation
\[
-y''(t)+\frac{q}{t^2}y(t)=0,\qquad 0<t<1.
\]
The endpoint $t=0$ is limit-circle precisely for $q<3/4$.
To describe this range with a single real parameter, including the
critical coefficient $q=-1/4$, we set
\begin{equation}\label{eq:q}
 q=q(r):=r|r|-\frac14,\qquad -\infty<r<1.
\end{equation}
This is a continuous bijection onto $(-\infty,3/4)$.
The cases $r>0$, $r=0$, and $r<0$ correspond, respectively, to
power-law, logarithmic, and oscillatory solutions near $t=0$.
A real fundamental system on $(0,1]$, normalized by
$W(u_r,w_r)=u_rw_r'-u_r'w_r=1$, is
\begin{equation}\label{eq:euler-frame}
 \begin{aligned}
 u_r(t)&=\tfrac12(t^{1/2+r}+t^{1/2-r}),&
 w_r(t)&=\frac{t^{1/2+r}-t^{1/2-r}}{2r},&&r>0,\\
 u_0(t)&=t^{1/2},&
 w_0(t)&=t^{1/2}\log t,&&r=0,\\
 u_r(t)&=t^{1/2}\cos(r\log t),&
 w_r(t)&=t^{1/2}\frac{\sin(r\log t)}{r},&&r<0.
 \end{aligned}
\end{equation}
On each compact set of $r<1$, these functions have a common
$L^2(0,1)$ majorant. In particular, near $r=0$ the divided difference
is bounded by $Ct^{1/2-\varepsilon}(1+|\log t|)$ for any fixed small
$\varepsilon>0$. The limiting frame is therefore continuous also
at the repeated indicial root.

\begin{remark}\label{rem:thresholds}
The two thresholds play distinct roles: semiboundedness of
$-d^2/dt^2+q/t^2$ holds for $q\ge-1/4$ (the Hardy threshold), whereas
the endpoint $t=0$ is limit-circle for $q<3/4$ and limit-point for
$q\ge3/4$. Note that $q(r)$ remains negative for $0\le r<1/2$, so
$r\ge0$ does not correspond to a repulsive potential. For $r<0$ the
minimal operator is not semibounded; our gap-continuity results still
cover this range, although finite-Morse-index formulas do not apply.
\end{remark}
Consider the minimal operator $G_{s,r}$ associated with
\begin{equation}\label{eq:quadratic}
 g_{s,r}=-\frac{d^2}{dt^2}+\frac{q(r)}{t^2+s^2}.
\end{equation}
At $s=0$ it gives the same indicial roots as \eqref{eq:euler-frame}.

\begin{theorem}[Gap continuity of the minimal Bessel operators]\label{thm:bessel2}
The map $(s,r)\mapsto G_{s,r}$ is jointly gap-continuous on
$[0,1]\times(-\infty,1)$. Every $G_{s,r}$ is symmetric Fredholm,
and its maximal graph domain embeds compactly in $L^2(0,1)$.
\end{theorem}
\begin{proof}
Appendix~\ref{app:quadratic-proof} constructs a real Wronskian-one
frame for \eqref{eq:quadratic}. The estimate \eqref{eq:Phi-bound}
controls both a solution and its derivative in the indicial parameter,
so the divided difference \eqref{eq:normalized-frame} is jointly
continuous in $L^2(0,1)$ through $s=r=0$.
The Volterra criterion, Theorem~\ref{lem:volterra}, now applies.
\end{proof}

\begin{proposition}[The limit-point obstruction]\label{prop:threshold}
Fix $q\ge 3/4$. For the quadratic family
$-d^2/dt^2+q/(t^2+s^2)$ on $(0,1)$, the minimal operators
for $s>0$ cannot converge in the gap topology, as $s\downarrow0$,
to the minimal operator at $s=0$.
\end{proposition}
\begin{proof}
For $s>0$ the two regular ends give Fredholm index $-2$.
At $s=0$, put $\nu=\sqrt{q+1/4}\ge1$. Of the two zero solutions
$t^{1/2\pm\nu}$, only $t^{1/2+\nu}$ belongs to $L^2(0,1)$.
On $C_c^\infty(0,1)$, positivity of $q$ and Poincar\'e's inequality
give $\langle T_0u,u\rangle\ge\pi^2\|u\|^2$.
By closure the minimal operator is injective with closed range.
Its index is therefore $-1$, whereas the Fredholm index is locally
constant in the gap topology. For $n$ diagonal copies the two
indices are $-2n$ and $-n$.
\end{proof}

We will also use the following consequence of compact graph embedding.

\begin{lemma}\label{lem:strong-compact}
Let $A_v$ be a family of closed densely defined operators on a Hilbert space $H$
that is gap-continuous at $v_0$. Suppose that the domain $\D(A_{v_0})$, equipped with the
graph norm, embeds compactly into $H$. Let $B_v\in\mathcal{B}(H)$ be locally uniformly
bounded operators converging strongly to $B_{v_0}$ as $v\to v_0$. Then $A_v+B_v$
is gap-continuous at $v_0$. If each $A_v$ is Fredholm of index $k$ and has compact
graph-domain embedding into $H$, then each $A_v+B_v$ is Fredholm of index $k$.
\end{lemma}
\begin{proof}
Let $P_v$ denote the orthogonal projection onto $\gr A_v$ and $\pi_1$
the projection onto the first coordinate. Since $\pi_1P_{v_0}$ is
compact, strong convergence $B_v\to B_{v_0}$ implies
$(B_v-B_{v_0})\pi_1P_{v_0}\to0$ in norm, which combined with
$P_v\to P_{v_0}$ gives
\begin{equation}\label{eq:strong-diff-graph}
 \|(B_v-B_{v_0})\pi_1P_v\|\longrightarrow0 .
\end{equation}
The uniformly invertible shears $(x,y)\mapsto(x,y+B_vx)$ carry
$\gr A_v$ onto $\gr(A_v+B_v)$, and on $\gr A_v$ their difference from
the shear for $B_{v_0}$ tends to zero in norm by
\eqref{eq:strong-diff-graph}; comparing the images of unit vectors in
both directions establishes gap continuity at $v_0$. At each $v$, $B_v$
is compact from the graph domain into $H$, so the Fredholm property and
index are preserved.
\end{proof}

\section{Prescribing the singular limit by boundary conditions}\label{sec:renormalization}
Gap continuity of the minimal operators does not by itself select a
self-adjoint extension. Write $(u_{s,r},w_{s,r})$ for the real Wronskian-one fundamental frame of
$g_{s,r}y=0$ constructed in Appendix~\ref{app:quadratic-proof}
(see \eqref{eq:normalized-frame}), built from the hypergeometric solution
$\Phi(z,s,t)$ in \eqref{eq:Phi} and reducing to \eqref{eq:euler-frame}
at $s=0$, and recall from Section~\ref{sec:bessel} that $G_{s,r}$ is the
minimal operator of $g_{s,r}$.
A line $\xi=[a:b]\in\mathbb{RP}^1$ selects the left endpoint condition
through the zero solution $h_{s,r,\xi}=a u_{s,r}+b w_{s,r}$. Because
$(u_{s,r},w_{s,r})$ varies continuously in $L^2(0,1)$ up to $s=0$,
holding $\xi$ fixed defines a continuous boundary condition on both sides
of the Hardy threshold and at the singular limit itself.

\begin{theorem}[Selection of singular boundary conditions]
\label{thm:boundary-selection}
For $s\in[0,1]$, $r<1$, and $\xi\in\mathbb{RP}^1$, define
\begin{equation}\label{eq:selected-domain}
 A_{s,r,\xi}=G_{s,r}^*\big|_{\D_{s,r,\xi}},\qquad
 \D_{s,r,\xi}=\bigl\{y\in\D(G_{s,r}^*):
 y(1)=0,\ [y,h_{s,r,\xi}](0)=0\bigr\}.
\end{equation}
These are self-adjoint operators with compact resolvent, depending
jointly continuously on $(s,r,\xi)$ in the norm-resolvent topology.
At $s=0$, they exhaust all separated self-adjoint realizations with
Dirichlet condition at $1$, and convergence as $s\downarrow0$ is
uniform for $(r,\xi)$ in compact subsets of
$(-\infty,1)\times\mathbb{RP}^1$.

For any continuous path $\theta\mapsto(r_\theta,\xi_\theta)$
whose two singular endpoint operators are invertible, one has,
for all sufficiently small $s>0$,
\begin{equation}\label{eq:selected-sf}
 \sf(A_{s,r_\theta,\xi_\theta};\theta\in[0,1])
 =\sf(A_{0,r_\theta,\xi_\theta};\theta\in[0,1]).
\end{equation}
The continuity and spectral-flow conclusions also hold after adding
jointly strongly continuous, locally uniformly bounded real
multiplication operators $B_{s,r,\xi}$, provided the perturbed
singular endpoint operators are invertible in \eqref{eq:selected-sf}.
\end{theorem}

\begin{proof}
The Wronskian-one frame $(u_{s,r},w_{s,r})$ gives left boundary
coordinates
\[
 a_{s,r}(y)=-[y,w_{s,r}](0),\qquad
 b_{s,r}(y)=[y,u_{s,r}](0).
\]
These traces depend continuously on $(s,r)$ in the maximal graph
topology, since Green's identity expresses each endpoint bracket through
the regular trace at $1$ and an $L^2$ pairing of $G_{s,r}^*y$ with a
frame vector, while the frames vary continuously in $L^2$ and in $C^1$
near $1$. Moreover, cutting off the frame vectors near $0$ (with cutoff
derivatives supported away from the singular endpoint) and adding
regular trace functions near $1$ supplies continuous right inverses for
the full boundary trace map.

In these coordinates the condition $[y,h_{s,r,\xi}](0)=0$ reads
\[
 a\,b_{s,r}(y)-b\,a_{s,r}(y)=0,
 \qquad \xi=[a:b],
\]
which specifies a continuous Lagrangian line; together with $y(1)=0$,
this yields continuous Lagrangian boundary data exhausting all separated
singular conditions at $s=0$. Theorems~\ref{thm:canonical}
and~\ref{lem:volterra} then give gap continuity, or equivalently
norm-resolvent continuity for self-adjoint operators \cite{BLP05}, and
compactness of the maximal graph embedding yields compact resolvent,
with joint continuity on compact parameter sets implying uniform
convergence as $s\downarrow0$.

For a fixed parameter path, the regularized paths deform continuously to
the singular path while keeping their endpoints invertible for
sufficiently small $s$, so homotopy invariance of spectral flow gives
\eqref{eq:selected-sf}. Finally, bounded multiplication perturbations
are handled by Lemma~\ref{lem:strong-compact}, since compactness of the
graph-domain embeddings strengthens strong continuity to gap continuity
while leaving both the domains and resolvent compactness unchanged.
\end{proof}

For fixed $(r,\xi)$ the spectral projections of the regularizations
therefore converge in norm on every bounded interval whose endpoints
avoid the limiting spectrum.

\paragraph{The physical lines in the Bessel frame.}
For $s>0$, write $h=h_{s,r,\xi}$. The boundary condition is
\begin{equation}\label{eq:quadratic-cauchy-bc}
 h(0)y'(0)-h'(0)y(0)=0.
\end{equation}
Where $h(0)\ne0$ its Robin form is
\begin{equation}\label{eq:quadratic-robin}
 y'(0)=\frac{a u'_{s,r}(0)+b w'_{s,r}(0)}
              {a u_{s,r}(0)+b w_{s,r}(0)}\,y(0).
\end{equation}
The scaling $\Phi(z,s,t)=s^{1/2+z}\Phi(z,1,t/s)$ in
\eqref{eq:Phi} gives $\Phi(z,s,0)=s^{1/2+z}U_z$ and
$s\partial_t\Phi(z,s,0)=s^{1/2+z}U'_z$, where
$U_z=\Phi(z,1,0)$ and $U'_z=\partial_t\Phi(z,1,0)$.
With $q=-1/4-\sigma^2$, $\sigma>0$, put
$\ell_{a,b}(w)=a\operatorname{Re}w+(b/\sigma)\operatorname{Im}w$.
The scaled Robin law is
\begin{equation}\label{eq:scaled-quadratic-robin}
 sy'(0)=
 \frac{\ell_{a,b}(s^{i\sigma}U'_{i\sigma})}
      {\ell_{a,b}(s^{i\sigma}U_{i\sigma})}\,y(0).
\end{equation}
It is periodic in $\log s$ with period $\pi/\sigma$, with Dirichlet
condition at its poles. At $q=-1/4$, differentiating the scaling law
at $z=0$ gives, for a fixed coefficient line $[a:b]$,
\[
 s y'(0)=
 \frac{(a+b\log s)U'_0+b\,\partial_z U'_0}
      {(a+b\log s)U_0+b\,\partial_z U_0}\,y(0)
\]
whenever the denominator is nonzero. Here $U_0>0$ by
\eqref{eq:U-Up-gamma}. For $b=0$ the coefficient is $U'_0/U_0$
for all $s$. For fixed $b\ne0$ it equals
$U'_0/U_0+O(1/|\log s|)$ as $s\downarrow0$.

\begin{lemma}\label{lem:bessel-inner-phase}
Let $q=-1/4-\sigma^2$ with $\sigma>0$. For $s>0$, let $\xi_{\mathrm{D}}(s)=[-w_{s,r}(0):u_{s,r}(0)]$ and $\xi_{\mathrm{N}}(s)=[-w'_{s,r}(0):u'_{s,r}(0)]$ denote the Lagrangian lines in the canonical Bessel frame $(u_{s,r},w_{s,r})$ corresponding to the regular Dirichlet condition $y(0)=0$ and Neumann condition $y'(0)=0$, respectively. In the symplectically rescaled coordinates $(\sqrt\sigma a,b/\sqrt\sigma)$, their continuous angular lifts $\alpha_{\mathrm{D}}(s)$ and $\alpha_{\mathrm{N}}(s)$ satisfy
\begin{equation}\label{eq:bessel-DN-angles}
 \alpha_{\mathrm{D}}(s)=\sigma\log s+\arg U_{i\sigma}+\frac\pi2,\qquad
 \alpha_{\mathrm{N}}(s)=\sigma\log s+\arg U'_{i\sigma}+\frac\pi2,
\end{equation}
where $U_z=\Phi(z,1,0^+)$ and $U'_z=\partial_t\Phi(z,1,0^+)$.
With the lifts agreeing at $\sigma=0$, their separation is independent
of $s$ and satisfies
\begin{equation}\label{eq:inner-gudermannian-shift}
 a(\sigma):=\alpha_{\mathrm N}(s)-\alpha_{\mathrm D}(s)
 =\arg(U'_{i\sigma}/U_{i\sigma})\in(0,\pi/2),
\end{equation}
with positive right derivative $a_1:=a'(0^+)>0$.
\end{lemma}
\begin{proof}
In the rescaled coordinates $(X,Y)=(\sqrt\sigma a,b/\sqrt\sigma)$ one has
$\sqrt\sigma\,\ell_{a,b}(w)=X\operatorname{Re}w+Y\operatorname{Im}w$.
By~\eqref{eq:scaled-quadratic-robin}, the Dirichlet condition
$h_{s,r,\xi}(0)=0$ is orthogonality of $(X,Y)$ to $s^{i\sigma}U_{i\sigma}$
in $\mathbb R^2$, while $h'_{s,r,\xi}(0)=0$ is orthogonality to
$s^{i\sigma}U'_{i\sigma}$; this yields~\eqref{eq:bessel-DN-angles}.
The common rotation $\sigma\log s$ cancels in the difference; the range
$0<a(\sigma)<\pi/2$ and positivity of $a_1=a'(0^+)$ are verified in
Appendix~\ref{app:scalar-phases}, \eqref{eq:euler-inner-shift-value}.
\end{proof}

\begin{proposition}[A fixed Dirichlet cutoff need not converge]
\label{prop:dirichlet-cycle}
For $q=-1/4-\sigma^2$, $\sigma>0$, let $D_s$, $s>0$, be the
quadratic Dirichlet regularization $-d^2/dt^2+q/(t^2+s^2)$ with $y(0)=y(1)=0$.
The family $D_s$ has no norm-resolvent limit as $s\downarrow0$, although its minimal
operators converge in the gap topology. Its norm-resolvent subsequential
limits are exactly all the separated singular extensions with $y(1)=0$.
\end{proposition}
\begin{proof}
By \eqref{eq:bessel-DN-angles}, the physical Dirichlet line makes one
turn in $\mathbb{RP}^1$ each time $\log s$ changes by $\pi/\sigma$.
Every line thus occurs along a sequence $s_k\downarrow0$.
Theorem~\ref{thm:boundary-selection} identifies its operator limit.
Conversely, compactness of $\mathbb{RP}^1$ supplies a convergent
subsequence of lines along any sequence of cutoffs. Distinct lines
give distinct singular domains, which proves both assertions.
\end{proof}
This periodic cycling of boundary phases as $s\downarrow0$ is the
operator-theoretic counterpart of the renormalization-group limit cycle
for supercritical inverse-square potentials \cite{BC03,BHWZ17}.

\subsection{A single crossing through the singular limit}
\begin{example}\label{ex:quadratic-boundary-selection}
Take $q=-5/4$ and normalize the two zero solutions at $t=1$ by
$(u_s(1),u'_s(1))=(1,1/2)$ and $(w_s(1),w'_s(1))=(0,1)$.
This is a continuous change of the Wronskian-one frame in
\eqref{eq:normalized-frame}, obtained by inverting its matrix at $1$.
In particular $u_0=\sqrt t\cos(\log t)$ and
$w_0=\sqrt t\sin(\log t)$. Use $y(1)=0$ and the left line generated
by $h_{s,\theta}=\cos\theta\,u_s+\sin\theta\,w_s$.
Theorem~\ref{thm:boundary-selection} continues each such realization
to the singular phase $h_{0,\theta}=\sqrt t\cos(\log t-\theta)$.

For $\pi/4\le\theta\le3\pi/4$, the right Dirichlet zero solutions
are multiples of $w_s$. Since $W(w_s,h_{s,\theta})=-\cos\theta$,
the unique crossing occurs at $\pi/2$ and its kernel is spanned by
$w_s$. Differentiating the left condition and using Green's identity
gives $\lambda'(\pi/2)\|w_s\|_2^2=1$. Thus the spectral flow
equals $+1$ for every $0\le s\le1$, including the singular problem.
For physical Dirichlet conditions at both ends, all these limiting
phases occur along subsequences by Proposition~\ref{prop:dirichlet-cycle}.
\end{example}

\section{Quantitative growth of the Morse index}\label{sec:counting}

The structure of the logarithmic count is already visible in
Example~\ref{ex:quadratic-boundary-selection}: zero is an eigenvalue
precisely when the left boundary line coincides with the Cauchy-data line
of homogeneous solutions selected by the right Dirichlet condition. For
general matrix operators with lower-order perturbations, changing the
boundary condition alters the negative index by at most $2n$, and
diagonalizing $Q$ while bracketing the perturbation between scalar
potentials reduces the problem to scalar comparison operators in each
eigendirection. For each supercritical eigenvalue ($q<-1/4$), the left
physical Dirichlet line rotates in Bessel coordinates at rate
$\sigma\log(1/s)$ while the Cauchy-data line selected by $y(1)=0$
converges as $s\downarrow0$, so Theorem~\ref{thm:sf}
converts this relative rotation directly into the leading count. The
critical comparison gives an $O(1)$ bound on the Morse index of each
scalar component with $q\ge-1/4$.

For a real symmetric $n\times n$ matrix $Q$ with eigenvalues
$q_1,\ldots,q_n$, define
\begin{equation}\label{eq:count-coefficient}
 \kappa(Q)=\frac1\pi\sum_{j=1}^n\sqrt{(-q_j-1/4)_+},
 \qquad x_+=\max\{x,0\}.
\end{equation}

\subsection{Half-interval index growth with Coulomb-type perturbations}

\begin{theorem}[Boundary-independent leading Morse-index growth]
\label{thm:counting}
Let $0<s\le1$, and let $H_s$ be any self-adjoint realization on
$L^2((0,1),\R^n)$ of
\begin{equation}\label{eq:count-model}
 -\frac{d^2}{dt^2}+\frac{Q}{\varphi_s(t)}+B_s(t),\qquad
 \varphi_s(t)=t^2+s^2,
\end{equation}
where $B_s$ is measurable, symmetric and satisfies
\begin{equation}\label{eq:coulomb-remainder}
 \|B_s(t)\|\le M\left(1+\frac1{t+s}\right)
 \quad\hbox{for almost every }t.
\end{equation}
In particular, uniformly bounded matrix perturbations are allowed. The endpoint conditions may be
coupled and may depend arbitrarily on $s$. Then
\begin{equation}\label{eq:count-asymptotic}
 \nm(H_s)=\kappa(Q)\log\frac1s+O(1),\qquad s\downarrow0.
\end{equation}
The remainder is bounded by a constant depending only on $Q,n,M$,
not on the endpoint conditions or the particular family $B_s$. Neither
the critical value $q_j=-1/4$ nor the limit-point range $q_j\ge3/4$
need be excluded.
The remainder is locally uniform in $Q$ on sets whose spectra stay
away from $-1/4$. The same estimate holds on complex $L^2$, with
complex multiplicities and arbitrary self-adjoint boundary conditions.
\end{theorem}
\begin{proof}
For fixed $s>0$ the potential is bounded. The form of any regular
self-adjoint boundary condition restricts to the Dirichlet form on
$H^1_0$, of codimension at most $2n$ in its form domain. Thus the
min--max principle gives, with $D_s$ denoting the Dirichlet realization,
\begin{equation}\label{eq:boundary-count-bound}
 0\le \nm(H_{s,\Lambda})-\nm(D_s)\le 2n.
\end{equation}
The exact comparison, including a Dirichlet zero eigenvalue, is
\begin{equation}\label{eq:boundary-count-bound-exact}
 \nm(H_{s,\Lambda})-\nm(H_{s,\widehat\Lambda})
  =\itriple(W_s,\Lambda,\Lambda_D)
   -\itriple(W_s,\widehat\Lambda,\Lambda_D),
\end{equation}
where $W_s$ is the Cauchy-data space of homogeneous solutions and
$\Lambda_D$ the Dirichlet space (see also \cite{BCLS26} for boundary-space
eigenvalue interlacing and index differences between self-adjoint
extensions). The regular form argument in
Appendix~\ref{app:regular-boundary-inertia} covers bounded measurable
potentials.

Put $m_s(t)=M(1+(t+s)^{-1})$. With Dirichlet boundary conditions,
the inequalities $-m_s I\le B_s\le m_s I$ and a constant orthogonal
diagonalization of $Q$ bound the required count between two sums of
scalar counts. For each eigenvalue $q$ these are the Dirichlet
realizations $D_s^\pm(q)$ of $-y''+(q/(t^2+s^2)\pm m_s)y$.
These inequalities hold at each $s$. Only the scalar comparison
potentials enter the spectral-flow calculation.

\smallskip
\noindent\emph{Convergence of the Cauchy-data line.}
First take $q=-1/4-\sigma^2$, $\sigma\ge0$.
Write $F_s=(u_s,w_s)$ for the row in \eqref{eq:normalized-frame}
and $Y_s=\binom{F_s}{F'_s}$, so $\det Y_s=1$. For each sign and
$s>0$, choose a pair of perturbed solutions with the same initial
values and derivatives at $t=0$ as $u_s,w_s$, respectively.
The left Dirichlet angle computed in
Lemma~\ref{lem:bessel-inner-phase} then applies to either comparison
problem when $\sigma>0$. It remains to find the line selected by
$y(1)=0$. Writing the perturbed frames as $Y_s Z_s^\pm$, variation
of constants gives
\begin{equation}\label{eq:count-transfer}
 (Z_s^\pm)'=\pm m_s\binom{-w_s}{u_s}F_s Z_s^\pm,
 \qquad Z_s^\pm(0)=I .
\end{equation}
At $s=0$ this is an integral equation with integrable coefficient.
Indeed, writing $\rho=(t^2+s^2)^{1/2}$, the bounded hypergeometric
factor in \eqref{eq:Phi} and its bounded parameter derivative give
$|F_s(t)|\le C\sqrt\rho(1+|\log\rho|)$, uniformly on compact
$\sigma$-intervals. For $w_s$ use the divided difference in
\eqref{eq:normalized-frame}, including its derivative at zero.
Hence $m_s|F_s|^2\le C(1+|\log t|^2)$ on $(0,1)$. Dominated
convergence in $L^1$ and the integral equation give joint uniform
convergence of $Z_s^\pm$ as the parameters vary, including $s=0$.
The coefficient has trace zero, so $\det Z_s^\pm=1$ and the
perturbed frames have Wronskian one.

In its coefficient plane, the solutions satisfying the right Dirichlet
condition form the line
\[
 C_s^\pm=\ker\bigl(F_s(1)Z_s^\pm(1)\bigr).
\]
The row remains nonzero at $s=0$. Hence $C_s^\pm$ has a continuous
angular lift $\theta^\pm(s)$ on the closed interval, and this lift is bounded.
At the regular left end, $Z_s^\pm(0)=I$ gives
$\xi_D(s)=[-w_s(0):u_s(0)]$ in both frames.
When $\xi_D(s)=C_s^\pm$, any nonzero vector on their common line
gives a zero solution satisfying both Dirichlet conditions.
Conversely, each zero eigenfunction gives such a vector.

\smallskip
\noindent\emph{Counting the crossings.}
Fix $0<s<s_0\le1$ and keep $y(1)=0$. To apply the spectral-flow
formula, first set both left traces to zero on the underlying symmetric
operator. Its adjoint leaves the two left traces free. Imposing
$y(0)=0$ on that adjoint recovers the physical Dirichlet operator.
On $[s,s_0]$ the underlying symmetric operators have a common
domain in $H^2$, are injective by
initial-value uniqueness, and have a one-dimensional adjoint kernel.
Their boundary quotient is the left coefficient plane with form
$\det$, so Theorem~\ref{thm:sf} gives
$\sf(D_p^\pm(q);p\in[s,s_0])=-\Mas(\xi_D(p),C_p^\pm;p\in[s,s_0])$.
Let $\alpha_D$ and $\theta^\pm$ be continuous angular lifts of $\xi_D$
and $C^\pm$ in the same symplectic coordinates. At an intersection
$\alpha_D=\theta^\pm\pmod\pi$, the Robbin--Salamon crossing form
evaluates on a unit vector $x=(\cos\alpha_D,\sin\alpha_D)$ to
$\det(x,\dot x_{C^\pm}-\dot x_{\xi_D})=(\theta^\pm-\alpha_D)'=-(\alpha_D-\theta^\pm)'$,
so by \eqref{eq:clm-convention} the negative Maslov index
$-\Mas(\xi_D,C^\pm)$ counts upward passages of $\alpha_D-\theta^\pm$
through multiples of $\pi$. Because the operators $D_p^\pm(q)$ are
uniformly bounded below on $[s,s_0]$, we also have
$\sf(D_p^\pm(q);p\in[s,s_0])=\nm(D_s^\pm(q))-\nm(D_{s_0}^\pm(q))$,
whence
\begin{equation}\label{eq:count-relative-angle}
 \nm(D_s^\pm(q))-\nm(D_{s_0}^\pm(q))
   =\frac{(\alpha_D-\theta^\pm)(s_0)
                 -(\alpha_D-\theta^\pm)(s)}{\pi}+O(1).
\end{equation}
Indeed, setting $d=(\alpha_D-\theta^\pm)/\pi$, rotating the second line
to the horizontal, and deforming the difference angle to the affine path
with the same endpoints gives the exact signed count
$\lfloor d(s_0)\rfloor-\lfloor d(s)\rfloor$ under \eqref{eq:sf-convention}.

For $\sigma>0$, working in $(\sqrt\sigma\,a,b/\sqrt\sigma)$ gives
$\alpha_D(s_0)-\alpha_D(s)=\sigma\log(s_0/s)$ by
\eqref{eq:bessel-DN-angles}, while the lift of $C_s^\pm$ remains
bounded down to $s=0$ (locally uniformly in $\sigma>0$ by joint
continuity on compact parameter rectangles).
At $\sigma=0$, working in $(a,b)$ and differentiating
$\Phi(z,s,0)=s^{1/2+z}U_z$ yields the Dirichlet line
$[-U_0\log s-\partial_z U_0:U_0]$, which lies in one affine chart if
$U_0\ne0$ and is constant otherwise, so its lift is bounded as well.
Substituting these bounds into \eqref{eq:count-relative-angle} (and
noting that $\nm(D_{s_0}^\pm(q))$ is bounded on compact sets of $q$)
gives $\nm(D_s^\pm(q))=\sigma\pi^{-1}\log(1/s)+O(1)$ in both cases.
Finally, for every $q\ge-1/4$, form ordering on $H^1_0$ implies
$0\le\nm(D_s^\pm(q))\le\nm(D_s^-(-1/4))=O(1)$. Summing over the
eigenvalues of $Q$ and applying \eqref{eq:boundary-count-bound} proves
\eqref{eq:count-asymptotic}; the same comparisons on complex form
domains count each scalar eigenspace with its complex dimension.
\end{proof}

\subsection{Two-ended systems and full-period gyroscopic Hill operators}

\begin{proposition}\label{prop:two-ends}
On the two-ended interval $(0,2)$, consider the differential expression
\begin{equation}\label{eq:two-ends-model}
 -\frac{d^2}{dt^2} + V_{s_-, s_+}(t) + B(t),
\end{equation}
where $V_{s_-, s_+}(t) = Q_-/(t^2+s_-^2)$ for $t\in(0,1)$ and $V_{s_-, s_+}(t) = Q_+/((2-t)^2+s_+^2)$ for $t\in(1,2)$, with $Q_\pm$ real symmetric $n\times n$ matrices and $B$ a measurable, uniformly bounded symmetric matrix potential. Match $y, y'$ across the regular interface $t=1$, and impose arbitrary self-adjoint boundary conditions at $t=0$ and $t=2$. Then, as $s_-, s_+ \downarrow 0$ independently,
\begin{equation}\label{eq:two-count}
 \nm = \kappa(Q_-)\log\frac1{s_-} + \kappa(Q_+)\log\frac1{s_+} + O(1).
\end{equation}
The remainder is bounded uniformly over all self-adjoint boundary conditions and bounded sets of perturbations $B$.
\end{proposition}
\begin{proof}
Replacing the outer boundary conditions by Dirichlet conditions changes the
negative index by at most $2n$ by \eqref{eq:boundary-count-bound}, and
imposing $y(1)=0$ restricts the form domain by at most $n$ further
dimensions. The resulting decoupled Dirichlet problems on $(0,1)$ and $(1,2)$
have additive negative indices, so applying Theorem~\ref{thm:counting} to the
left half and, after reflection, to the right half gives \eqref{eq:two-count}.
\end{proof}

\begin{theorem}[Morse-index growth for gyroscopic Hill operators]
\label{thm:gyroscopic-master}\label{cor:gyroscopic-count}
Let $K^T=-K$ and $L^T=L$ be constant real $n\times n$ matrices, and let $B_e$ be a measurable, uniformly bounded, symmetric matrix family on $(0,2\pi)$. For any self-adjoint realization of the gyroscopic Hill operator
\begin{equation}\label{eq:gyro-hill}
 \mathcal G_e = -\frac{d^2}{dt^2} - 2K\frac{d}{dt} + \frac{L}{1-e\cos t} + B_e(t), \qquad 0<t<2\pi, \quad 0\le e<1,
\end{equation}
under arbitrary self-adjoint boundary conditions in the quasi-derivative $y'+Ky$ (including periodic, antiperiodic, or Floquet conditions $y(2\pi)=\omega y(0), y'(2\pi)+Ky(2\pi)=\omega(y'(0)+Ky(0))$ with $|\omega|=1$), the Morse index satisfies
\begin{equation}\label{eq:gyro-count}
 \nm(\mathcal G_e) = \kappa(2L)\log\frac1{1-e} + O(1), \qquad e\uparrow1.
\end{equation}
The remainder is uniform across all self-adjoint boundary conditions,
all Floquet multipliers $|\omega|=1$, and bounded sets of perturbations
$B_e$. It is also uniform on compact sets of $(K,L)$ whose matrices
$L$ have spectra separated from $-1/8$.
In particular, for the periodic essential operator $\mathcal A_e = -y'' - 2Jy' + R(1+e\cos t)^{-1}y$ in celestial mechanics \cite[Section 2.1]{HO16}, with constant symmetric matrix $R$ and standard complex structure $J$, the formula is
\begin{equation}\label{eq:ere-count}
 \nm(\mathcal A_e)
 =\frac1\pi\sum_j\sqrt{(-2\lambda_j(R)-1/4)_+}
                  \log\frac1{1-e}+O(1).
\end{equation}
\end{theorem}
\begin{proof}
Set $O(t)=\exp(-Kt)$ and write $y=Oz$. Then $y'+Ky=Oz'$ and the
transformed potential is
\[
 K^2+O^TB_e O+\frac{O^T L O}{1-e\cos t}.
\]
The gauge transports every self-adjoint boundary condition, and by
\eqref{eq:boundary-count-bound} replacing it by Dirichlet conditions at $0,2\pi$
and cutting with $z(\pi)=0$ changes the Morse index by at most $3n$.

On either half let $x$ be the distance from the singular end $t_*$,
where $t_*=0$ or $2\pi$, put $L_*=O(t_*)^T L O(t_*)$, and set
$s_e^2=2(1-e)/e$. For $e$ near one,
\[
 \frac{O(t)^T L O(t)}{1-e\cos t}
   =\frac{2L_*}{x^2+s_e^2}+E_e(x),\qquad
 \|E_e(x)\|\le C\bigl(1+(x+s_e)^{-1}\bigr).
\]
Indeed, since $(O^T L O)'=O^T[K,L]O$, we have
$\|O(t)^T L O(t)-L_*\|\le\|[K,L]\|x$, while
$2(1-\cos x)=x^2+O(x^4)$ and $1-\cos x\ge c x^2$ on
$[0,\pi]$; thus the numerator variation contributes $O(x/(x^2+s_e^2))$,
the Taylor remainder contributes $O(1)$, and replacing $2L_*/e$ by
$2L_*$ contributes $O((1-e)/(x^2+s_e^2))=O(1)$.
Rescaling by $x=\pi u$ and multiplying the operator by $\pi^2$ preserves
the Morse index while transforming the principal potential on $(0,1)$ into
$2L_*/(u^2+(s_e/\pi)^2)$ with a remainder satisfying
\eqref{eq:coulomb-remainder} for the cutoff $s_e/\pi$.
Theorem~\ref{thm:counting} therefore yields
$\kappa(2L_*)\log(\pi/s_e)+O(1)$ on each half.
Because both matrices $L_*$ are orthogonally conjugate to $L$ and
$2\log(1/s_e)=\log(1/(1-e))+\log(e/2)$, summing the two halves proves
\eqref{eq:gyro-count} with the stated uniformity, and translating by
$\pi$ gives \eqref{eq:ere-count} for the periodic essential operator.
\end{proof}

\begin{corollary}\label{cor:hill-count}
Setting $K=0$ in Theorem~\ref{thm:gyroscopic-master}, for any self-adjoint realization of the symmetric Hill family
\begin{equation}\label{eq:periodic-hill}
 \mathcal H_e = -\frac{d^2}{dt^2} + \frac{L}{1-e\cos t} + B_e(t), \qquad 0<t<2\pi, \quad 0\le e<1,
\end{equation}
under arbitrary self-adjoint boundary conditions, one has
\begin{equation}\label{eq:hill-count}
 \nm(\mathcal H_e) = \kappa(2L)\log\frac1{1-e} + O(1), \qquad e\uparrow1.
\end{equation}
\end{corollary}

\begin{example}\label{ex:noncommuting}
For a concrete noncommuting example, take
$K=\left(\begin{smallmatrix}0&-1\\1&0\end{smallmatrix}\right)$,
$L=\operatorname{diag}(-1,-1/8)$ and $B_e=0$.
The eigenvalues of $2L$ are $-2$ (supercritical, $-(-2)-1/4 = 7/4$) and $-1/4$ (critical, $-(-1/4)-1/4 = 0$).
By Theorem~\ref{thm:gyroscopic-master},
\[
 \nm(\mathcal G_e)=\frac{\sqrt7}{2\pi}\log\frac1{1-e}+O(1).
\]
Only the eigenvalue $-2$ contributes to the logarithmic coefficient, even with the coupling $-2Ky'$.
\end{example}

If $\kappa(2R)>0$, \eqref{eq:ere-count} gives
$\nm(\mathcal A_e)\ge a\log(1/(1-e))$ for each $0<a<\kappa(2R)$
and $e$ sufficiently close to one. Thus a fixed essential matrix with
an eigenvalue below $-1/8$ excludes a sequence of relative equilibria
with bounded Morse index and $e\uparrow1$.

\begin{remark}\label{rem:duality}
In \eqref{eq:count-relative-angle}, the spectral flow is evaluated on
intervals $[s,s_0]$ of positive cutoff where the operators remain
semibounded; the singular limit $s=0$ enters only to guarantee that
the right Dirichlet Cauchy line $C_s^\pm$ has a continuous limit and
hence a bounded angular lift as $s\downarrow0$.
\end{remark}

\section{Elliptic relative equilibria and Euler degeneracy curves}
\label{sec:ere-applications}

For a planar central configuration the full-period estimate gives a
Floquet count in terms of the essential Hessian. For Euler configurations
we then pass to the half-period sectors. Their indices determine the
endpoints of the degeneracy curves, and their relative phases give
the finer locations within the ordering of \cite{ZL17}.

\subsection{Normalization and full-period Floquet counts}
\label{subsec:ere-brake}

Let $N\ge3$ point masses $m_1,\ldots,m_N>0$ move in the plane,
with gravitational constant normalized to one. Define
\[
 M=\operatorname{diag}(m_1I_2,\ldots,m_NI_2),\qquad
 U(x)=\sum_{i<j}\frac{m_im_j}{|x_i-x_j|},\qquad
 \mathcal X=\{x\in(\R^2)^N:\textstyle\sum_i m_ix_i=0\}.
\]
Here $I_k$ is the $k\times k$ identity matrix, $U$ is the positive
Newtonian potential and $\langle Mx,y\rangle$ is the mass inner
product. A collision-free normalized central configuration
$a\in\mathcal X$ satisfies $\langle Ma,a\rangle=1$ and
$\nabla U(a)+U(a)Ma=0$. Let
$\mathbb J_N=\operatorname{diag}(J_2,\ldots,J_2)$, where
$J_2=\left(\begin{smallmatrix}0&-1\\1&0\end{smallmatrix}\right)$.
Choose a $2N\times d$ matrix $P$ whose columns form a mass-orthonormal
basis of the complement of $a,\mathbb J_Na$ inside $\mathcal X$.
The essential dimension is $d=2N-4$. Choosing the basis in complex
pairs makes $J=P^TM\mathbb J_NP=\operatorname{diag}(J_2,\ldots,J_2)$
a fixed complex structure.
This paired ordering differs from the grouped standard convention
only by a fixed permutation; $R$ below is expressed in the same
paired basis.
In the normalization of \cite[Section 2.1]{HO16},
\begin{equation}\label{eq:ere-normalized-R}
 R=I_d+\frac{P^TD^2U(a)P}{U(a)},\qquad P^TMP=I_d.
\end{equation}
This is the constrained Hessian on the inertia ellipsoid, normalized
by $U(a)$ and restricted to the essential directions. The translation
and Kepler blocks have already been removed. In true anomaly
$t\in[0,2\pi]$, the essential operator is
\begin{equation}\label{eq:ere-family}
 A_{R,e,\omega}y=-y''-2Jy'+\frac{R}{1+e\cos t}y,
 \qquad J^T=-J,\quad J^2=-I_d,
\end{equation}
on the complex self-adjoint domain
\begin{equation}\label{eq:ere-floquet-domain}
 D_\omega=\{y\in H^2((0,2\pi),\C^d):
 y(2\pi)=\omega y(0),\quad y'(2\pi)=\omega y'(0)\},
 \qquad |\omega|=1.
\end{equation}
The second condition can equivalently be imposed on $y'+Jy$.
The corresponding Hamiltonian matrix, in coordinates $(y'+Jy,y)$,
is generated by
\[
 B(t)=\begin{pmatrix}I_d&-J\\J&I_d-R/(1+e\cos t)\end{pmatrix}.
\]
The fundamental solution is normalized by
$\gamma_{R,e}'=\mathbb J_{2d}B\gamma_{R,e}$ and
$\gamma_{R,e}(0)=I_{2d}$, with
$\mathbb J_{2d}=\left(\begin{smallmatrix}0&-I_d\\I_d&0\end{smallmatrix}\right)$.
Writing $M_{R,e}=\gamma_{R,e}(2\pi)$ for its monodromy, one has
\begin{equation}\label{eq:ere-index-identification}
 \dim\ker A_{R,e,\omega}=\dim_\C\ker(M_{R,e}-\omega I),\qquad
 n_-(A_{R,e,\omega})=i_\omega(\gamma_{R,e}).
\end{equation}
Here $i_\omega$ is the Maslov-type index at the multiplier $\omega$,
with the Morse normalization fixed by the second equality \cite[Lemma~5.3]{HO16}.
If the endpoint convention uses $y(0)=\omega y(2\pi)$, it exchanges $\omega$ and $\bar\omega$;
complex conjugation gives identical counts for this real family.
For $\omega=\pm1$, the complex multiplicities equal the real ones;
they are not doubled.

For periodic and antiperiodic boundary conditions, an orthogonal
involution $S$ satisfying $SJS=-J$ and $SR=RS$ gives an orthogonal
decomposition into two self-adjoint problems on $[0,\pi]$.
Appendix~\ref{app:reversible} gives the full construction and the
endpoint conditions in each symmetry sector. The full-period
counting theorem below does not require this symmetry.

\begin{theorem}[Uniform Floquet and iterated counts]
\label{thm:ere-uniform-floquet}
Define
\begin{equation}\label{eq:ere-C}
 C(R)=\frac1\pi\sum_{j=1}^d
       \sqrt{(-2\lambda_j(R)-1/4)_+},\qquad
 L(e)=\log\frac1{1-e}.
\end{equation}
For a fixed real symmetric $R$, there is a constant $C_0$ such that
\begin{equation}\label{eq:ere-uniform-count}
 \left|n_-(A_{R,e,\omega})-C(R)L(e)\right|\le C_0
\end{equation}
for all $|\omega|=1$ and all $e$ sufficiently close to $1$.
The constant can be chosen uniformly for $R$ in a compact set of
symmetric matrices whose spectra stay a positive distance from
$-1/8$. At a critical eigenvalue $-1/8$ the bound
\eqref{eq:ere-uniform-count} still holds for fixed $R$, though
uniformity across the threshold is not claimed.

Let $N_m(R,e)$ be the negative index on the $m$-fold period,
with periodic endpoint conditions. Then
\begin{equation}\label{eq:ere-iteration}
 N_m(R,e)=\sum_{\omega^m=1}n_-(A_{R,e,\omega}),\qquad
 |N_m(R,e)-mC(R)L(e)|\le mC_0.
\end{equation}
The mean negative index exists and satisfies
\begin{equation}\label{eq:ere-mean}
 \lim_{m\to\infty}\frac{N_m(R,e)}m
 =\frac1{2\pi}\int_0^{2\pi}n_-(A_{R,e,e^{i\theta}})\,d\theta
 =C(R)L(e)+O(1).
\end{equation}
\end{theorem}
\begin{proof}
Translation by $\pi$ on Floquet functions gives
\eqref{eq:ere-uniform-count} from Theorem~\ref{thm:gyroscopic-master}
with $K=J$, $L=R$ and $B_e=0$.

On the $m$-fold period, translation by $2\pi$ is unitary and its
$m$th power is the identity. Its eigenspaces, indexed by
$\omega^m=1$, are mutually orthogonal and reduce the differential
operator; scaling restriction to a single period by $\sqrt m$ yields a
unitary equivalence with \eqref{eq:ere-floquet-domain}.
This proves the first identity in \eqref{eq:ere-iteration}, and summing
\eqref{eq:ere-uniform-count} gives the bound. For fixed $R$ and $e$, the
Floquet count is bounded by endpoint comparison and locally constant
except at the finitely many phases that are eigenvalues of the finite
matrix $M_{R,e}$. Hence it is Riemann integrable with equally spaced
Riemann sums $N_m/m$, proving \eqref{eq:ere-mean}.
\end{proof}

For each fixed $e<1$, the $2\pi$-periodic orthogonal gauge
$O(\theta)=\exp(-J\theta)$ transforms $A_{R,e,\omega}$ into a $d\times d$
$\theta$-periodic Schr\"odinger operator while preserving every Floquet
multiplier $\omega$, so \cite[Corollary~3.4(i)]{JLS17} already yields the
uniform cross-multiplier bound
$|n_-(A_{R,e,\omega_1})-n_-(A_{R,e,\omega_2})|\le2d$.
The logarithmic leading term $C(R)L(e)$ and its bounded remainder as
$e\uparrow1$ come from the singular Bessel boundary analysis above.

For a direct comparison with \cite[Theorem~1.1, Eq.~(1.5)]{HO16}, write
$\lambda_{\min}(R)=-1/8-r_1$ with $r_1>0$, and set
\[
 \varepsilon=\frac12\min\left\{\frac{r_1}{2r_1+5},\frac18\right\},
 \qquad \widehat e=\frac{1-e^2}{2}.
\]
Hu and Ou prove that
\[
 i_1(\gamma_e)\ge
 \frac{2\sqrt{r_1}}{\pi}
 \log\frac{\varepsilon^2}{\sqrt{\widehat e}}-6,
 \qquad \widehat e<\varepsilon^3,
\]
where $i_1(\gamma_e)=n_-(A_{R,e,1})$ under \eqref{eq:ere-index-identification}.
Theorem~\ref{thm:ere-uniform-floquet} sharpens this single-eigenvalue lower bound to
an asymptotic equality summing the contributions from all supercritical
eigenvalues, with a bounded remainder uniform in $|\omega|=1$.

\begin{theorem}[Asymptotic Floquet spectral flow]\label{thm:global-sf}
Fix $R=R^T$ with $C(R)>0$, and let $|\omega|=1$.
For the full-period family \eqref{eq:ere-family}, and fixed $0\le e_0<1$,
\begin{equation}\label{eq:ere-net-sf}
 \sf(A_{R,\cdot,\omega};[e_0,e])=-C(R)L(e)+O(1),\qquad e\uparrow1.
\end{equation}
For $e$ in the near-collision range of \eqref{eq:ere-uniform-count},
$n_-(A_{R,e,\omega})\le B$ implies
\begin{equation}\label{eq:ere-index-exclusion}
 1-e\ge\exp\left(-\frac{B+C_0}{C(R)}\right).
\end{equation}
For $\omega=\pm1$, each reversible half-interval sector has Morse index
$C(R)L(e)/2+O(1)$ and spectral flow $-C(R)L(e)/2+O(1)$ on $[e_0,e]$.
\end{theorem}
\begin{proof}
On each compact interval $[e_0,e]\subset[0,1)$ the domain is fixed
and the operators have a uniform lower bound. Hence, with the stated
endpoint convention,
\[
 \sf(A_{R,\cdot,\omega};[e_0,e])
 =n_-(A_{R,e_0,\omega})-n_-(A_{R,e,\omega}).
\]
Apply Theorem~\ref{thm:ere-uniform-floquet}; its lower estimate also
gives \eqref{eq:ere-index-exclusion}. When $SR=RS$ and $SJS=-J$, the
reversible symmetry decomposition on $[0,\pi]$ splits
$n_-(A_{R,e,\omega}) = n_-(A_{R,e,\omega}^{\rm even}) + n_-(A_{R,e,\omega}^{\rm odd})$
for $\omega=\pm1$, and the one-ended version of the same counting argument
gives Morse index $C(R)L(e)/2+O(1)$, and hence spectral flow
$-C(R)L(e)/2+O(1)$ on $[e_0,e]$, for each half-interval sector.
\end{proof}

\begin{remark}
When $C(R)>0$, the physical path has no norm-resolvent-continuous
self-adjoint Fredholm extension to $e=1$. Indeed, a
continuous Fredholm family on a compact parameter interval has
bounded spectral flow on its initial subintervals, contradicting
\eqref{eq:ere-net-sf} if such a physical extension existed.
\end{remark}

\subsection{Collinear Euler configurations: endpoints and bulk distribution}
\label{subsec:euler-bulk}

For an Euler collinear configuration, use the mass parameter
$\delta\in(0,7)$ (with closure $[0,7]$) of \cite[Section 5.2]{HO16}. It is the parameter
called $\beta$ in \cite{ZL17}, and is different from the Lagrange
mass parameter used below. To define it directly from the masses,
put $\widehat m_i=m_i/(m_1+m_2+m_3)$ and order the bodies on the
line as $1,2,3$. The ratio
$x=|a_1-a_2|/|a_2-a_3|$ is the unique positive root of
\begin{multline}\label{eq:euler-quintic}
 (\widehat m_3+\widehat m_2)x^5
 +(3\widehat m_3+2\widehat m_2)x^4
 +(3\widehat m_3+\widehat m_2)x^3\\
 -(3\widehat m_1+\widehat m_2)x^2
 -(3\widehat m_1+2\widehat m_2)x
 -(\widehat m_1+\widehat m_2)=0.
\end{multline}
In terms of this ratio, \cite[equations (1.4) and (2.1)]{ZL17} gives
\begin{equation}\label{eq:euler-mass-parameter}
 \delta=\frac{\widehat m_1(3x^2+3x+1)
              +\widehat m_3x^2(x^2+3x+3)}
 {x^2+\widehat m_2((x+1)^2(x^2+1)-x^2)}.
\end{equation}
The reduced Hessian is
\begin{equation}\label{eq:euler-R}
 R_\delta=\operatorname{diag}(-\delta,2\delta+3).
\end{equation}
For strictly positive masses the Euler parameter ranges over $(0,7)$ (with
closure $[0,7]$); however, none of our analytical results use the upper
bound $7$, and we work throughout with the extended family
\eqref{eq:euler-R} for all $\delta>0$. Only the first component can
be supercritical, so Theorem~\ref{thm:ere-uniform-floquet} gives
\begin{equation}\label{eq:euler-index}
 i_\omega(\gamma_{\delta,e})
   =c(\delta)L(e)+O(1),\qquad
 c(\delta)=\frac{\sqrt{(2\delta-1/4)_+}}\pi.
\end{equation}
The remainder is uniform in $\omega$ and locally uniform in
$\delta>1/8$. For fixed $\delta\le1/8$ the essential Floquet indices
are bounded. Above $1/8$, the coefficient refines the logarithmic
lower bound in \cite[Theorems~1.1 and~5.7]{HO16}.

For example, equal masses give $x=1$, $\delta=7/5$ and
$R=\operatorname{diag}(-7/5,29/5)$. Therefore, for every unit
Floquet multiplier,
\[
 i_\omega(\gamma_{7/5,e})
   =\frac{\sqrt{255}}{10\pi}\log\frac1{1-e}+O(1),\qquad
 N_m(7/5,e)=\frac{m\sqrt{255}}{10\pi}\log\frac1{1-e}+O(m).
\]
Here $N_m(7/5,e)$ abbreviates $N_m(R_{7/5},e)$, and the constant
in $O(m)$ is independent of the positive integer $m$.

\begin{remark}[Homothetic collision ends and comparison with \cite{HOY21}]
For a physical Euler configuration $a$ and a homothetic collision
$q(t)=r(t)a$, the constrained Hessian
eigenvalue $\lambda=-U(a)\delta$ gives the spiral threshold $\delta>1/8$
in \cite[Theorem~1.1(a)]{HOY21}. Fix such a $\delta$ and a finite-energy
orbit, and let $h=|t-T|$ be the physical time to collision.
The radial energy integral gives $r^3=(9U(a)/2)h^2(1+O(h^{2/3}))$.
In mass-orthonormal essential coordinates the potential in physical time
is $U(a)(R_\delta-I)/r^3$ by \eqref{eq:ere-normalized-R}.
Put $u=h^{1/3}$ and $y=\sqrt3\,u\,w$. The kinetic integral becomes
$\int |w'+w/u|^2\,du$. Integrating the cross term and using the two
Dirichlet conditions gives $\int (|w'|^2+2|w|^2/u^2)\,du$.
The transformed potential, including this kinetic contribution, is
\[
 \frac{2I}{u^2}+\frac{9U(a)u^4}{r^3}(R_\delta-I)
 =\frac{2R_\delta}{u^2}+O(1).
\]
To count on $a_\varepsilon\le u\le h_0^{1/3}$ with
$a_\varepsilon=\varepsilon^{1/3}$, we bound the remainder between $\pm C I$.
The positive component has bounded index by comparison with $-d^2/du^2-C$
on the full interval. In the weak component, put $\sigma=\sqrt{2\delta-1/4}$:
because the singular frame \eqref{eq:euler-frame} is $O(\sqrt u)$, the
variation-of-constants equation \eqref{eq:count-transfer} with coefficient $\pm C$ gives
$Z(a_\varepsilon)=I+O(a_\varepsilon^2)$, so in the rescaled coefficient
plane the left Dirichlet angle is
$\sigma\log a_\varepsilon+\pi/2+O(a_\varepsilon^2)$ while the right
Dirichlet line is fixed. To apply \eqref{eq:count-relative-angle} on a
fixed interval, set $\ell=h_0^{1/3}-a_\varepsilon$ and
$u=a_\varepsilon+\ell x$. For either scalar comparison potential $V$, the
unitary change of variable followed by multiplication by $\ell^2$ yields
$-\partial_x^2+\ell^2 V(a_\varepsilon+\ell x)$ on $(0,1)$, preserving the
index and sending the two Dirichlet conditions to $x=0,1$. Since these
form continuous regular families on each compact subinterval of
$0<a_\varepsilon<h_0^{1/3}$ and the logarithmic interval in $u$ has length
$\frac13\log(h_0/\varepsilon)$, their relative angle gives
\begin{equation}\label{eq:euler-physical-end}
 n^-_{\rm end}(\varepsilon)
 =\frac{\sqrt{2\delta-1/4}}{3\pi}\log\frac{h_0}{\varepsilon}
 +O(1),\qquad \varepsilon\downarrow0.
\end{equation}
For this single collision end, substituting $\lambda/U(a)=-\delta$ into
Eq.~(7) of the preprint version of \cite{HOY21} produces
$\sqrt{\delta-1/8}/(3\sqrt{2}\,\pi)$, which is half the coefficient in
\eqref{eq:euler-physical-end}; the factor of two stems from the radial
normalization in their Lemma~2.1(a), where the energy law gives
$r^{3/2}/h\to 3\sqrt{U(a)/2}$ whereas the printed limit for $r/h^{2/3}$
gives $r^{3/2}/h\to 3\sqrt{2U(a)}$.
\end{remark}

\paragraph{Known Euler curve structure.}
In the $(\delta,e)$-parameter plane, the periodic ($\omega=1$) and
antiperiodic ($\omega=-1$) degeneracy curves are the loci where the
monodromy matrix acquires a Floquet multiplier $+1$ or $-1$,
equivalently where the corresponding boundary value problem has a
nontrivial kernel ($\ker A_{R_\delta,e,\pm1}\ne\{0\}$). We use the
global curve structure of Long--Zhou \cite[Theorems~1.1 and~1.3]{ZL17},
recorded here in the notation used below.

\begin{proposition}[Interlacing and nullity of Euler curves, Long--Zhou]
\label{prop:euler-interlacing}
For the extended Euler operator family with parameter $\delta>0$ and eccentricity $e\in[0,1)$, there exist continuous families of periodic degeneracy curves $\delta=g_j(e)$ and antiperiodic degeneracy curves $\delta=h_j^-(e), h_j^+(e)$, $j\ge1$, satisfying:
\begin{enumerate}[label=(\roman*)]
\item \emph{Strict non-intersection}: A periodic degeneracy curve and an antiperiodic degeneracy curve never intersect in the interior $e\in[0,1)$:
\begin{equation}\label{eq:euler-non-intersection}
 g_j(e) \ne h_k^\pm(e) \qquad \text{for all } j, k\ge1 \text{ and all } e\in[0,1).
\end{equation}
\item \emph{Global interlacing order}: The degeneracy curves satisfy
\begin{equation}\label{eq:euler-interlacing}
 0 < h_1^-(e) \le h_1^+(e) < g_1(e) < h_2^-(e) \le h_2^+(e) < g_2(e) < \cdots.
\end{equation}
The weak inequalities occur only within antiperiodic pairs.
\item \emph{Nullity bounds}: For each periodic curve, the nullity satisfies $\nu_1(g_j(e))=2$.
For the antiperiodic curves, $\nu_{-1}(h_j^\pm(e)) = 1$ when $h_j^-(e) \ne h_j^+(e)$, and $\nu_{-1} \le 2$ at coincidences $h_j^-(e) = h_j^+(e)$.
\end{enumerate}
\end{proposition}
\begin{proof}
These are the curve-structure results of \cite[Theorems~1.1 and~1.3]{ZL17},
with $g_j=\beta_{2j-1}(1,e)=\beta_{2j}(1,e)$ and
$h_j^-=\beta_{2j-1}(-1,e)$, $h_j^+=\beta_{2j}(-1,e)$.
The exceptional degeneracy at $\delta=0$ is excluded by this numbering.
The periodic index jumps by two; the antiperiodic index jumps by
one per label, or by two when the labels coincide, by the same
curve-structure theorem and its monotonicity result.
Completeness of this
numbering follows from \cite[definition (4.26), Lemma~4.2(iii),
equation (4.13), and Corollary~4.3]{ZL17}. In particular, each interval
$(g_{j-1},g_j)$ between consecutive positive periodic roots contains exactly
two antiperiodic labels, counted with multiplicity.
\end{proof}

For $0<e<1$, Hu--Wu--Zhu \cite[Theorem~1.1 and Remark~1.2]{HWZ26}
give $h_j^-(e)<h_j^+(e)$ for every $j\ge1$, with one-dimensional
full antiperiodic kernel at each root. They place the weak coordinate
second in \cite[Eq.~(2.3)]{HWZ26}. It is odd about $t=0$ on the
lower branch and even on the upper branch. Antiperiodicity reverses
these parities at $\pi$, so in our half-interval notation $h_j^-$
belongs to $DN$ and $h_j^+$ to $ND$. At $e=0$ each pair coincides.

\paragraph{The homogeneous equation at collision.}
By the time-reversal symmetry on $[0,2\pi]$ detailed in
Appendix~\ref{app:reversible}, the periodic ($\omega=1$) and
antiperiodic ($\omega=-1$) problems split orthogonally into four
self-adjoint sectors on $[0,\pi]$, labeled by
$\alpha\beta\in\{NN,DD,ND,DN\}$ according to the conditions at $t=0$
and $t=\pi$ (the periodic sectors are $NN,DD$ and the antiperiodic
sectors are $ND,DN$). Write $q_e=1+e\cos t$,
$M=2\delta+3$, and $\tau=\pi-t$. In the forward coordinate $t$, the
homogeneous equations are
\[
 -y_1''+2y_2'-\delta y_1/q_e=0,\qquad
 -y_2''-2y_1'+My_2/q_e=0.
\]
Introducing the strong quasi-derivative $p=y_2'+2y_1$ gives
$p'=My_2/q_e$; at a regular endpoint the physical conditions read
$y_1'=y_2=0$ for $N$ and $y_1=p=0$ for $D$. At collision ($e=1$), the
strong component $y_2$ is limit-point: imposing $y_2\in L^2(0,\pi)$
eliminates the growing branch at collision, and combining this with the
regular $N$ or $D$ boundary condition at $t=0$ eliminates the boundary term in
the energy identity so that the defect vanishes and the auxiliary
variable $u=p/M$ satisfies a decoupled scalar equation, from which the
full coupled solution is reconstructed explicitly.

\begin{lemma}\label{lem:euler-limit-kernel}
Let $e=1$, $\delta>0$, and prescribe $\alpha=N$ or $D$ at $t=0$.
The Euler solutions with $y_2\in L^2(0,\pi)$ form a line. They are
given by
\begin{equation}\label{eq:euler-exact-reconstruction}
 u''+\left(1+\frac\delta q\right)u=0,\qquad
 y_2=qu',\qquad
 2y_1=(3\delta+4+\cos t)u+\sin t\,u',\quad q=1+\cos t,
\end{equation}
where $u'(0)=0$ for $N$ and $u(0)=0$ for $D$. Every nonzero
solution has a nonzero pair of weak Bessel coefficients.
\end{lemma}
\begin{proof}
For an Euler solution put $u=p/M$ and $W=u''+(1+\delta/q)u$.
Then $y_2=qu'$ and $2y_1=Mu-(qu')'$. Substituting these expressions into
$-y_1''+2y_2'-\delta y_1/q=0$ and using the identities
$q''=1-q$ and $(q')^2=q(2-q)$ gives
\begin{equation}\label{eq:euler-defect-energy}
 -(q^3W')'+2\delta q^2W=0.
\end{equation}
The regular conditions imply $W'(0)=0$ for $N$ and $W(0)=0$ for $D$.
Near collision ($\tau=\pi-t\downarrow0$), $q=\tau^2/2+O(\tau^4)$, so
\eqref{eq:euler-defect-energy} has indicial equation
$\alpha(\alpha+5)=4\delta$, with roots $-5/2\pm\nu$, where
$\nu=\sqrt{4\delta+25/4}>5/2$. If the growing branch $C\tau^{-5/2-\nu}$
has $C\ne0$, the inhomogeneous equation $u''+(1+\delta/q)u=W$ forces the
power $\kappa=(-5/2-\nu)+2=-1/2-\nu$. Substitution gives
$\kappa(\kappa-1)+2\delta=6\delta+7+2\nu>0$, and hence
\[
 u\sim\frac{C}{6\delta+7+2\nu}\tau^{-1/2-\nu},\qquad
 y_2\asymp\tau^{1/2-\nu}.
\]
Both homogeneous exponents of $u$ have positive real part, so this
term contradicts $y_2\in L^2$. For the remaining branch,
$W=O(\tau^{-5/2+\nu})$ and $q^3\overline W W'\to0$.
Integrating \eqref{eq:euler-defect-energy} gives
$\int_0^\pi(q^3|W'|^2+2\delta q^2|W|^2)\,dt=0$, hence $W=0$.

Conversely, \eqref{eq:euler-exact-reconstruction} reconstructs a
solution from either regular scalar line. If
$\gamma=1/2\pm\sqrt{1/4-2\delta}$ are distinct, reconstruction
multiplies the coefficient of $\tau^\gamma$ by
$(3\delta+3-\gamma)/2\ne0$. At $\delta=1/8$, writing
$u=\sqrt\tau(A+B\log\tau)+o(\sqrt\tau)$ gives
\begin{equation}\label{eq:euler-critical-trace-map}
 y_1=\sqrt\tau(a+b\log\tau)+o(\sqrt\tau),\qquad
 a=\frac{23}{16}A-\frac12B,\quad b=\frac{23}{16}B.
\end{equation}
In both cases $y_2'$ and $y_2/\sqrt q$ belong to $L^2$ and
the weak coefficients determine the solution. In particular
$-y_1''-\delta y_1/q=-2y_2'\in L^2$ and $p=Mu\to0$ at collision.
\end{proof}

\begin{lemma}[The Euler Friedrichs spectral gap]\label{lem:euler-friedrichs-gap}
Let $e=1$ and let $F_\delta^\alpha$ be the Euler realization with
regular condition $\alpha=N,D$ and the Friedrichs condition at collision.
Then
\begin{equation}\label{eq:euler-friedrichs-obstruction}
 \ker F_\delta^\alpha=\{0\}\qquad (0<\delta\le1/8).
\end{equation}
At $\delta=1/8$ the logarithmic coefficient $b$ in
\eqref{eq:euler-critical-trace-map} is nonzero for each nonzero
solution of Lemma~\ref{lem:euler-limit-kernel}.
\end{lemma}
\begin{proof}
Under $x=t/2$, the scalar equation $u''+(1+\delta/(2\cos^2(t/2)))u=0$
becomes $-v''-2\delta\sec^2 x\,v=4v$ for $v(x)=u(2x)$. Reflecting $v$
evenly for $N$ and oddly for $D$, by \eqref{eq:euler-exact-reconstruction}
a nonzero element of $\ker F_\delta^\alpha$ would give an eigenfunction with eigenvalue $4$ of
\[
 H_\delta=-\frac{d^2}{dx^2}-2\delta\sec^2 x,
 \qquad -\pi/2<x<\pi/2,
\]
with Friedrichs conditions at both ends. Indeed, reconstruction
preserves the principal power below the threshold and the logarithm-free
condition at the threshold.

Number the eigenvalues from one. On the common smooth core,
$H_{1/8}\le H_\delta\le H_0-2\delta$ in the sense of forms.
Writing $h_\delta$ for the closed Friedrichs form, completion of the
core in the shifted form norms gives
$\D(h_0)\subseteq\D(h_\delta)\subseteq\D(h_{1/8})$.
Each inequality extends to the smaller of its two form domains.
The second Dirichlet eigenvalue of $H_0$ is $4$. For $H_{1/8}$,
direct substitution gives
the eigenfunction $\sqrt{\cos x}\,(1-3\sin^2 x)$ at $25/4$.
It is logarithm-free at both ends and has two interior zeros, so it
is the third eigenfunction by Sturm oscillation \cite{Zet05}.
The min--max principle gives
\[
 \lambda_2(H_\delta)\le 4-2\delta<4<25/4\le\lambda_3(H_\delta).
\]
Hence $\ker F_\delta^\alpha=\{0\}$ for both parities. At $\delta=1/8$, the
condition $b=0$ is precisely the Friedrichs condition.
\end{proof}

\begin{corollary}[Endpoints of fixed-label Euler degeneracy curves]
\label{cor:euler-first-pair}
As $e\uparrow1$, $h_1^-(e)\to0$. Every other fixed label tends to
$1/8$ from above: these are $g_j(e),h_j^+(e)$ for $j\ge1$ and
$h_j^-(e)$ for $j\ge2$. In particular, $h_1^+-h_1^-\to1/8$.
\end{corollary}
\begin{proof}
Write $c=1/8$. At fixed $e<1$ let $i_\omega(x)$ and $\nu_\omega(x)$
be the Morse index and nullity at mass $x$. The indices just to the
right of zero are $B_1=3$ and $B_{-1}=2$
\cite[Sections~3.1 and~4]{ZL17}. As the mass increases, each positive
root adds its nullity to the index. Thus the numbers of crossings,
counted with multiplicity, in $(0,x)$ and $(0,x]$ are respectively
\[
 i_\omega(x)-B_\omega,\qquad
 i_\omega(x)+\nu_\omega(x)-B_\omega .
\]

Fix $0<a<c$. By Lemma~\ref{lem:euler-friedrichs-gap}, zero stays
outside both Friedrichs spectra on $(0,c]$. Continuing their indices
from small positive mass, passing to $e<1$ through the restricted
realizations, and freeing the weak collision value in the terminal
Neumann sectors gives \eqref{eq:euler-critical-sector-count} in
Appendix~\ref{app:euler-critical-index}.
Thus for $e$ close enough to one, both multipliers have index three and
zero nullity at $a$ and at $c$. Since $B_1=3$ and $B_{-1}=2$, there is
no periodic crossing in $(0,c]$, while the single antiperiodic crossing
in $(0,c]$ already lies in $(0,a)$. Its label is $h_1^-$, so $h_1^-\to0$.
Every periodic root and every remaining antiperiodic root is strictly
above $c$.

For each fixed $d>c$, \eqref{eq:euler-index} gives
$i_1(d),i_{-1}(d)\to\infty$. The number of roots below $d$ therefore
eventually exceeds any fixed label. Taking $d\downarrow c$ gives all
the remaining limits.
\end{proof}

In either enumeration the index at the $j$th curve is $2j+O(1)$,
with an absolute bound. This gives the distribution of labels growing with $L(e)$.

\begin{theorem}[Bulk distribution of Euler degeneracy curves]
\label{thm:euler-density}
Let $1/8<a<b<\infty$ be fixed. As $e\uparrow1$,
\begin{align}
 \#\{j:a<g_j(e)<b\}
 &=\frac{L(e)}{2\pi}
       \left(\sqrt{2b-1/4}-\sqrt{2a-1/4}\right)+O(1),
       \label{eq:euler-periodic-density}\\
 \#\{(j,\epsilon):a<h_j^\epsilon(e)<b,\ \epsilon\in\{-,+\}\}
 &=\frac{L(e)}{\pi}
       \left(\sqrt{2b-1/4}-\sqrt{2a-1/4}\right)+O(1).
       \label{eq:euler-antiperiodic-density}
\end{align}
Changing either open endpoint to a closed one changes only the bounded
remainder.

If $I$ is a compact interval contained in
$(0,\infty)$, then uniformly for integers $j$ with
$j/L(e)\in I$,
\begin{equation}\label{eq:euler-quantiles}
 \begin{split}
 g_j(e)&=\frac18+2\pi^2\left(\frac j{L(e)}\right)^2
                          +O\bigl(L(e)^{-1}\bigr),\\
 h_j^\pm(e)&=\frac18+2\pi^2\left(\frac j{L(e)}\right)^2
                          +O\bigl(L(e)^{-1}\bigr).
 \end{split}
\end{equation}
In particular, every interval between successive labeled curves
whose endpoints remain in a fixed compact subinterval of
$(1/8,\infty)$ has length $O(L(e)^{-1})$.
\end{theorem}
\begin{proof}
At fixed $e$, Proposition~\ref{prop:euler-interlacing} converts the
index difference into a count of roots. Either endpoint kernel has
dimension at most two, so
\[
 i_1(\gamma_{b,e})-i_1(\gamma_{a,e})
 =2\#\{j:a<g_j(e)<b\}+O(1).
\]
For $i_{-1}$ the factor is one when each label is counted.
Substituting \eqref{eq:euler-index} proves the density formulas.

Write $I=[\eta_-,\eta_+]$, with $0<\eta_-\le\eta_+
<\infty$. Choose $a_0,b_0\in(1/8,\infty)$ so that
$c(a_0)<2\eta_-$ and $c(b_0)>2\eta_+$.
The known index jumps and \eqref{eq:euler-index} at $a_0,b_0$
confine all three curves with $j/L(e)\in I$ to $(a_0,b_0)$ for
$e$ sufficiently close to one. Local uniformity of the remainder
can now be applied at the varying curve parameter. At each such
curve its relevant index equals $2j+O(1)$, hence
\[
 \sqrt{2g_j(e)-1/4}=\frac{2\pi j}{L(e)}+O(L(e)^{-1}),
\]
and the same equality holds with $g_j$ replaced by $h_j^\pm$.
All quantities are uniformly bounded on this compact interval;
squaring and dividing by two proves \eqref{eq:euler-quantiles}.
The leading expression changes by $O(L(e)^{-1})$ when $j$
changes by one and by zero between the two labels with the same
$j$. Interlacing and the uniform errors prove the final assertion.
\end{proof}

On the physical mass range, the limiting upper endpoint for $j/L(e)$
in \eqref{eq:euler-quantiles} is $c(7)/2=\sqrt{55}/(4\pi)$.
On compact supercritical intervals the density can also be written as
\[
 \frac1{L(e)}\sum_j\delta_{g_j(e)}
 \ \longrightarrow\ \frac{dx}{2\pi\sqrt{2x-1/4}},
\]
where $\delta_x$ denotes the unit point mass at $x\in(1/8,\infty)$. The convergence
is in the sense of integration against continuous compactly supported
test functions; for the antiperiodic labels the density is twice
as large. This formulation follows by approximating test functions
by step functions and using \eqref{eq:euler-periodic-density}.

\begin{figure}[!tbp]
 \centering
 \includegraphics[width=\textwidth]{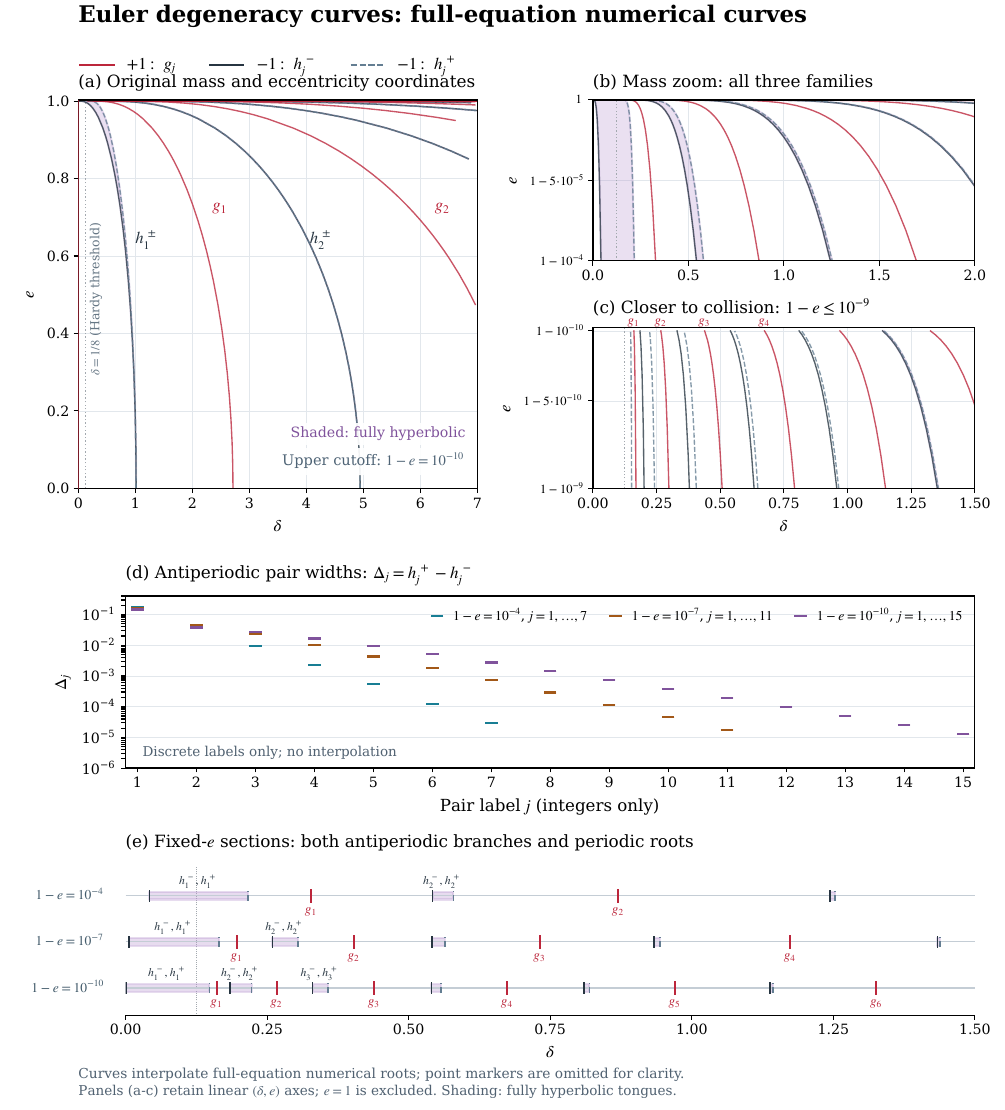}
 \caption{Numerical Euler degeneracy curves from the full coupled
 system, with point markers omitted. (a) The sampled range is
 $0\le\delta\le7$, $0\le e\le1-10^{-10}$; the dotted line marks
 $\delta=1/8$ and $\delta=0$ is also periodic-degenerate.
 (b,c) Successive enlargements in both mass and eccentricity,
 retaining linear $(\delta,e)$ axes and all three curve families.
 Uppermost plotted points correspond to the finite-$e$ cutoffs.
 Shading marks the open fully hyperbolic regions between antiperiodic
 pairs classified in \cite{ZL17}; unresolved strip widths are omitted.
 (d) Pair widths $\Delta_j=h_j^+-h_j^-$, including the first pair,
 on a logarithmic vertical scale; only integer labels are evaluated,
 with unconnected short ticks and the label ranges given in the legend.
 (e) Sections at $1-e=10^{-4},10^{-7},10^{-10}$: dark and dashed
 ticks mark $h_j^-,h_j^+$, red ticks mark $g_j$, and shaded segments
 show the pair widths in the original mass units.
 Curves connect computed roots for illustration only (the proofs do not
 use them). Independent integration and grid refinement agree within
 $10^{-9}$ on the sections in~(e) and within $10^{-8}$ elsewhere; at the
 two smallest sampled positive $e$, the second antiperiodic pair width is
 below tabulation accuracy.}\label{fig:euler-numerical-overview}
\end{figure}

Figure~\ref{fig:euler-numerical-overview} displays computed degeneracy
curves for the full coupled Euler system. The enlargements retain
both antiperiodic branches, and fixed-$e$ sections resolve their
pair widths alongside the periodic reference curves.

\subsection{Euler boundary phases and curve locations}
The Morse count fixes the limiting endpoints of all curves with fixed
label, whereas resolving the spacing between neighboring curves requires
the relative angle between the Cauchy-data line of homogeneous solutions
and the physical endpoint line at collision. We place the weak component
in the continuous Bessel coordinates of Theorem~\ref{thm:boundary-selection}
and control the strong component by its weighted energy form: uniformly on
compact mass intervals crossing $\delta=1/8$, this Cauchy-data line
converges as $e\uparrow1$ to a limit $C_1^\alpha(\delta)$ depending on the
mass $\delta$ and regular left condition $\alpha\in\{N,D\}$, independent
of the terminal strong condition. In the supercritical range, two sectors
with the same left condition ($NN$ versus $ND$, or $DD$ versus $DN$)
share this limiting outer phase while their physical weak lines differ at
finite $e<1$ by the inner shift $a(\rho)$ of
Lemma~\ref{lem:bessel-inner-phase}. For each fixed $\delta>1/8$,
$\rho=\sqrt{2\delta/e-1/4}$ tends to $\sigma=\sqrt{2\delta-1/4}$ as
$e\uparrow1$, so $a(\rho)\to a(\sigma)$. Estimating their finite-$e$
outer-phase difference then isolates the gaps $g_j-h_j^+$ and
$h_{j+1}^--g_j$.

To vary the weak boundary line at $\pi$, keep $y_2(\pi)=0$ for terminal
$N$ and $p(\pi)=0$ for terminal $D$. These choices recover the physical
sectors when the weak condition is respectively Neumann or Dirichlet.

\begin{lemma}\label{lem:euler-full-family}
Fix the regular condition $\alpha$ and the terminal strong condition
$\beta$ above. Prescribing a weak coefficient line
$\xi\in\mathbb{RP}^1$ gives full Euler realizations
$A_{\delta,e,\xi}^{\alpha\beta}$ in $L^2((0,\pi);\mathbb C^2)$
with compact resolvent. For $e$ near one they are jointly
norm-resolvent continuous in $(\delta,e,\xi)$, uniformly on compact
mass intervals in $(0,\infty)$ and all coefficient lines. At $e=1$
the realization is independent of $\beta$.
\end{lemma}
\begin{proof}
Take the weak component $f=y_1$ in the domain of
$T=-\partial_t^2-\delta/q_e$ with the weak condition $\alpha$ at
$0$ and the line $\xi$ at $\tau=0$. Set
\[
 s^2=2(1-e)/e,\qquad r|r|=1/4-2\delta/e,\qquad
 d(\tau)=2(1-\cos\tau).
\]
Then $-\delta/q_e=q(r)/(s^2+\tau^2)+R$, with
\[
 q(r)=r|r|-1/4,\qquad
 R=-\frac{2\delta}{e}\,
 \frac{\tau^2-d(\tau)}{(s^2+d(\tau))(s^2+\tau^2)}.
\]
Since $d(\tau)\asymp\tau^2$ and $\tau^2-d(\tau)=O(\tau^4)$,
$R$ is a locally uniformly bounded, strongly continuous multiplier.
Cutoff functions away from collision give right inverses for the two
regular traces. The proof of Theorem~\ref{thm:boundary-selection}
therefore applies on $(0,\pi)$ with either fixed regular condition.
Thus $T$ has compact resolvent and varies in norm resolvent, including
at $r=0$.

For the strong component $v=y_2$ set
\[
 t(v,w)=\int_0^\pi (v'\overline{w'}+Mv\bar w/q_e)\,dt .
\]
Its form domain $V$ consists of $H^1$ functions of finite $t$-energy,
with zero value at each regular $N$ endpoint. At $e=1$ finite energy
forces $v(\pi)=0$. Define the coupled operator on pairs
$(f,v)\in\D(T)\times V$ for which there is a $g\in L^2$ such that
\[
 t(v,w)+2\int_0^\pi f\overline{w'}\,dt=\int_0^\pi g\bar w\,dt
 \qquad(w\in V).
\]
Its value is $(Tf+2v',g)$. This identity gives the natural condition
$p=0$ at a regular $D$ endpoint and makes the operator symmetric.
Smooth pairs of compact support give a dense domain.

For $e<1$ these are regular self-adjoint boundary conditions. For
$z=\pm i$ and data $H=(H_1,H_2)$ their resolvent satisfies
$\|(f,v)\|\le\|H\|$. Taking $w=v$ in the strong equation and
using $g=H_2+zv$ gives, on compact mass intervals,
\[
 \|v'\|_2+\|v/\sqrt{q_e}\|_2\le C\|H\|,\qquad
 f=(T-z)^{-1}(H_1-2v').
\]
Suppose $(\delta_n,e_n,\xi_n)\to(\delta_*,1,\xi_*)$ and
$H_n\rightharpoonup H$. After taking a subsequence, $v_n$ converges
weakly in $H^1$ and strongly in $L^2$. Set $R_n=(T_n-z)^{-1}$ and
$k_n=H_{1,n}-2v_n'$, so that $k_n\rightharpoonup k=H_1-2v'$.
With $R_*=(T_*-z)^{-1}$ we have
\[
 \|f_n-R_*k\|_2\le\|R_n-R_*\|\,\|k_n\|_2+\|R_*(k_n-k)\|_2
 \longrightarrow0 .
\]
The first term tends to zero by norm convergence of the scalar
resolvents and the second by compactness of $R_*$, while Fatou's lemma
gives $v/\sqrt{q_1}\in L^2$. Passing to the limit in the strong equation
first for tests supported away from collision extends it to all of $V_*$:
for a smooth cutoff vanishing on $\tau<h$ and equal to $1$ on $\tau\ge2h$,
the squared $V_*$-norm of the cutoff error is bounded by a constant times
$E_h(w)=\int_0^{2h}(|w'|^2+|w|^2/q_1)\,d\tau$, so continuity of the weak
formulation gives an error of order $E_h(w)^{1/2}\to0$.
Applying this argument to fixed $H$ proves surjectivity of the
symmetric limit operator at $z=\pm i$ and hence self-adjointness, with
compact resolvent by the same energy bounds. Uniqueness of the limit
identifies every subsequential limit (including along sequences at
$e=1$). If norm-resolvent convergence failed, choosing unit data $H_n$
along which the resolvent difference stayed bounded away from zero and
passing to a weakly convergent subsequence would yield a contradiction,
since both the varying resolvents and the compact limit resolvent send
$H_n$ to the same strong limit; at regular parameter values $e<1$,
continuity follows directly from continuous regular boundary data.
\end{proof}

\begin{lemma}\label{lem:euler-uniform-phases}
Fix a compact mass interval $J\subset(0,\infty)$. For $e$ sufficiently
close to one, the homogeneous solutions satisfying the regular condition
$\alpha$ and the terminal strong condition $\beta$ have a single weak
coefficient line $C_e^{\alpha\beta}(\delta)$. These lines converge
uniformly on $J$ to the line $C_1^\alpha(\delta)$ of
Lemma~\ref{lem:euler-limit-kernel}. On each compact interval
$K\subset(1/8,\infty)$ there are continuous real phases
$F_{\alpha\beta}$ such that the physical sector has a kernel exactly
when $F_{\alpha\beta}\in\pi\mathbb Z$, and
\begin{equation}\label{eq:four-quantization-euler}
 F_{\alpha\beta}(e,\delta)
 =\frac{L(e)\sigma(\delta)}2+\Theta_\alpha(\delta)
       -\arg V_\beta(i\sigma(\delta))+r_{\alpha\beta}(e,\delta),
 \qquad \max_{\alpha,\beta}\sup_K|r_{\alpha\beta}|\longrightarrow0,
\end{equation}
where $\sigma(\delta)=\sqrt{2\delta-1/4}$, $V_N=U'$, $V_D=U$ as in
Lemma~\ref{lem:bessel-inner-phase}, and $\Theta_\alpha$ is continuous.
\end{lemma}
\begin{proof}
The singular solution reconstructed in Lemma~\ref{lem:euler-limit-kernel}
belongs to the realization with its own weak coefficient line, since
its first component lies in the scalar maximal domain and integrating
by parts against the strong form core yields the second equation.
At $(1,\delta_*)$, choose any line $\eta\ne C_1^\alpha(\delta_*)$; by
Lemmas~\ref{lem:euler-limit-kernel} and~\ref{lem:euler-full-family},
the realization with line $\eta$ is invertible at $(1,\delta_*)$ and
hence remains invertible in a neighborhood. For $e<1$ near one, the two
regular conditions at $t=0$ cut the four-dimensional solution space down
to dimension two, and the single strong condition at $t=\pi$ leaves a
subspace of dimension at least one. Because the realization with line
$\eta$ is invertible, the weak trace map on this subspace is injective
and its image meets $\eta$ trivially, so the image is a unique line
$C_e^{\alpha\beta}(\delta)\in\mathbb{RP}^1$. Along any sequence
tending to $(1,\delta_*)$, every subsequential limit of these lines in
$\mathbb{RP}^1$ produces a limiting realization with a zero eigenvalue
(as invertibility is open in the norm-resolvent topology), which
Lemma~\ref{lem:euler-limit-kernel} identifies as $C_1^\alpha(\delta_*)$.
By compactness of $J$, convergence is uniform in all four sectors.

Set $\rho^2=2\delta/e-1/4$ where this is positive, and let
$\theta_e^{\alpha\beta}$ be a continuous angular lift of
$C_e^{\alpha\beta}$ in the rescaled coordinates
$(\sqrt\rho\,a,b/\sqrt\rho)$. With $\alpha_\beta(s)=\rho\log s+\arg V_\beta(i\rho)+\pi/2$
denoting the physical boundary angle of Lemma~\ref{lem:bessel-inner-phase},
the angle of the Cauchy line relative to the physical line is
\[
 F_{\alpha\beta}=\theta_e^{\alpha\beta}-\alpha_\beta(s)
 =\rho\log(1/s)+\theta_e^{\alpha\beta}-\arg V_\beta(i\rho)-\frac\pi2,
\]
so the two lines coincide precisely when $F_{\alpha\beta}\in\pi\mathbb Z$.
On $K$, the coordinate rescaling is uniformly invertible and compatible
lifts converge uniformly to $\theta_1^\alpha$. Substituting
$\rho-\sigma=O_K(1-e)$ and $\log(1/s)=L/2-\tfrac12\log(2/e)$ yields
\eqref{eq:four-quantization-euler} with
$\Theta_\alpha=\theta_1^\alpha-\tfrac\sigma2\log2-\tfrac\pi2$.
\end{proof}

For $e<1$ and $\rho>0$, the weak Dirichlet and Neumann lines share the
rotation $\rho\log s$ and differ by $a(\rho)$ by
\eqref{eq:inner-gudermannian-shift}, which tends to $a(\sigma)$ as
$e\uparrow1$ for each fixed $\delta>1/8$. Since
$a(\sigma)=a_1\sigma+o(\sigma)$ with $a_1=a'(0^+)>0$ near $\sigma=0$, an
$o(1)$ phase error could exceed the shift itself; distinguishing the
neighboring Euler roots therefore requires an $o(\sigma)$ error across
the range where those roots lie.

\begin{remark}[The repeated-root coefficient chart]\label{rem:euler-critical-chart}
By Lemma~\ref{lem:euler-friedrichs-gap}, $b\ne0$ on each limiting
Cauchy line at $\delta=1/8$. Thus $h_e=a_e/b_e$ is a valid coordinate
near the critical point, and Lemma~\ref{lem:euler-uniform-phases}
gives $h_e\to h_1$ uniformly there. In these coordinates
$\theta_e=\pi/2-\arctan(\rho h_e)$, so uniform convergence of $h_e$
alone (without derivative bounds on the Cauchy lines) gives
\[
 \bigl|\arctan(\rho h_e)-\arctan(\rho h_1)\bigr|
 \le \rho\,|h_e-h_1|=o(\rho).
\]
The same chart explains why the limiting $NN$ slow phase
$B(\sigma):=\Theta_N(\delta)-\arg U'_{i\sigma}$ (normalized by
$\theta_1^N\to\pi/2$ and $\arg U'_0=0$, so $B(0)=0$) admits a
first-order expansion at $\sigma=0$. By the critical Friedrichs gap and
continuity of the limiting Cauchy line
(Lemmas~\ref{lem:euler-friedrichs-gap} and~\ref{lem:euler-uniform-phases}),
the coordinate $h_1^N$ tends to a finite limit $H_N$ as
$\delta\downarrow1/8$, giving $\theta_1^N=\pi/2-H_N\sigma+o(\sigma)$.
Together with regularity of the inner phase $\arg U'_{i\sigma}$ at
$\sigma=0$ (from \eqref{eq:U-Up-gamma}), this yields
$B(\sigma)=b_0\sigma+o(\sigma)$ with $b_0=B'(0^+)$, evaluated in closed
form in Appendix~\ref{app:scalar-phases}, \eqref{eq:euler-b0}.

The largest error in replacing $\rho$ by $\sigma>0$ comes from
the term $L\rho/2$. On a bounded mass interval,
\[
 \frac{L|\rho-\sigma|}{\sigma}
 =\frac{2\delta(1-e)L}{e\,\sigma(\rho+\sigma)}
 \le C\,\frac{(1-e)L}{\sigma^2},
\]
so $(1-e)L/\sigma^2\to0$ preserves the $o(\sigma)$ accuracy, as verified
below on the full range of the localization theorem.
\end{remark}

\begin{proposition}[Determination of the periodic phase integers]
\label{prop:euler-calibration}
Fix $1/8<d<\infty$ and let $e<1$ be sufficiently close to one.
Normalize the Cauchy angles by continuation from $\theta=\pi/2$
at $(e,\delta)=(1,1/8)$ and take $\arg U_0=\arg U'_0=0$.
Then the exact phases of Lemma~\ref{lem:euler-uniform-phases} satisfy
\begin{equation}\label{eq:euler-calibrated-integers}
 0<F_{NN}(e,1/8),F_{DD}(e,1/8)<\pi,\qquad
 F_{NN}(e,g_j)=F_{DD}(e,g_j)=j\pi
\end{equation}
for every periodic root $g_j(e)\le d$.
\end{proposition}
\begin{proof}
In the critical chart $h_e=a_e/b_e$ is bounded and
$\theta_e=\pi/2-\arctan(\rho h_e)$. At $\delta=1/8$,
$\rho=s/(2\sqrt2)$, where $s^2=2(1-e)/e$. Thus
\[
 F_{\alpha\alpha}(e,1/8)
 =\rho\log(1/s)-\arctan(\rho h_e)-\arg V_\alpha(i\rho)
 =\rho\bigl(\log(1/s)+O(1)\bigr)\in(0,\pi).
\]

Fix $e<1$ and isolate a periodic crossing $g_j$ between two invertible
masses. Imposing the regular left and strong right conditions while
setting both weak right traces to zero defines a symmetric operator with
fixed domain as $\delta$ varies, Fredholm index $-1$ by the regular
boundary count, and trivial kernel by the proof of
Lemma~\ref{lem:euler-uniform-phases}. Its adjoint frees the two weak
traces, and Green's form on their plane reduces to $W_\tau$ because the
strong boundary contribution vanishes under either $y_2=0$ or $p=0$.
Since the Bessel frame has Wronskian one and the rescaling of $(a,b)$
has determinant one, the relation
$F_{\alpha\alpha}=\theta_e^{\alpha\alpha}-\alpha_\alpha(s)$ shows that
the Robbin--Salamon crossing form \eqref{eq:clm-convention} for the
ordered pair (physical line, Cauchy line) evaluates on a unit crossing
vector to $F'_{\alpha\alpha}$; its local Maslov index therefore counts
upward passages of $F_{\alpha\alpha}$ through multiples of $\pi$.

By Proposition~\ref{prop:euler-interlacing}, the total periodic nullity
at $g_j$ is two and the periodic index increases by two. Because each
sector has nullity at most one, each sector index increases by one.
At fixed $e$, the physical quadratic forms share a common lower bound on
the mass interval, so their local spectral flows are $-1$ and
Theorem~\ref{thm:sf} yields Maslov index $+1$ in each sector. Since
$g_1$ is the first periodic root above $1/8$ by
Corollary~\ref{cor:euler-first-pair} and both phases start in
$(0,\pi)$, they pass upwards through $\pi,2\pi,\ldots$ at the successive
roots.
\end{proof}

\begin{theorem}[Uniform localization from the Hardy threshold to the bulk]
\label{thm:euler-fine-splitting}\label{thm:euler-uniform-localization}
Put $L=\log(1/(1-e))$ and $\sigma(\delta)=\sqrt{2\delta-1/4}$.
Fix $1/8<d<\infty$ and set $s_j=\sigma(g_j(e))$, $a_j=a(s_j)$,
with $a$ as in \eqref{eq:inner-gudermannian-shift}. The limiting $NN$ slow phase
$B$, normalized by $B(0)=0$ as in Remark~\ref{rem:euler-critical-chart}, satisfies
\begin{equation}\label{eq:euler-absolute-position}
 \frac L2s_j+B(s_j)=j\pi+o(s_j)
\end{equation}
uniformly for $j\ge1$ with $g_j(e)\le d$. Uniformly for $j\ge2$ in
this range,
\begin{align}
 s_j-s_{j-1}&=\frac{2\pi}{L}+o(L^{-1}),
       \label{eq:euler-uniform-spacing}\\
 \sigma(h_j^+)&=s_j-\frac{2a_j}{L}+o(s_j/L),\qquad
 \sigma(h_j^-)=s_j-\frac{2(\pi-a_j)}L+o(L^{-1}).
       \label{eq:euler-phase-splitting}
\end{align}
Their mass gap is
\begin{equation}\label{eq:euler-uniform-large-gap}
 h_j^+-h_j^-=
 \frac2L\left(s_j-\frac\pi L\right)(\pi-2a_j)
 +o\!\left(\frac{s_j+L^{-1}}L\right).
\end{equation}
If also $g_{j+1}(e)\le d$, the same periodic root gives
\begin{equation}\label{eq:euler-uniform-small-gap}
 \sigma(h_{j+1}^-)=s_j+\frac{2a_j}L+o(s_j/L),\qquad
 h_{j+1}^--h_j^+=\frac{4s_ja_j}L+o(s_j^2/L).
\end{equation}
\end{theorem}
\begin{proof}
\noindent\emph{Step 1: critical-to-bulk phase control.}
In the critical chart supplied by Lemma~\ref{lem:euler-friedrichs-gap}
and Remark~\ref{rem:euler-critical-chart}, the exact phase is
\[
 F_{\alpha\beta}=\frac{L\sigma}{2}+H_{\alpha\beta},\qquad
 H_{\alpha\beta}=\frac{L(\rho-\sigma)}2
 -\frac\rho2\log(2/e)-\arctan(\rho h_e^{\alpha\beta})
 -\arg V_\beta(i\rho),
\]
where $h_e^{\alpha\beta}\to h_1^\alpha$ uniformly.
Corollary~\ref{cor:euler-first-pair} gives $g_1\to1/8$. The
normalized phase at this root is
$\pi=\rho(L/2+O(1))$. Hence $L\rho\to2\pi$ there.
Since $L^2(\rho^2-s_1^2)=O((1-e)L^2)\to0$, we get
$s_1=2\pi/L+o(L^{-1})$.
On all subsequent periodic intervals $\sigma\ge c/L$, so
$(1-e)L/\sigma^2\le C(1-e)L^3\to0$.
The conversion to $\sigma$ thus contributes $o(\sigma)$, uniformly.
Combining the critical chart with uniform phase convergence away
from zero, for each fixed $c>0$ we obtain
\[
 \sup_{\delta:\ c/L\le\sigma(\delta)\le\sigma(d)}
 \frac{\left|F_{NN}(e,\delta)-L\sigma(\delta)/2
                    -B(\sigma(\delta))\right|}{\sigma(\delta)}
 \longrightarrow0\qquad(e\uparrow1).
\]
This proves \eqref{eq:euler-absolute-position}.

\smallskip
\noindent\emph{Step 2: outer-phase comparison and root bracketing.}
The slow phases are uniformly bounded. Subtracting the $NN$
phases at consecutive roots first gives $s_j-s_{j-1}=O(L^{-1})$.
Uniform continuity of the limiting slow phase improves this to
\eqref{eq:euler-uniform-spacing}. At $g_j$ compare $ND$ with $NN$,
which have the same limiting Cauchy line $C_1^N$, and $DN$ with $DD$,
which have limiting line $C_1^D$. Within each pair the outer angles
differ by $o(s_j)$, by the critical coefficient chart near zero and
by uniform phase convergence elsewhere, while the inner boundary shift is
$a(\rho)=a_j+o(s_j)$ by Lemma~\ref{lem:bessel-inner-phase} and the
$\rho$-to-$\sigma$ estimate above.
Since both periodic phases equal $j\pi$, we obtain
\[
 F_{ND}(e,g_j)-F_{NN}(e,g_j)=a_j+o(s_j),\qquad
 F_{DN}(e,g_j)-F_{DD}(e,g_j)=-a_j+o(s_j).
\]

For $|\Delta\sigma|\le C\sigma/L$, the slow-phase increment is
$o(\sigma)$. Indeed, uniform convergence to continuous $h_1^\alpha$
bounds the change of $\arctan(\rho h_e^{\alpha\beta})$ by
$O(|\Delta\sigma|)+\sigma o(1)$; the inner arguments are smooth and
the conversion term is already $o(\sigma)$. Away from the critical chart
ordinary uniform continuity suffices. Apply these estimates up to a
fixed mass slightly larger than $d$ to allow for the endpoints of the
test intervals. For each sufficiently small fixed $\varepsilon>0$, shifting $\sigma$
from $s_j$ by $\Delta\sigma=-2a_j/L\pm\varepsilon s_j/L$ changes the
fast phase by $\frac{L}{2}\Delta\sigma=-a_j\pm\frac12\varepsilon s_j$
while altering the slow phase by $o(s_j)$. Combining this with
$F_{ND}(e,g_j)-j\pi=a_j+o(s_j)$ cancels $a_j$ and leaves
$F_{ND}-j\pi=\pm\frac12\varepsilon s_j+o(s_j)$ uniformly in the
admitted labels, so the opposite signs bracket a root. For $DN$, the
same test about $s_{j-1}+2a_{j-1}/L$ brackets a root at level
$(j-1)\pi$ with half-width $\varepsilon s_{j-1}/L$. On
$0<s\le\sigma(d)$, $a(s)/s$ has a positive lower bound and
$a(s)<\pi/2$, which together with \eqref{eq:euler-uniform-spacing}
places both test intervals strictly between $g_{j-1}$ and $g_j$ for
sufficiently small $\varepsilon$. The sector labels of
\cite[Remark~1.2]{HWZ26} identify these roots as $h_j^+$ and $h_j^-$,
and letting $\varepsilon\downarrow0$ yields the uniform errors
$o(s_j/L)$ and $o(s_{j-1}/L)$. Using \eqref{eq:euler-uniform-spacing}
and continuity of $a$ to express $h_j^-$ in terms of $s_j$ gives the second
formula in \eqref{eq:euler-phase-splitting}, while the corresponding test
interval above $g_j$ yields \eqref{eq:euler-uniform-small-gap} whenever
$g_{j+1}\le d$.

\smallskip
\noindent\emph{Step 3: conversion to mass gaps.}
Finally, $\delta=1/8+\sigma^2/2$ exactly. The mean frequency of
$h_j^\pm$ is $s_j-\pi/L+o(L^{-1})$, giving
\eqref{eq:euler-uniform-large-gap}. For $h_j^+$ and $h_{j+1}^-$ it is
$s_j+o(s_j/L)$, giving the mass formula in
\eqref{eq:euler-uniform-small-gap}.
\end{proof}

\begin{corollary}\label{cor:euler-gap-regimes}
If $j/L$ stays in a fixed compact subinterval of
$(0,\infty)$, then uniformly
\begin{equation}\label{eq:euler-fine-splitting}
 \begin{aligned}
 g_j-h_j^+&=\frac{2s_ja_j}{L}+o(L^{-1}),&
 h_{j+1}^--g_j&=\frac{2s_ja_j}{L}+o(L^{-1}),\\
 h_j^+-h_j^-&=\frac{2s_j(\pi-2a_j)}{L}+o(L^{-1}).
 \end{aligned}
\end{equation}
For each fixed $j\ge2$,
\begin{equation}\label{eq:euler-fixed-j-limit}
 \begin{gathered}
 s_j\sim\frac{2\pi j}{L},\\
 h_j^+-h_j^-\sim\frac{2\pi^2(2j-1)}{L^2},\qquad
 h_{j+1}^--h_j^+\sim\frac{16\pi^2a_1j^2}{L^3}.
 \end{gathered}
\end{equation}
\end{corollary}
\begin{proof}
Theorem~\ref{thm:euler-density} places all adjacent periodic roots in a fixed
supercritical compact interval for the bulk range, where
Theorem~\ref{thm:euler-uniform-localization} gives \eqref{eq:euler-fine-splitting}.
For fixed $j$, Corollary~\ref{cor:euler-first-pair} gives $s_j\to0$.
Thus \eqref{eq:euler-absolute-position} yields $s_j\sim2\pi j/L$, and
Lemma~\ref{lem:bessel-inner-phase} gives $a(s_j)\sim a_1s_j$.
Substitution proves \eqref{eq:euler-fixed-j-limit}; its $L^{-2}$
coefficient retains the $-\pi/L$ correction in
\eqref{eq:euler-uniform-large-gap}.
\end{proof}

\paragraph{Comparison with the Floquet stability parameters.}
In \cite{MSS06b} the Euler mass parameter is $\beta_c=\delta$ and the
small energy parameter is $\varepsilon=(1-e^2)/2$. At fixed
$1/8<\delta<7$ in the physical parameter range, their Theorem~1
describes the full-period quantity
$\operatorname{tr}_2=\mu+\mu^{-1}$ for a reciprocal Floquet pair,
assuming $d_g e_g\ne0$. Lemma~\ref{lem:euler-mss-matching} verifies this hypothesis. Their imaginary exponent is $i\gamma_2$ with
$\gamma_2=\sqrt{8\delta-1}=2\sigma(\delta)$.
Write $\Sigma$ for the transition scale defined after Eq.~(26)
in the authors' version of \cite{MSS06b}.
Integrating $q'/q=-Q/2$ between the fixed sections and using their
Lemma~3(b) gives $2\log\Sigma=\log(1/q_0)+O(1)$,
where $q_0=\sqrt{1-e}$. The factor $\Sigma^{2i\gamma_2}$ in
\cite[Eq.~(33)]{MSS06b} thus has phase $\sigma(\delta)L+O(1)$,
twice the leading half-interval phase in \eqref{eq:euler-absolute-position}.
With this energy parameter, we obtain $\gamma_2/2$ as the coefficient
of $\log(1/\varepsilon)$ from Eq.~(33), whereas the last line of
\cite[Eq.~(7)]{MSS06b} displays $\gamma_2$.

The phase determination in
Theorem~\ref{thm:euler-uniform-localization} also identifies the absolute
curve labels and resolves the two distinct gap scales in
\eqref{eq:euler-fixed-j-limit} as the mass approaches $1/8$.

\subsection{Euler geometry near the Hardy threshold}
The same relative phase control also locates individual curves near
the reference rays. The following coordinates put the mass parameter in
the horizontal direction, as in Figure~\ref{fig:euler-quantile}.

\begin{corollary}[Quadratic displacements from the reference rays]\label{cor:euler-weighted-rays}
Set $\xi=\sigma(\delta)$, $\eta=L^{-1}$,
$r=(\xi^2+\eta^2)^{1/2}$ and $k_j=2\pi j$. With $a_1=a'(0^+)>0$ and
$b_0=B'(0^+)$, put
\[
 c=b_0+a_1,\ b_0,\ b_0-a_1
\]
for the branches $h_j^+,g_j,h_{j+1}^-$, respectively.
There exist $r_0>0$ and one modulus $\varepsilon(r)\to0$, independent
of $j\ge2$ and these three choices, such that every corresponding
point on these branches with $0<r<r_0$ satisfies
\begin{equation}\label{eq:euler-weighted-position}
 \left|\eta-\frac\xi{k_j}-\frac{2c\xi^2}{k_j^2}\right|
 \le\varepsilon(r)\frac{\xi^2}{k_j^2}.
\end{equation}
In particular the signed displacement from the reference ray
$R_j:\eta=\xi/k_j$ has the uniform normalized limit
\begin{equation}\label{eq:euler-scaled-ray-deviation}
 D_j:=\frac{k_j^2}{2\xi^2}\left(\eta-\frac\xi{k_j}\right)
 =c+o(1).
\end{equation}
\end{corollary}
\begin{proof}
By Remark~\ref{rem:euler-critical-chart} and
Lemma~\ref{lem:bessel-inner-phase}, we have $B(s)=b_0s+o(s)$ and
$a(s)=a_1s+o(s)$ (with closed forms in \eqref{eq:euler-b0} and
\eqref{eq:euler-inner-shift-value}).
Equation~\eqref{eq:euler-absolute-position} and the gap formulas at $g_j$
in Theorem~\ref{thm:euler-uniform-localization} give
\begin{equation}\label{eq:euler-joint-parabolas}
 \xi=k_j\eta-2k_jc\eta^2+o(k_j\eta^2),
 \qquad \eta\to0,\quad k_j\eta\to0,
\end{equation}
uniformly on the range $g_{j+1}\le d$.

To check this range before using the expansion, choose fixed bounds
$1/8<d_*<d_1<d_0<\infty$ within the critical chart. Counting and periodic
spacing supply a periodic root $g_K\in(d_*,d_1)$ with $g_{K+1}<d_0$ for
large $L$. Points in a sufficiently small disk have $\delta<d_*$;
interlacing places the corresponding periodic roots below $d_0$. The coarser
critical relation is already $k_j=\xi/\eta+O(\xi)$, so
$k_j\eta\asymp\xi$ uniformly. Thus $r\to0$ puts these points in
the joint regime of \eqref{eq:euler-joint-parabolas}.

Put $\alpha=\xi/k_j$. Dividing \eqref{eq:euler-joint-parabolas} by
$k_j$ gives $\alpha=\eta-2c\eta^2+o(\eta^2)$ uniformly.
Since $\eta/\alpha=1+O(\eta)$, inversion yields
$\eta=\alpha+2c\alpha^2+o(\alpha^2)$. The coarse comparison above
turns this joint remainder into one radial modulus, proving
\eqref{eq:euler-weighted-position} and
\eqref{eq:euler-scaled-ray-deviation}.
\end{proof}

The uniform $o(r^2/j^2)$ remainder in \eqref{eq:euler-weighted-position}
resolves the separation of branches along a common ray even for large
$j$; both the near-critical and bulk regimes are illustrated in
Figure~\ref{fig:euler-quantile} in the introduction.

\subsection{Floquet indices and stability boundaries for the Lagrange triangle}

For the Lagrange triangle the parameter and reduced matrix are
\begin{equation}\label{eq:lagrange-R-varying}
 \beta=\frac{27(m_1m_2+m_2m_3+m_3m_1)}{(m_1+m_2+m_3)^2},\qquad
 R_\beta=\operatorname{diag}\left(
 \frac{3+\sqrt{9-\beta}}2,\frac{3-\sqrt{9-\beta}}2\right).
\end{equation}
Both eigenvalues are nonnegative for $0\le\beta\le9$.
Consequently $C(R_\beta)=0$, and all essential Floquet indices
are bounded as $e\uparrow1$, uniformly in $\beta\in[0,9]$ and
$|\omega|=1$. For $\beta>0$, \cite{HLS14} gives the stronger conclusion
$i_1=0$ and periodic nondegeneracy.

Specifically, \cite[Theorem 1.2]{HLS14} gives
$\beta_s(e)\le\beta_m(e)\le\beta_k(e)$: two antiperiodic curves
and a hyperbolicity envelope. Their open regions, in increasing
$\beta$, are strongly stable, elliptic--hyperbolic, strongly
stable and hyperbolic, with empty intervals omitted at coincidences.
All three boundaries tend to $0$ as $e\uparrow1$. Their locations
depend on bounded-index behavior, invisible in a leading
logarithmic coefficient that vanishes identically.

\subsection{Symmetry blocks and stability criteria}\label{subsec:universality}

Formula \eqref{eq:ere-C} applies directly to every planar central
configuration for which the essential reduction
\eqref{eq:ere-family} has been established. If a mass-orthogonal
symmetry decomposition gives
$R=\bigoplus_{\ell=1}^q R_\ell$, then
\begin{equation}\label{eq:ere-block-coefficient}
 C(R)=\sum_{\ell=1}^q C(R_\ell).
\end{equation}
The identity uses only the spectrum of $R$. If the subspaces are also
$J$-invariant, the differential operator decomposes and the indices
add exactly for each $e<1$. For a first-order symplectic reduction,
the induced kinetic form and domain determine whether the second-order
realizations are unitarily equivalent.

For regular $(1+n)$-gon configurations, the symmetry reduction
of Hu--Long--Ou \cite{HLO20} supplies mass-dependent essential
blocks. Their spectra can be inserted into
\eqref{eq:ere-block-coefficient}, giving the sum of all
supercritical contributions and the bounded-index exclusion
\eqref{eq:ere-index-exclusion} whenever that sum is positive.
On a compact mass interval with every eigenvalue separated from
$-1/8$, the error is uniform by
Theorem~\ref{thm:ere-uniform-floquet}. If all blocks are at or above
the Hardy threshold, the index stays bounded for each fixed mass.

Hu--Sun
\cite{HS10} develop the Morse-index approach to elliptic Lagrange
motions; the large-central-mass results of \cite{HLO20} and the
trace-formula criteria of \cite{HOW15,HOS25} give sufficient stability
or nondegeneracy conditions. To pass from an index count to curve
locations for another configuration, one needs the mass dependence of
its essential blocks and the multiplicities of its crossings. The graph-continuity
criterion of Theorem~\ref{lem:volterra} requires global $L^2$ control
of the solution frame. As
Proposition~\ref{prop:threshold} shows, a component becoming limit-point
can change the minimal Fredholm index.

\appendix

\section{Continuous fundamental frame for the regularized Bessel equation}
\label{app:quadratic-proof}

\begin{proof}[Proof of Theorem~\ref{thm:bessel2}]
We apply Theorem~\ref{lem:volterra}. It is enough to construct two
real homogeneous solutions with Wronskian one that depend jointly
continuously on $(s,r)$ in $L^2(0,1)$. The difficulty is to keep this
continuity both at the singular limit $s=0$ and at the critical value
$r=0$, where the two power solutions coalesce.

We use the Gauss hypergeometric function and Pfaff's transformation
with the branch issued from $x=0$ \cite[Sections~15.8 and~15.10]{DLMF}.

\smallskip
\noindent\emph{Step 1: a common representation of the solutions.}
For $z\in\mathbb C$ with $|\operatorname{Re}z|<1$, define, for
$0<t\le1$ and $0\le s\le1$,
\begin{equation}\label{eq:Phi}
\Phi(z,s,t)=(t^2+s^2)^{1/4+z/2}
{}_2F_1\left(-\frac14-\frac z2,\frac34-\frac z2;
             1-z;\frac{s^2}{t^2+s^2}\right).
\end{equation}
The powers use the positive real base. The representation is useful
because the argument of the hypergeometric function lies in $[0,1)$,
even when $s$ tends to zero. Pfaff's transformation gives the
equivalent expression
\begin{equation}\label{eq:Phi-infinity}
\Phi(z,s,t)=t^{1/2+z}
{}_2F_1\left(-\frac14-\frac z2,\frac14-\frac z2;
             1-z;-\frac{s^2}{t^2}\right).
\end{equation}
Indeed, for $s>0$ the substitution
$y(t)=t^{1/2+z}F(-s^2/t^2)$ transforms
\begin{equation}\label{eq:Phi-ode}
-y''+\frac{z^2-1/4}{t^2+s^2}y=0
\end{equation}
into the Gauss equation with parameters
$a=-1/4-z/2$, $b=1/4-z/2$, and $c=1-z$;
see \cite[Section~15.10]{DLMF}. Formula~\eqref{eq:Phi-infinity}
then follows, and Pfaff's identity
\cite[Section~15.8]{DLMF} yields~\eqref{eq:Phi}.
At $s=0$, both representations reduce to $t^{1/2+z}$.

\smallskip
\noindent\emph{Step 2: uniform $L^2$ control.}
Let $K$ be a compact subset of $|\operatorname{Re}z|<1$.
For the hypergeometric series $H_z(x)$ in~\eqref{eq:Phi},
write
\[
H_z(x)=\sum_{m=0}^\infty c_m(z)x^m,
\qquad
c_m(z)=\frac{(a)_m(b)_m}{(c)_m m!},
\]
where now $a=-1/4-z/2$, $b=3/4-z/2$ and $c=1-z$.
The relation $c-a-b=1/2$ and the gamma-ratio estimates
\cite[Section~5.11]{DLMF} give
\[
\sup_{z\in K}|c_m(z)|\le C_K(m+1)^{-3/2}.
\]
At terminating parameters the estimate follows by continuity of
the coefficients. Thus $H_z$ is uniformly bounded for $z\in K$
and $0\le x\le1$. Applying Cauchy's estimate in $z$ on a
slightly larger compact subset of the strip also bounds
$\partial_zH_z$ uniformly.

Choose $\rho$ with
$\max\{1/2,\sup_{z\in K}|\operatorname{Re}z|\}<\rho<1$.
The power factor in~\eqref{eq:Phi}, and its $z$-derivative,
then give
\begin{equation}\label{eq:Phi-bound}
|\Phi(z,s,t)|+|\partial_z\Phi(z,s,t)|
\le C_K\bigl(1+t^{1/2-\rho}(1+|\log t|)\bigr),
\end{equation}
uniformly for $z\in K$, $0\le s\le1$, and $0<t\le1$.
The right-hand side belongs to $L^2(0,1)$ because $\rho<1$.
Pointwise continuity for $t>0$, followed by dominated convergence,
therefore gives joint $L^2$ continuity of $\Phi$ and
$\partial_z\Phi$ in $(z,s)$, including at $s=0$.

\smallskip
\noindent\emph{Step 3: a real normalized frame.}
Set $z=r$ for $r\ge0$ and $z=ir$ for $r<0$, so that
$z^2-1/4=q(r)$. Define
\begin{equation}\label{eq:normalized-frame}
\begin{aligned}
u_{s,r}&=\frac{\Phi(z,s,\cdot)+\Phi(-z,s,\cdot)}2,\\
w_{s,r}&=\begin{cases}
\displaystyle\frac{\Phi(z,s,\cdot)-\Phi(-z,s,\cdot)}{2z},&z\ne0,\\[4pt]
\partial_z\Phi(0,s,\cdot),&z=0.
\end{cases}
\end{aligned}
\end{equation}
These solutions are real: for imaginary $z$, the terms with
$z$ and $-z$ are complex conjugates. At $z=0$, the
parameter derivative also solves~\eqref{eq:Phi-ode}, since the
coefficient depends on $z^2$. Furthermore,
\[
\frac{\Phi(z,s,t)-\Phi(-z,s,t)}{2z}
=\frac12\int_{-1}^1\partial_z\Phi(\theta z,s,t)\,d\theta.
\]
Together with~\eqref{eq:Phi-bound}, this proves that
$(s,r)\mapsto(u_{s,r},w_{s,r})$ is continuous in
$L^2(0,1)^2$, uniformly across $r=0$ where the two indicial roots coalesce.

It remains to check the normalization. Since~\eqref{eq:Phi-ode}
contains no first derivative, the Wronskian is constant in $t$.
For fixed $s$ and $z\ne0$, formula~\eqref{eq:Phi-infinity}
shows, as $t\to\infty$, that
$\Phi(\pm z,s,t)\sim t^{1/2\pm z}$; their derivatives
have the corresponding asymptotics. The Wronskian of these
solutions is therefore $-2z$. The combinations
in~\eqref{eq:normalized-frame} consequently satisfy
\[
W(u_{s,r},w_{s,r})=1.
\]
The same identity holds at $z=0$ by continuity. Thus the
associated phase-space frame is symplectic. Theorem~\ref{lem:volterra}
applies and yields gap continuity, Fredholmness, and compactness
of the maximal graph embedding.
\end{proof}

\section{Reversible involution symmetry and the half-interval formulation}
\label{app:reversible}

We use the full-period operator and normalization of
Section~\ref{subsec:ere-brake}; in particular, $R$, $J$, and
$D_\omega$ are defined in~\eqref{eq:ere-normalized-R}--\eqref{eq:ere-floquet-domain}.

While \eqref{eq:ere-family} is formulated on the full period $[0,2\pi]$,
the coefficient becomes singular as $e\uparrow1$ at $t=\pi$,
an interior point of $[0,2\pi]$.
In the true-anomaly convention used in \eqref{eq:ere-family},
the Kepler radius is proportional to $(1+e\cos t)^{-1}$:
for $e>0$, $t=0$ is periastron and $t=\pi$ is apoastron.
The singularity of the normalized coefficient at $t=\pi$ must therefore
be distinguished from the physical periastron collision.
The half-interval reduction comes from the discrete time-reversal
involution symmetry ($\mathcal{R}$-symmetry) of the linearized
Hamiltonian system in the rotating pulsating frame \cite{HS10,HO16,ZL17}.
Under the reflection $t\mapsto-t$ (or $t\mapsto 2\pi-t$), the Keplerian radius
and the potential factor $1+e\cos t$ are even.
For the reversible subclass, assume there is an orthogonal involution
$S$ with $SJS=-J$ and $SR=RS$. In coordinates where this is possible, let $S=\operatorname{diag}(S_2,\ldots,S_2)$ with
$S_2=\left(\begin{smallmatrix}1&0\\0&-1\end{smallmatrix}\right)$.
Then $S$ is symmetric and orthogonal ($S^T=S=S^{-1}$), reverses the complex structure,
$SJS=-J$, and commutes with $R$ whenever $R$ preserves the coordinate axes
(which holds for the Euler diagonal matrix in Section~\ref{subsec:euler-bulk}, but is an additional hypothesis for general central configurations).
The reflection involution on $H^2((0,2\pi),\C^d)$ defined by
\[
 (\mathcal R y)(t) = S y(2\pi-t)
\]
satisfies $\mathcal R^2=I$ and commutes with the differential operator:
$A_{R,e}\mathcal R = \mathcal R A_{R,e}$.
For periodic ($\omega=1$) and antiperiodic ($\omega=-1$) solutions, $\mathcal R$
preserves the domain $D_\omega$, decomposing the problem orthogonally into
the even eigenspace ($\mathcal R y = y$) and the odd eigenspace ($\mathcal R y = -y$).
If $\mathcal R y=\varepsilon y$, $\varepsilon\in\{1,-1\}$, then
\[
\begin{aligned}
 Sy(0)&=\varepsilon\omega y(0),&
 Sy'(0)&=-\varepsilon\omega y'(0),\\
 Sy(\pi)&=\varepsilon y(\pi),&
 Sy'(\pi)&=-\varepsilon y'(\pi).
\end{aligned}
\]
Thus the endpoint conditions at $t=0$ interchange between periodic
and antiperiodic boundary conditions, whereas at $t=\pi$ they depend
only on $\varepsilon$.

Restriction to the half-interval, multiplied by $\sqrt2$, is unitary
on each symmetry sector. This symmetry reduction converts
the periodic and antiperiodic spectral problems into self-adjoint boundary value
problems on $[0,\pi]$ for the differential operator
\begin{equation}\label{eq:ere-half-operator}
 A_{R,e}^\pm y = -y'' - 2Jy' + \frac{R}{1+e\cos t}y,\qquad 0<t<\pi.
\end{equation}
For $d=2$, we label the two resulting endpoint conditions by the
condition imposed on the first component $y_1$: $N$ (Neumann in $y_1$,
Dirichlet in $y_2$) and $D$ (Dirichlet in $y_1$, Neumann in $y_2$).
On $[0,\pi]$, the regular endpoint $t=0$ is separated from the singular
endpoint $t=\pi$:
\begin{enumerate}[label=(\roman*)]
\item \emph{Regular endpoint at $t=0$}: At $t=0$, $1+e\cos 0=1+e\in[1,2]$. The denominator is smooth and bounded away from zero, so the potential
is uniformly bounded there for fixed $R$.
For $d=2$, the involution relations at $t=0$ impose condition $N$ when
$\varepsilon\omega=1$ and condition $D$ when $\varepsilon\omega=-1$
(the general block-diagonal reversible case repeats these conditions in
each pair):
\begin{equation}\label{eq:apoastron-bc}
 \begin{aligned}
 N:\quad & y_1'(0)=0,\quad y_2(0)=0,\\
 D:\quad & y_1(0)=0,\quad y_2'(0)=0.
 \end{aligned}
\end{equation}
In canonical coordinates, each condition defines a constant Lagrangian
subspace of $\R^{2d}$ (Neumann in the first components and Dirichlet in
the second for $N$, or vice versa for $D$).
\item \emph{Singular coefficient at the endpoint $t=\pi$}:
At $t=\pi$, $1+e\cos\pi = 1-e\to0$ as $e\uparrow1$.
In the half-interval formulation, the singularity lies at the boundary
endpoint $t=\pi$ rather than in the interior.
Introducing the local distance to the boundary $\tau=\pi-t\in[0,\pi]$,
\begin{equation}\label{eq:periastron-tau}
 1+e\cos t = 1 - e\cos\tau = (1-e) + 2e\sin^2\frac{\tau}{2} = (1-e) + \frac{e}{2}\tau^2 + O(\tau^4).
\end{equation}
At $e=1$, the limiting potential on $[0,\pi]$ is
\begin{equation}\label{eq:poschl-teller-pot}
 \frac{R}{1+\cos t} = \frac{R}{2\cos^2(t/2)} = \frac{R}{2\sin^2(\tau/2)} = \frac{2R}{\tau^2} + O(1),
\end{equation}
which has an inverse-square endpoint singularity at $\tau=0$ ($t=\pi$).
For $e<1$, the boundary conditions at $t=\pi$ ($\varepsilon=1$ for $N$
and $\varepsilon=-1$ for $D$) define two Lagrangian subspaces in the
4-dimensional phase space $\R^4$ in coordinates $(y_1,y_2,y_1',y_2')$:
\begin{equation}\label{eq:brake-anti-brake-subspaces}
\begin{aligned}
 L_N &= \{(y_1,y_2,y_1',y_2')\in\R^4 : y_1' = 0,\ y_2 = 0\},\\
 L_D &= \{(y_1,y_2,y_1',y_2')\in\R^4 : y_1 = 0,\ y_2' = 0\}.
\end{aligned}
\end{equation}
Both are Lagrangian for the Green symplectic form expressed through
$(y'+Jy,y)$ (and for $y_1=0$, the condition $y_2'=0$ is equivalent to
$p=y_2'+2y_1=0$ in Section~\ref{subsec:euler-bulk}). In $(y,y')$
coordinates this is the pullback of the canonical form, including the
gyroscopic term.
They are orthogonal with respect to the standard Euclidean metric on $\R^4$ and hence transverse:
\begin{equation}\label{eq:boundary-transversality}
 L_N \cap L_D = \{0\}, \qquad L_N \oplus L_D = \R^4.
\end{equation}
For the Euler case $d=2$, pairing the conditions at $t=0$ and $t=\pi$
yields the four self-adjoint sectors $\alpha\beta\in\{NN,DD,ND,DN\}$ used
in Section~\ref{sec:ere-applications}:
\begin{itemize}
\item \emph{Periodic ($\omega=1$)}: sector $NN$ ($\varepsilon=+1$: $L_N$
at $t=0$ and $L_N$ at $t=\pi$) and sector $DD$ ($\varepsilon=-1$: $L_D$
at $t=0$ and $L_D$ at $t=\pi$).
\item \emph{Antiperiodic ($\omega=-1$)}: sector $ND$ ($\varepsilon=-1$:
$L_N$ at $t=0$ and $L_D$ at $t=\pi$) and sector $DN$ ($\varepsilon=+1$:
$L_D$ at $t=0$ and $L_N$ at $t=\pi$).
\end{itemize}
\end{enumerate}
Folding fixes the half-period $[0,\pi]$ and the Hilbert space
$L^2((0,\pi),\C^d)$, with regular endpoint $t=0$ and singular endpoint
$t=\pi$. For Euler, the second component has coefficient
$4\delta+6>3/4$ and is limit-point for every $\delta\ge0$, while for
$\delta>1/8$ the first component is supercritical and limit-circle.
Lemma~\ref{lem:euler-full-family} keeps the strong component in its
energy domain and prescribes the weak boundary line. The Cauchy lines
of these coupled realizations enter Theorem~\ref{thm:euler-fine-splitting}.

\section{Connection formulas, phase slopes, and matching coefficients}
\label{app:scalar-phases}
\label{subsec:euler-angles-splitting}
\paragraph{The inner Dirichlet and Neumann constants.}
In the boundary angle formulas \eqref{eq:bessel-DN-angles} of
Lemma~\ref{lem:bessel-inner-phase}, both angles contain the rotation
$\sigma\log s$. It cancels in their difference, so
$a(\sigma)=\arg(U'_{i\sigma}/U_{i\sigma})$ by
\eqref{eq:inner-gudermannian-shift}, with the argument chosen to vanish
at $\sigma=0$. Here $U_z$ and $U'_z$ are the value and right derivative
of $\Phi(z,1,t)$ at $t=0$. To compute them we use Gauss's connection
formula \cite[equation~15.10.21]{DLMF},
\begin{multline}\label{eq:gauss-connection}
 {}_2F_1(a,b;c;z)=\frac{\Gamma(c)\Gamma(c-a-b)}{\Gamma(c-a)\Gamma(c-b)}
 {}_2F_1(a,b;a+b-c+1;1-z)\\
 +(1-z)^{c-a-b}\frac{\Gamma(c)\Gamma(a+b-c)}{\Gamma(a)\Gamma(b)}
 {}_2F_1(c-a,c-b;c-a-b+1;1-z).
\end{multline}

At $s=1$ the hypergeometric argument in \eqref{eq:Phi} is
$x=(1+t^2)^{-1}$. Its parameters are $a=-1/4-z/2$, $b=3/4-z/2$
and $c=1-z$. Since $c-a-b=1/2$ and $(1-x)^{1/2}=t+O(t^3)$
as $t\downarrow0$, the first and second terms of
\eqref{eq:gauss-connection} give, respectively, the constant and linear
coefficients in $\Phi(z,1,t)=U_z+tU'_z+O(t^2)$. For $z=i\sigma$ they are
\begin{equation}\label{eq:U-Up-gamma}
 U_{i\sigma}=\frac{\Gamma(1-i\sigma)\Gamma(1/2)}{\Gamma\left(\frac54-\frac{i\sigma}2\right)\Gamma\left(\frac14-\frac{i\sigma}2\right)},\qquad
 U'_{i\sigma}=\frac{\Gamma(1-i\sigma)\Gamma(-1/2)}{\Gamma\left(-\frac14-\frac{i\sigma}2\right)\Gamma\left(\frac34-\frac{i\sigma}2\right)}.
\end{equation}
Let $\zeta=1/4+i\sigma/2$ and
$Z=\Gamma(1-\bar\zeta)/\Gamma(\zeta)$. The recurrence
$\Gamma(1+w)=w\Gamma(w)$ and the identity
$\Gamma(-1/2)=-2\Gamma(1/2)$ give
\[
 \frac{U'_{i\sigma}}{U_{i\sigma}}
 =\frac{1/4+\sigma^2}{2}
   \left(\frac{\Gamma(\bar\zeta)}{\Gamma(1-\zeta)}\right)^2
 =\frac{1/4+\sigma^2}{2\overline Z^{\,2}} .
\]
The prefactor is positive, so $a(\sigma)=2\arg Z$. Euler's reflection
formula gives
\[
 Z=\frac{\pi}{|\Gamma(\zeta)|^2\sin(\pi\bar\zeta)},\qquad
 \arg Z=-\arg\sin\!\left(\frac\pi4-\frac{i\pi\sigma}2\right)
       =\arctan\!\left(\tanh\frac{\pi\sigma}2\right).
\]
Doubling this angle on the branch through zero gives
\begin{equation}\label{eq:euler-inner-shift-value}
 a(\sigma)=\arctan(\sinh\pi\sigma),\qquad a_1=a'(0^+)=\pi.
\end{equation}
In particular $\pi-2a(\sigma)=4\arctan(e^{-\pi\sigma})$ and
the coefficient $16\pi^2a_1j^2$ in \eqref{eq:euler-fixed-j-limit}
is $16\pi^3j^2$.

\paragraph{The scalar comparison phases.}
As an explicitly solvable reference model, we record the even and odd
boundary phases of the decoupled P\"oschl--Teller equation, where the
same Gudermannian function $\operatorname{gd}(\pi\sigma)$ governing the
inner shift $a(\sigma)$ in \eqref{eq:euler-inner-shift-value} also gives
the even--odd phase separation. For $\delta>1/8$, put
$\sigma=\sqrt{2\delta-1/4}$ and consider the decoupled equation
\begin{equation}\label{eq:euler-limit-ode}
 -y''-\frac{\delta}{2\cos^2(t/2)}y=0,\qquad 0<t<\pi.
\end{equation}
With $\gamma=1/2+i\sigma$ and $z=\sin^2(t/2)$, the substitution
$y=(1-z)^{\gamma/2}F(z)$ gives the unit-normalized even and odd
solutions with $(y(0),y'(0))=(1,0)$ and $(0,1)$, respectively,
\[
 \begin{aligned}
 y_{\rm even}(t)&=(\cos(t/2))^\gamma {}_2F_1(\gamma/2,\gamma/2;1/2;z),\\
 y_{\rm odd}(t)&=2\sin(t/2)(\cos(t/2))^\gamma
 {}_2F_1((\gamma+1)/2,(\gamma+1)/2;3/2;z).
\end{aligned}
\]
As $\tau=\pi-t\downarrow0$, one has $1-z=\tau^2/4+O(\tau^4)$. Applying
\eqref{eq:gauss-connection} with $c-a-b=-i\sigma$ yields
$y_\epsilon=2\operatorname{Re}[C_{-,\epsilon}(\tau/2)^{1/2-i\sigma}]
+O(\tau^{5/2})$, with differentiated error $O(\tau^{3/2})$, where
\[
 C_{-,\mathrm{even}}=\frac{\sqrt\pi\Gamma(i\sigma)}
 {\Gamma(1/4+i\sigma/2)^2},\qquad
 C_{-,\mathrm{odd}}=\frac{\sqrt\pi\Gamma(i\sigma)}
 {\Gamma(3/4+i\sigma/2)^2}.
\]
In the Wronskian-one frame $(u,w)$ of \eqref{eq:euler-frame},
\[
 (\tau/2)^{1/2-i\sigma}=2^{-1/2}2^{i\sigma}(u-i\sigma w).
\]
The two remainder estimates give zero limiting brackets for the
remainder against both $u$ and $w$. Choose the arguments of
$C_{-,\epsilon}$ with common limit $-\pi/2$ as $\sigma\downarrow0$,
and put $\Phi_\epsilon=\arg C_{-,\epsilon}+\sigma\log2$.
The coefficient line is $[\cos\Phi_\epsilon:\sigma\sin\Phi_\epsilon]$,
so its angle is taken modulo $\pi$.
To add the two angles, use Legendre's duplication formula
\[
 \Gamma(1/4+i\sigma/2)\Gamma(3/4+i\sigma/2)
 =2^{1/2-i\sigma}\sqrt\pi\,\Gamma(1/2+i\sigma).
\]
It gives $\Phi(\delta)=\Phi_{\rm even}+\Phi_{\rm odd}$ in the form
\begin{equation}\label{eq:total-phi-closed}
 \Phi(\delta) = 2\arg\left(\frac{\Gamma(i\sigma)}{\Gamma\left(\frac12+i\sigma\right)}\right) + 4\sigma\log 2 \pmod\pi.
\end{equation}

\begin{proposition}\label{prop:scalar-phase-separation}
For the scalar equation \eqref{eq:euler-limit-ode}, choose continuous
phase lifts agreeing at $\sigma=0$. Then
\begin{equation}\label{eq:gudermannian-exact}
\begin{aligned}
 \Phi_{\rm even}-\Phi_{\rm odd}
 &=2\arg\frac{\Gamma(3/4+i\sigma/2)}{\Gamma(1/4+i\sigma/2)}
 =2\arctan\!\left(\tanh\frac{\pi\sigma}{2}\right)\\
 &=\operatorname{gd}(\pi\sigma)
 =\arctan(\sinh\pi\sigma)\in(0,\pi/2).
\end{aligned}
\end{equation}
Here $\operatorname{gd}$ denotes the classical Gudermannian function.
\end{proposition}
\begin{proof}
In the difference $\arg C_{-,\rm even}-\arg C_{-,\rm odd}$, the
common numerator cancels. The remaining angle is $2\arg Z$, evaluated
above on the lift through zero.
\end{proof}
\begin{remark}
Abel's identity and the initial data give
$W_\tau(y_{\rm odd},y_{\rm even})=1$ in the coordinate $\tau=\pi-t$.
Evaluating at the singular end yields
\[
 2\sigma |C_{-,\rm even}|\,|C_{-,\rm odd}|
 \sin(\Phi_{\rm even}-\Phi_{\rm odd})=1.
\]
Thus the two scalar lines are transverse. For the Euler spacings we
compare the two terminal conditions in a single regular sector.
They have the same limiting outer phase by
Lemma~\ref{lem:euler-uniform-phases}.
\end{remark}

\paragraph{Critical slope of the coupled Euler phase.}
We now return to the coupled Euler system and use its reconstruction variable $u=p/M$ to compute the coefficient $b_0$ in the quadratic displacement from the reference rays.
Let $B(\sigma)=\theta_1^N-\sigma\log2/2-\arg U'_{i\sigma}-\pi/2$,
where $\theta_1^N\to\pi/2$ and $\arg U'_0=0$. This is the limiting
$NN$ slow phase of Remark~\ref{rem:euler-critical-chart}, with $B(0)=0$.
At the critical point the regular $N$ solution in
\eqref{eq:euler-exact-reconstruction} is proportional to
\[
 u(t)=\sqrt{\cos(t/2)}\left[
 \mathsf P_{3/2}(\sin(t/2))+\frac2\pi
 \mathsf Q_{3/2}(\sin(t/2))\right].
\]
Indeed, since $\mathsf Q_{3/2}(z)=(\pi/2)\mathsf P_{3/2}(-z)$, the bracket
equals $\mathsf P_{3/2}(\sin(t/2))+\mathsf P_{3/2}(-\sin(t/2))$, so $u$ is
even in $t$ and satisfies $u'(0)=0$.
The Ferrers endpoint formulas \cite[Sections~14.8--14.9]{DLMF14} give
$\mathsf P_{3/2}(\sin(t/2))=1+o(1)$ and
$\mathsf Q_{3/2}(\sin(t/2))=-\log\tau+4\log2-8/3+o(1)$.
Using the trace map \eqref{eq:euler-critical-trace-map} yields
\[
 H_N:=\frac ab=-\frac\pi2-4\log2+\frac{160}{69}.
\]
The critical chart therefore gives
$\theta_1^N=\pi/2-\sigma H_N+o(\sigma)$.
Logarithmic differentiation of \eqref{eq:U-Up-gamma}, using the
digamma reflection and recurrence identities, gives
$\left.\frac{d}{d\sigma}\arg U'_{i\sigma}\right|_0
=2+\pi/2-3\log2$. Consequently
\begin{equation}\label{eq:euler-b0}
 B(\sigma)=b_0\sigma+o(\sigma),\qquad
 b_0=-H_N-\tfrac12\log2-(2+\pi/2-3\log2)
 =\frac{13}{2}\log2-\frac{298}{69}.
\end{equation}
Thus $b_0=0.18661609\ldots$, and the three coefficients in
Corollary~\ref{cor:euler-weighted-rays} are $b_0+\pi,b_0,b_0-\pi$.
The $+1$ term of the coupled reconstruction enters through $H_N$.

\begin{lemma}[Nonvanishing of the Euler matching coefficients]\label{lem:euler-mss-matching}
Fix $1/8<\delta<\infty$. For the Euler operator family, the coefficients
$d_g$ and $e_g$ of \cite[Lemma~5]{MSS06b} are both nonzero for every
sufficiently small matching cutoff.
\end{lemma}
\begin{proof}
We first show $\Omega(U,\mathcal RU)\ne0$ via the defect energy
\eqref{eq:euler-defect-energy}, next prove nondegeneracy of
$\Omega(\cdot,\mathcal R\cdot)$ on $E=\Span_{\mathbb C}(U,V)$ by a
Hermitian signature argument, and finally pass from the matching cutoff
of \cite{MSS06b} to the midpoint using exterior powers.

Translate the singularities to $0,2\pi$ and put $g=1-\cos t$.
In the order $Y=(y_2,y_1)^T$ of \cite{MSS06b}, the equation is
$Y''+2J_*Y'=g^{-1}\operatorname{diag}(2\delta+3,-\delta)Y$,
where $J_*=\left(\begin{smallmatrix}0&1\\-1&0\end{smallmatrix}\right)$.
Write $S=\operatorname{diag}(1,-1)$ and
$(\mathcal RY)(t)=S Y(2\pi-t)$. This reflection reverses the conserved
bilinear form
$\Omega(Y,Z)=Y^TZ'-Y'^TZ+2Y^TJ_*Z$.
At the midpoint, in coordinates $(Y,\sqrt g\,Y')$, the matrix of
$\sqrt{g(\pi)}\,\Omega(\cdot,\mathcal R\cdot)$ is precisely $G_M$
in \cite[Eq.~(27)]{MSS06b}. Write $X_1=(x_{ij})_{i,j=1,2}$ for
the matrix of this pairing on the first two columns transported from
the matching cutoff. The coefficients in \cite[Lemma~5]{MSS06b} are
$d_g=x_{11}$ and $e_g=\det X_1$.

Let $p=1/2+\sqrt{25/4+4\delta}$ and take the real solution
$U=t^p(e_1+O(t))$ at the left endpoint. Its defect $W$ from
Lemma~\ref{lem:euler-limit-kernel} has leading term
$\frac{p}{2\delta+3}[(p-1)(p-2)+2\delta]t^{p-3}$.
At the right endpoint, pairing with $\mathcal RU$ detects the
coefficient of $(2\pi-t)^{1-p}e_1$: the two exponents sum to one
and their derivative factors differ by $2p-1$. If
$\Omega(U,\mathcal RU)=0$, the growing defect branch is therefore
absent there, by the reconstruction in Lemma~\ref{lem:euler-limit-kernel}.
The boundary term in \eqref{eq:euler-defect-energy} then vanishes at
both ends. Integration gives $W=0$, contradicting its left leading term.

Take the oscillatory solution $V$ from \eqref{eq:euler-exact-reconstruction},
normalized to have leading term $t^{1/2+i\sigma}e_2$, where
$\sigma=\sqrt{2\delta-1/4}$. The left endpoint orders give
$\Omega(U,V)=0$ and $\Omega(\overline V,V)=2i\sigma$.
Thus $E=\operatorname{span}_{\mathbb C}(U,V)$ is Lagrangian, and the
Hermitian form $i\Omega(\bar z,z)$ on $E$ is semidefinite with null
line $\mathbb C U$. On $\mathcal RE$ its sign is reversed. A vector in
$E\cap\mathcal RE$ would therefore lie in both $\mathbb C U$ and
$\mathbb C\mathcal RU$; their pairing is nonzero, so the intersection
is trivial. Hence $\Omega(\cdot,\mathcal R\cdot)$ restricted to $E$
is nonsingular.

It remains to pass to the cutoff in \cite[Eqs.~(21), (27)]{MSS06b}.
In their time $dt=\sqrt g\,d\tau$, let $\Phi(0,-T)$ carry the cutoff
to the midpoint. As $\tau\to-\infty$ the coefficient matrix differs
from its diagonalizable limit by $O(e^{c\tau})$ for some $c>0$.
In degrees $k=1,2$, the leading data of $U$ and $U\wedge V$ are,
up to normalization, the limiting eigenvectors $P_1$ and $P_1\wedge P_2$
in the notation of \cite{MSS06b}. Their exponents $\lambda_k$ have
maximal real part in the respective degrees. Variation of constants
gives $\|\bigwedge^k\Phi(0,-T)\|\le C e^{T\Re\lambda_k}$.
Applying this bound to the $o(1)$ error in the normalized leading data
shows that the transported eigenvectors, multiplied by
$e^{-T\lambda_k}$, converge to the corresponding midpoint data.
Since $\det X_1$ is the induced pairing on the exterior product of
the two columns, the normalized limits of $d_g$ and $e_g$ are the
two nonzero midpoint pairings above, up to nonzero factors.
\end{proof}

\section{Subcritical and critical Euler sector indices}\label{app:euler-critical-index}
Let $c=1/8$ and write $(f,v)=(y_1,y_2)$. At $e<1$ the full form is
\[
 Q_{\delta,e}[f,v]=\int_0^\pi
 \left(|f'|^2-\frac{\delta}{q_e}|f|^2+|v'|^2+
       4\operatorname{Re}(v'\bar f)+\frac{2\delta+3}{q_e}|v|^2\right)dt .
\]
The essential condition is $v=0$ at an $N$ endpoint and $f=0$ at a
$D$ endpoint. Keep these conditions and impose $f(\pi)=0$ in the
form domain. Denote this restricted realization by
$\widehat A_{\delta,e}^{\alpha\beta}$. At $e=1$, take the closure in
form norm (after a scalar shift) of the same form
on smooth pairs satisfying the regular essential condition and
vanishing near collision. We write $F_\delta^\alpha$ for the
associated Friedrichs operator.

\paragraph{Continuity of the restricted forms.}
At $e=1$ and $0\le\delta<c$, Hardy's inequality gives the usual
$H^1$ form domain with $f(\pi)=v(\pi)=0$. At $\delta=c$ the weak form
$\mathfrak t_c[f]=\int(|f'|^2-c|f|^2/q_1)$ is understood on its
Friedrichs closure. Near collision $c/q_1=1/(4\tau^2)+O(1)$.
For a core function, write $f(\pi-\tau)=\sqrt\tau\,z(\tau)$.
Integration of the cross term gives
\[
 \int_0^\pi\left(|\partial_\tau(\sqrt\tau\,z)|^2-
               \frac{|\sqrt\tau\,z|^2}{4\tau^2}\right)d\tau
 =\int_0^\pi \tau |z'|^2\,d\tau+\tfrac12|z(\pi)|^2 .
\]
The shifted critical domain embeds compactly in $L^2$ by radial
Rellich compactness. Logarithmic cutoffs, linear in $\log\tau$ on
$[h^2,h]$, have radial energy $1/|\log h|$. Thus the critical
condition admits the $\sqrt\tau$ branch and excludes
$\sqrt\tau\log\tau$. This is the relaxed limit of the decreasing
subcritical weak forms \cite[Theorem~3.2]{Sim78}.

For $(\delta,e)\to(a,1)$, $0\le\delta\le c$, the weak form is
bounded below by
\[
 \mathfrak t_c[f]+(c-\delta)\int_0^\pi |f|^2/q_1\,dt
                  -C(1-e)\|f\|_2^2 .
\]
At $\delta=c$ the term with coefficient $c-\delta$ is omitted,
since its integral can diverge on the critical domain. The bound
follows from $q_e\ge q_1$ on the collision half and an $O(1-e)$
difference of reciprocals on the other half. The mixed term is bounded in modulus
by $\tfrac12\|v'\|_2^2+8\|f\|_2^2$. After a common shift, bounded form energy controls the critical norm of $f$,
the $H^1$ norm of $v$, and its weighted potential energy, so both components
are precompact in $L^2$. Along a convergent subsequence, Fatou's lemma gives
$v/\sqrt{q_1}\in L^2$; when $a<c$, it also gives $f/\sqrt{q_1}\in L^2$,
which combined with the radial energy implies $f\in H^1$ and $f(\pi)=0$.
Since the positive terms are weakly lower semicontinuous and the mixed term
converges by weak convergence of $v'$ and strong convergence of $f$,
we obtain lower semicontinuity of the forms, while recovery follows by
approximation on the common core away from collision. Together these yield
strong resolvent convergence; at a real point below the common lower bound
the resolvents are self-adjoint and collectively compact, and hence converge
in norm \cite[Proposition~2.1 and Theorem~3.4(a)]{AP68}. In particular,
\[
 \widehat A_{\delta,e}^{\alpha\beta}\longrightarrow F_a^\alpha
 \quad\text{in norm resolvent}\qquad
 ((\delta,e)\to(a,1),\ 0\le\delta\le c).
\]
The same argument on the edge $e=1$ gives continuity of
$F_\delta^\alpha$ up to $c$.

\paragraph{The singular limit and its continuation.}
At $\delta=0$ the periodic form is nonnegative \cite{ZL17}. Reflect
a half-interval core pair across $t=0$ with its prescribed parity.
It vanishes near both collision endpoints and is a periodic test
pair. Passing to $e=1$ on its support and then taking the form
closure gives $F_0^\alpha\ge0$ for both regular conditions. With $q=1+\cos t$,
integrating the first zero equation $-y_1''+2y_2'=0$ yields
$y_1'=2y_2+C$, reducing the system to
\[
 y_1'=2y_2+C,\qquad
 \left(-\partial_t^2+\frac3q-4\right)y_2=2C.
\]
The unique collision-admissible homogeneous line is $q\sin t$; the other
solution has exponent $-2$. For $\alpha=N$, $C=0$. For $\alpha=D$, the
admissible homogeneous solution fails the regular Neumann condition,
so the inhomogeneous problem is uniquely solvable for each $C$.
The kernels are exactly one-dimensional, generated by
\[
 v_N=\tfrac14(-q^2,q\sin t),\qquad
 v_D=((2+\cos t)\sin t,q\cos t).
\]
Their first variations are
$Q'_0[v]=\int_0^\pi q^{-1}(-|v_1|^2+2|v_2|^2)\,dt$, with
$Q'_0[v_N]=-3\pi/32$ and $Q'_0[v_D]=-3\pi/2$.
Analytic form perturbation on the fixed subcritical domain and the
positive gap above these simple zero eigenvalues give one negative eigenvalue
in each sector for small positive $\delta$.

Lemma~\ref{lem:euler-friedrichs-gap} excludes zero from both
Friedrichs spectra for $0<\delta\le c$. Norm-resolvent continuity
and the common lower bound therefore give $n_-(F_\delta^\alpha)=1$
throughout this interval.

\paragraph{The physical boundary correction.}
Fix $0<a\le c$. For $e<1$ close to one, the restricted realization
$\widehat A_{a,e}^{\alpha\beta}$ has negative index one and trivial
kernel. For terminal condition $D$, this already coincides with the
physical operator. For terminal condition $N$, the physical form domain
frees the single scalar trace $f(\pi)$, so any nonpositive subspace
intersects $\{f(\pi)=0\}$ in codimension at most one; hence the sum of
the physical negative index and nullity is at most two. To see that
this upper bound is attained by the negative index alone, choose a
negative test pair for $F_a^\alpha$ in the form core supported away
from collision (whose form value remains negative for $e$ close to one)
together with a smooth cutoff $\chi$ supported in the remaining
collision neighborhood and equal to one near $\pi$. Then
\[
 Q_{a,e}[\chi,0]=\|\chi'\|_2^2-a\int_0^\pi |\chi|^2/q_e\,dt
        \longrightarrow-\infty
\]
as $e\uparrow1$. Because these two vectors have disjoint supports,
their span is negative definite, which forces index two and zero
nullity. Equivalently, let $W$ be the Cauchy data of the homogeneous
solutions, $\Lambda$ and $\widehat\Lambda$ the physical and restricted
conditions, and $\Lambda_F$ the full Dirichlet (Friedrichs) condition in
the boundary space of the common regular minimal operator at $e<1$; the
triple-index formula \cite[Eq.~(6.2)]{HPW25} then writes this correction as
$\iota(W,\Lambda,\Lambda_F)-\iota(W,\widehat\Lambda,\Lambda_F)$, which
equals one for $\beta=N$ and zero for $\beta=D$. Consequently, for each
fixed $0<a\le c$ and all $e$ sufficiently close to one,
\begin{equation}\label{eq:euler-critical-sector-count}
 n_-(A^{\alpha\beta}_{a,e})=1+\mathbf{1}_{\{\beta=N\}},\qquad
 \ker A^{\alpha\beta}_{a,e}=0,\qquad \alpha,\beta\in\{N,D\}.
\end{equation}
By the orthogonal decomposition of Appendix~\ref{app:reversible}:
\begin{itemize}
\item \emph{Periodic ($\omega=1$)}: the sectors $NN\oplus DD$ give Morse
index $2+1=3$ with trivial kernel;
\item \emph{Antiperiodic ($\omega=-1$)}: the sectors $ND\oplus DN$ give
Morse index $1+2=3$ with trivial kernel.
\end{itemize}

\section{Graph coordinates and the regular boundary comparison}\label{app:domain-graph-example}

For the Robin family of Example~\ref{ex:robin-boundary-family}, take
$H=L^2(0,1;\mathbb R)$ and let
$\tau u=(u(0),u'(0),u(1),u'(1))$ on $H^2(0,1)$. The minimal operator
$Au=-u''$ has domain $\ker\tau$. For $\alpha\in\mathbb R$ consider
\begin{equation}\label{eq:app-robin-domain}
 B_\alpha u=-u'',\qquad
 \D(B_\alpha)=\{u\in H^2(0,1):u'(0)=\alpha u(0),\ u(1)=0\}.
\end{equation}
Its free traces are $c=u(0)$ and $d=u'(1)$. Letting
$z_1,\ldots,z_4\in\mathcal P_3$ be the cubic Hermite basis with unit
traces $\tau z_j=e_j\in\mathbb R^4$, the cubics
\[
 h_\alpha=z_1+\alpha z_2=1-3t^2+2t^3+\alpha(t-2t^2+t^3),
 \qquad k=z_4=t^3-t^2
\]
have traces $(1,\alpha,0,0)$ and $(0,0,0,1)$. Subtracting
$c h_\alpha+d k$ from $u$ therefore annihilates all four traces.
The values of $c,d$ are unique, giving
\begin{equation}\label{eq:app-robin-domain-sum}
 \D(B_\alpha)=\D(A)\oplus\Span\{h_\alpha,k\}.
\end{equation}
Green's identity also verifies self-adjointness directly. For
$v\in\D(B_\alpha^*)\subset\D(A^*)=H^2(0,1)$ the boundary term is
\[
 \langle -u'',v\rangle-\langle u,-v''\rangle
 =-d\,v(1)-c\,(v'(0)-\alpha v(0)).
\]
The two cubics let $c,d$ vary independently. Thus the boundary term
vanishes for every $u\in\D(B_\alpha)$ exactly when
$v(1)=0$ and $v'(0)=\alpha v(0)$, which gives
$\D(B_\alpha^*)=\D(B_\alpha)$.

Applying $u\mapsto(u,-u'')$ to the domain decomposition gives
\begin{equation}\label{eq:app-robin-graph-sum}
 \gr B_\alpha=\gr A\oplus
 \Span\{(h_\alpha,-h_\alpha''),(k,-k'')\}.
\end{equation}
In the cubic graph complement $E$ of Example~\ref{ex:robin-boundary-family},
tracing the first component identifies the last span $N_\alpha$ with
$\Lambda_\alpha=\{(c,\alpha c,0,d):c,d\in\mathbb R\}$.
Green's identity makes this two-dimensional plane isotropic in the
four-dimensional trace space, hence Lagrangian. Thus
$\D(B_\alpha)=\tau^{-1}\Lambda_\alpha$ and
\begin{equation}\label{eq:app-robin-canonical-sum}
 \gr B_\alpha=\gr A\oplus N_\alpha
\end{equation}
describe the same realization.

\subsection{The regular boundary comparison}\label{app:regular-boundary-inertia}
To verify \eqref{eq:boundary-count-bound-exact}, let $\mathfrak t_\Lambda$
be the quadratic form of the regular realization, set $G=H^1_0$, and let
$P$ be the form-orthogonal complement of $G$ within the form domain.
Testing against compactly supported functions shows that every element of
$P$ is a weak homogeneous solution; because the potential is bounded and
measurable, this implies $H^2$ regularity, so ODE uniqueness and the
regular form-domain boundary conditions identify $P$ via the trace map
with $W_s\cap(\Lambda+\Lambda_D)$. Writing $\gamma u=\beta_u+d_u$ with
$\beta_u\in\Lambda$ and $d_u\in\Lambda_D$, Green's identity gives
$\mathfrak t_\Lambda(u,v)=-\omega(\beta_u,d_v)$ on $P$.
Since the regular form domain embeds compactly in $L^2$, the operator
representing $\mathfrak t_\Lambda$ in a shifted form norm is Fredholm, and
the finite-codimension inertia identity \cite[Lemma~D.10]{HPW25} expresses
$\nm(H_{s,\Lambda})-\nm(D_s)$ as the negative index of this form plus
$\dim(W_s\cap\Lambda_D)-\dim(W_s\cap\Lambda\cap\Lambda_D)$.
On $W_s\cap(\Lambda+\Lambda_D)$ the triple-index form is
$Q(h)=\omega(\beta_h,d_h)$, whose positive index together with the
intersection terms equals $\itriple(W_s,\Lambda,\Lambda_D)$ by
\cite[Appendix~A, Remark~A.7, and Theorem~6.1]{HPW25}. Subtracting the
resulting identities for $\Lambda$ and $\widehat\Lambda$ yields
\eqref{eq:boundary-count-bound-exact}, with the complex case following
identically via the sesquilinear Green form.

\section{Two applications of the boundary coordinates}\label{app:boundary-applications}
The scalar Hill family illustrates passage from a quadratic cutoff
to the eccentricity parameter. The transmission example allows two
ends to become singular independently with coupled outer conditions.

\subsection{A singular Hill-type family}\label{sec:perturbations}

Consider
\begin{equation}\label{eq:hill}
 T_{e,l}^{\min}=-\frac{d^2}{dt^2}+\frac{l}{1-e\cos t},
 \qquad 0<t<1,\quad 0\le e\le1.
\end{equation}
At $e=1$ the left end becomes singular and the right end stays regular.

\begin{theorem}\label{thm:hill}
The minimal realizations in \eqref{eq:hill} form a jointly
gap-continuous family of symmetric Fredholm operators of index $-2$ for
\[
 (e,l)\in[0,1]\times(-\infty,3/8).
\]
Their maximal graph domains embed compactly in $L^2(0,1)$.
\end{theorem}
\begin{proof}
For $e<1$ this is ordinary norm-continuity of bounded potentials on
compact parameter neighborhoods. Near $e=1$, write
\[
 d_e=\frac{1-e}{e},\qquad s_e=\sqrt{2d_e},\qquad q_{e,l}=\frac{2l}{e}.
\]
Then the potential decomposes exactly as
\begin{equation}\label{eq:hill-decomposition}
 \frac{l}{1-e\cos t}
 =\frac{q_{e,l}}{t^2+s_e^2}+B_{e,l}(t),
\end{equation}
where
\[
 B_{e,l}(t)=\frac le
 \frac{t^2/2-(1-\cos t)}{(d_e+1-\cos t)(d_e+t^2/2)}.
\]
The bounds $1-\cos t\ge c t^2$ and
$|t^2/2-(1-\cos t)|\le Ct^4$ give a uniform bound for $B_{e,l}$
near each $(1,l_0)$. Pointwise convergence for $t>0$ gives strong
convergence of the multiplication operators. At $e=1$ the remainder
extends continuously to $t=0$ with value $l_0/6$.

Since $l_0<3/8$, a sufficiently small parameter neighborhood satisfies
$q_{e,l}<3/4$. Setting $r_{e,l}=q^{-1}(q_{e,l})<1$ via the inverse of \eqref{eq:q},
Theorem~\ref{thm:bessel2} ensures that the unperturbed minimal operator
$G_{s_e,r_{e,l}}=-d^2/dt^2+q_{e,l}/(t^2+s_e^2)$ is gap-continuous,
Fredholm of index $-2$, and has compact graph-domain embedding.
Lemma~\ref{lem:strong-compact} handles the remainder. Bounded perturbations leave
minimal and maximal domains unchanged at each fixed parameter and give equivalent
graph norms, which proves the compactness assertion.
\end{proof}

At $e=1$, the Hardy threshold is $2l=-1/4$, or $l=-1/8$.
For $l\ge-1/8$, the Friedrichs realization of the limiting
expression in \eqref{eq:hill}, with Dirichlet condition at $t=1$,
is bounded below and has compact resolvent. For $l<-1/8$, the
quadratic form is unbounded below already on $C_c^\infty(0,1)$.
The graph-continuity range of Theorem~\ref{thm:hill} includes both
regimes and ends at the limit-point threshold $l=3/8$
(Proposition~\ref{prop:threshold}).

\subsection{A coupled family with two singular limits}\label{sec:coupled}
We now allow both ends to become singular independently. Set
\[
H=L^2(0,2)=L^2(0,1)\oplus L^2(1,2),
\qquad
v=(s_-,r_-,s_+,r_+)\in
\bigl([0,1]\times(-\infty,1)\bigr)^2,
\]
and define the piecewise potential
\begin{equation}\label{eq:two-end-potential}
V_v(t)=
\begin{cases}
\displaystyle\frac{q(r_-)}{t^2+s_-^2},&0<t<1,\\[6pt]
\displaystyle\frac{q(r_+)}{(2-t)^2+s_+^2},&1<t<2,
\end{cases}
\end{equation}
where $q(r)$ is defined in \eqref{eq:q}. The potential may jump at
$t=1$, where the two subintervals are joined by the transmission
conditions $y(1^-)=y(1^+)$ and $y'(1^-)=y'(1^+)$.
The cut minimal operator
\[
T_v^\oplus=T_{v,-}^{\min}\oplus T_{v,+}^{\min}
\]
is the direct sum of the two one-sided Bessel regularizations.
By Theorem~\ref{thm:bessel2} and reflection $t\mapsto2-t$,
$T_v^\oplus$ is a jointly gap-continuous family of symmetric,
injective Fredholm operators of index $-4$, whose maximal graph
domains embed compactly into $H$.

\paragraph{Endpoint traces.}
On the left interval $(0,1)$, let $(f_-,g_-)=(u_{s_-,r_-},w_{s_-,r_-})$ be the normalized fundamental
solution frame from \eqref{eq:normalized-frame}. On the right interval $(1,2)$, the reflected frame
\[
 f_+(t)=u_{s_+,r_+}(2-t),\qquad
 g_+(t)=-w_{s_+,r_+}(2-t)
\]
satisfies $W(f_+,g_+) = f_+g_+' - f_+'g_+ = 1$ in the increasing $t$ coordinate. Both frames vary
continuously in $L^2$ and locally in $C^1$ up to the regular cut $c=1$.
For functions in the maximal domain, the outer singular boundary traces and inner regular Cauchy traces are
\begin{equation}\label{eq:outer-traces}
 \begin{split}
 \Gamma_{0,v}y&=([y,g_-](0),[y,f_-](0)),\\
 \Gamma_{2,v}y&=([y,g_+](2),[y,f_+](2)),\\
 z_y(1^\pm)&=(y'(1^\pm),y(1^\pm)).
 \end{split}
\end{equation}
The endpoint limits exist by Green's identity: for any homogeneous frame vector $h$,
\begin{equation}\label{eq:singular-trace-continuity}
 \begin{aligned}
 [y,h](0)&=[y,h](1^-)+\langle (T_v^\oplus)^*y,h\rangle_{(0,1)},\\
 [y,h](2)&=[y,h](1^+)-\langle (T_v^\oplus)^*y,h\rangle_{(1,2)}.
 \end{aligned}
\end{equation}
With the standard symplectic form $\Omega((p,x),(q,y))=py-xq$ on $\R^2$, the Wronskian identity
$W(f_\pm, g_\pm) = 1$ ensures that the endpoint brackets push forward to $\Omega$ on the outer coordinates.
The full trace map
\begin{equation}\label{eq:full-trace}
 \mathcal T_vy=(\Gamma_{0,v}y,z_y(1^-),z_y(1^+),\Gamma_{2,v}y)
\end{equation}
is the concrete trace map $\Tr_v$ constructed in Section~\ref{sec:sl} with $m_\pm = 2$ and $n=1$,
identifying the boundary quotient of $T_v^\oplus$ with $(\R^8, \Omega_v)$, where
\begin{equation}\label{eq:split-form}
 \Omega_v = -\Omega\oplus\Omega\oplus(-\Omega)\oplus\Omega.
\end{equation}

\begin{theorem}\label{thm:coupled}
Let $v=v(\theta)$, $0\le\theta\le1$, be any continuous parameter path for \eqref{eq:two-end-potential}, and let
\[
 \Lambda_\theta\in\Lag(\R^2\oplus\R^2,-\Omega\oplus\Omega)
\]
be a continuous path of outer boundary conditions. Define the self-adjoint realization $B_\theta$ by
the differential expression $-y''+V_{v(\theta)}y$ and the conditions
\begin{equation}\label{eq:coupled-domain}
 \begin{split}
 &(\Gamma_{0,v(\theta)}y,\Gamma_{2,v(\theta)}y)\in\Lambda_\theta,\\
 &y(1^-)=y(1^+),\qquad y'(1^-)=y'(1^+).
 \end{split}
\end{equation}
Then $B_\theta$ is a norm-resolvent-continuous family of self-adjoint Fredholm operators with compact resolvent.
In the coordinates \eqref{eq:full-trace}, let
\[
 \mathcal L_\theta
 =\bigl\{(\Gamma_0,z^-,z^+,\Gamma_2)\in\R^8:
 (\Gamma_0,\Gamma_2)\in\Lambda_\theta,\ z^-=z^+\bigr\}
\]
be the Lagrangian subspace defined by \eqref{eq:coupled-domain}, and let $\mathcal C_\theta = \mathcal T_{v(\theta)}\ker(T_{v(\theta)}^\oplus)^*$.
Then
\begin{equation}\label{eq:coupled-index}
 \sf(B_\theta)=-\Mas(\mathcal L_\theta,\mathcal C_\theta).
\end{equation}
The conclusions hold jointly as either or both $s_\pm \to 0$ and as $r_\pm$ pass through zero.
\end{theorem}
\begin{proof}
At the regular cut, continuity of $(y',y)$ corresponds to the diagonal
Lagrangian $\Delta_1$ for $\Omega\oplus(-\Omega)$, which combined with
the prescribed outer condition gives the continuous Lagrangian family
$\mathcal L_\theta$ in $(\R^8,\Omega_v)$. Since the cut minimal family
$T_{v(\theta)}^\oplus$ is already gap-continuous, injective, and
Fredholm, the continuous frames identify its boundary bundle
symplectically with $(\R^8,\Omega_v)$ via the trace map
\eqref{eq:full-trace}. Pulling back $\mathcal L_\theta$ yields a
continuous family of Lagrangian subspaces of the graph complement, so
Theorem~\ref{thm:canonical}(ii) gives gap continuity (equivalently,
norm-resolvent continuity) of the realizations $B_\theta$. Applying
Theorem~\ref{thm:sf} and expressing the bundle Maslov index in the
symplectic trace coordinates \eqref{eq:full-trace}--\eqref{eq:split-form}
then yields \eqref{eq:coupled-index}. Finally, compact embedding of the
cut maximal graph implies that every realization has compact resolvent.
\end{proof}

\begin{remark}\label{rem:coupled-morse}
If both singular coefficients are at least $-1/4$, each realization
has finite Morse index, and a uniform lower bound along the path
identifies spectral flow with the endpoint index difference. Without
such a bound eigenvalues may escape to $-\infty$; in particular, a
supercritical singular end yields infinite Morse index, while
\eqref{eq:coupled-index} still counts crossings through zero.
\end{remark}

\section*{Declaration of generative AI and AI-assisted technologies in the manuscript preparation process}
During the preparation of this work, the authors used generative AI tools
to assist with language editing and selected special-function calculations.
After using these tools, the authors reviewed and edited the content as
needed and take full responsibility for the content of the publication.

\par\medskip\noindent
\begingroup\small
\begin{minipage}[t]{0.31\textwidth}
\textsc{Xijun Hu}\\
School of Mathematics\\
Shandong University\\
Jinan 250100, P. R. China\\
\href{mailto:xjhu@sdu.edu.cn}{\texttt{xjhu@sdu.edu.cn}}
\end{minipage}\hfill
\begin{minipage}[t]{0.39\textwidth}
\textsc{Alessandro Portaluri}\\
DISAFA, Universit\`a di Torino\\
Largo Braccini 2, Grugliasco, Italy\\
\href{mailto:alessandro.portaluri@unito.it}{\texttt{alessandro.portaluri@unito.it}}\\
New York University Abu Dhabi\\
\href{mailto:ap9453@nyu.edu}{\texttt{ap9453@nyu.edu}}
\end{minipage}\hfill
\begin{minipage}[t]{0.28\textwidth}
\textsc{Li Wu}\\
School of Mathematics\\
Shandong University\\
Jinan 250100, P. R. China\\
\href{mailto:vvvli@sdu.edu.cn}{\texttt{vvvli@sdu.edu.cn}}
\end{minipage}
\endgroup
\end{document}